\documentclass{amsart}

\usepackage{amsmath,amsxtra,amssymb,amsthm,amsfonts,color}

\allowdisplaybreaks

\newtheorem{theorem}{Theorem}[section]
\newtheorem{lemma}[theorem]{Lemma}
\newtheorem{proposition}[theorem]{Proposition}
\newtheorem{corollary}[theorem]{Corollary}
\newtheorem{remark}[theorem]{Remark}

\newtheorem{theoA}{Theorem}
\newtheorem{theoB}{Theorem}

\newtheorem{coroAA}{Corollary}
\newtheorem{coroAB}{Corollary}
\newtheorem{coroB}{Corollary}

\theoremstyle{definition}

\newcommand{\N}{\mathbb{N}}

\newcommand{\R}{\mathbb{R}}
\newcommand{\C}{\mathbb{C}}

\newcommand{\E}{\mathcal{E}}

\newcommand{\U}{\mathcal{U}}
\newcommand{\oM}{\otimes_{\mathcal{M}}}
\newcommand{\la}{\mathbb{\lambda}}
\newcommand{\al}{\mathbb{\alpha}}

\newcommand{\ten}{\otimes}

\newcommand{\eps}{\varepsilon}

\newcommand{\dem}{\noindent {\bf Proof. }}

\newcommand{\fin}{\hspace*{\fill} $\square$ \vskip0.2cm}
\newcommand{\prodd}{\prod\nolimits}
\newcommand{\limm}{\lim\nolimits}

\newcommand{\summ}{\sum\nolimits}

\addtocounter{tocdepth}{-1}

\begin{document}

\null

\vskip-1.5cm

\title[A transference method]{A transference method in
quantum probability}

\author[Junge and Parcet]
{Marius Junge$^{\ast}$ and Javier Parcet$^{\dag}$}

\footnote{$^{\ast}$Partially supported by the NSF DMS-0556120.}

\footnote{$^{\dag}$Partially supported by \lq Programa Ram{\'o}n y
Cajal, 2005\rq${}$ and Grant MTM2007-60952, Spain.}


\maketitle


\tableofcontents

\section*{Introduction}

\renewcommand{\theequation}{$\tau_{p}$}
\addtocounter{equation}{-1}

The notion of independent random variables is central in
probability theory and has many applications in analysis.
Independence is also a fundamental concept in quantum probability,
where it can occur in many different forms. In terms of norm
estimates for sums of independent variables, free probability
often plays the role of the best of all worlds. This is
particularly true for applications in the theory of operator
spaces. We refer to the so-called Grothendieck's program for
operator spaces \cite{HM,PS,X2} and also to the noncommutative
$L_p$ embedding theory \cite{J2,JP3,JP4} due to the authors. On
the other hand, other notions of independence weaker than freeness
are often enough in the context of noncommutative Khintchine or
Rosenthal type inequalities \cite{J3,JX4,LuP}. A first motivation
for this paper was to remove a singularity at $p=1$ for the
classical Rosenthal's inequality \cite{Ros} and its noncommutative
form \cite{JX,JX4}, which is also related to the recent work by
Haagerup and Musat \cite{HM2} on a direct proof of Khintchine
inequalities for the generators of the CAR algebras. This easily
follows from our main result in this paper, a general transference
method which allows us to compare the norm of sums of independent
copies with the norm of sums of freely independent copies.

\vskip3pt

Let us illustrate our transference method. If $\mathcal{M}$ is a
von Neumann algebra, let $\mathcal{M}^{\ten_n}$ be the $n$-fold
tensor product and $\pi_{tens}^k: \mathcal{M} \to
\mathcal{M}^{\ten_n}$ the canonical $k$-th coordinate
homomorphism. It is standard to extend $\pi_{tens}^k:
L_p(\mathcal{M}) \to L_p(\mathcal{M}^{\ten_n})$ for $1 \le p \le
\infty$, see e.g. \cite{JX}. Similarly, we have $k$-th coordinate
homomorphisms $\pi_{free}^k: \mathcal{M} \to
(\mathcal{M},\phi)^{*_n}$ for free products. Our first result
implies that
\begin{equation} \label{main}
\frac{1}{cp} \, \Big\| \sum_{k=1}^n \pi_{tens}^k(x) \Big\|_p \le
\Big\| \sum_{k=1}^n \pi_{free}^k(x) \Big\|_p \le cp \, \Big\|
\sum_{k=1}^n \pi_{tens}^k(x) \Big\|_p
\end{equation}
for $x \in L_p(\mathcal{M})$ and $n \ge 1$, with uniformly bounded
constants when $p$ is close to $1$.

\vskip3pt

\renewcommand{\theequation}{\arabic{equation}}

In the following we will write $a \sim_c b$ when $\frac1c \le a/b
\le c$. We also need an operator-valued version of \eqref{main}
for further applications to the theory of operator spaces. This
requires the notion of independence over a given subalgebra. Let
us formally introduce the notion of \lq independent copies\rq${}$
that we will work with. Given a noncommutative probability space
$(\mathcal{A},\varphi)$ equipped with a normal faithful state
$\varphi$, a von Neumann subalgebra of $\mathcal{A}$ is called
conditioned if it is invariant under the action of the modular
group $\sigma_t^{\varphi}$. By Takesaki \cite{Ta2}, this holds if
and only if there is a $\varphi$-invariant normal faithful
conditional expectation. Let $\mathcal{N}$ be conditioned in
$\mathcal{A}$ with faithful conditional expectation
$\mathsf{E}_\mathcal{N}: \mathcal{A} \to \mathcal{N}$. Let us
consider two von Neumann subalgebras $\mathcal{M}_1,
\mathcal{M}_2$ of $\mathcal{A}$ satisfying $\mathcal{N} \subset
\mathcal{M}_1 \cap \mathcal{M}_2$. Then, $\mathcal{M}_1$ and
$\mathcal{M}_2$ are called \emph{independent over} $\mathcal{N}$
if
$$\mathsf{E}_\mathcal{N}(a_1 a_2) = \mathsf{E}_\mathcal{N}(a_1)
\mathsf{E}_\mathcal{N}(a_2)$$ holds for all $a_1 \in
\mathcal{M}_1$ and $a_2 \in \mathcal{M}_2$. Now, if
$(\mathcal{M}_k)_{k \ge 1}$ are conditioned subalgebras of
$\mathcal{A}$ with $\mathcal{N} \subset \mathcal{M}_k \subset
\mathcal{A}$, we shall say that the system $(\mathcal{M}_k)_{k \ge
1}$ is \emph{increasingly independent} if
\begin{itemize}
\item[\textbf{a)}] $\big\langle \mathcal{M}_1, \mathcal{M}_2,
\ldots, \mathcal{M}_{k-1} \big\rangle$ and $\mathcal{M}_k$ are
independent over $\mathcal{N}$.
\end{itemize}
We also need a technical notion. Given a von Neumann algebra
$\mathcal{M}$ containing $\mathcal{N}$ and a family of
$*$-isomorphisms $\pi_k: \mathcal{M} \to \mathcal{M}_k$ with
$\mathcal{N} \subset \mathcal{M}_k \subset \mathcal{A}$, we say
that $(\mathcal{M}_k)_{k \ge 1}$ is a system of
\emph{top-subsymmetric copies of $\mathcal{M}$ over} $\mathcal{N}$
if
\begin{itemize}
\item[\textbf{b)}] ${\pi_k}_{|_\mathcal{N}} = id$ and
$$\mathsf{E}_\mathcal{N} \big( \pi_{f(1)}(x_1) \cdots \pi_{f(m)}(x_m)
\big) = \mathsf{E}_\mathcal{N} \big( \pi_{g(1)}(x_1) \cdots
\pi_{g(m)}(x_m) \big)$$ holds for all $f,g: \{1,2,\ldots,m\} \to
\N$ satisfying
\begin{itemize}
\item[$\bullet$] $f_{|_{\{1,2,\ldots,m\} \setminus A}} =
g_{|_{\{1,2,\ldots,m\} \setminus A}}$,

\item[$\bullet$] $|A| \le 2$ and $A = \big\{ k \ | \, f(k) = \max
f \big\} = \big\{ k \ | \, g(k) = \max g \big\}$.
\end{itemize}
\end{itemize}

Of course, when no subalgebra $\mathcal{N}$ is specified, we shall
work with $\mathcal{N} = \langle \mathbf{1}_{\mathcal{A}} \rangle$
and $\mathsf{E}_\mathcal{N} = \varphi$. Using the assumptions of
conditioned subalgebras allows us to provide $L_p$ generalizations
of the conditional expectation $\mathsf{E}_\mathcal{N}$ and the
isomorphisms $\pi_k$, see \cite{JX} for details. Intuitively
speaking, condition b) means that we are allowed to exchange the
top element in the range of $f$ by the top element in the range of
$g$, but only if the top element does not occur more than once or
twice. Top-subsymmetry is exactly the technical assumption which
makes the argument in \cite{J3} work. Nevertheless, as in
\cite{J3}, we can consider two alternative stronger conditions:

\vskip5pt

\begin{itemize}
\item[\textbf{b2)}] \emph{Subsymmetry}: $g(k) = \varphi \circ
f(k)$ for any strictly increasing $\varphi: \N \to \N$.

\vskip5pt

\item[\textbf{b3)}] \emph{Symmetry}: $g(k) = \sigma \circ f(k)$
for any permutation $\sigma$ of the positive integers.
\end{itemize}

\vskip5pt

\noindent It is clear that the implications below hold
$$\mbox{Symmetry} \Rightarrow \mbox{Subsymmetry} \Rightarrow
\mbox{Top-subsymmetry}.$$

\noindent \textbf{Example 1.} Tensor product copies. Let
$$\Big( \mathcal{A}_{tens}, \mathcal{M}, \mathcal{M}_{tens}^k,
\mathcal{N}; \mathsf{E}_\mathcal{N} \Big) = \Big( \mathcal{N}
\bar\otimes \mathcal{R}^{\otimes_n}, \mathcal{N} \bar\ten
\mathcal{R}, \pi_{tens}^k(\mathcal{M}), \mathcal{N}; id \ten
\varphi_\mathcal{R}^{\ten_n} \Big)$$ where $n \in \N \cup
\{\infty\}$ and the homomorphisms $\pi_{tens}^k$ are given by
$$\pi_{tens}^k: n \ten x \in \mathcal{M} \mapsto n \ten
\mathbf{1}_{\mathcal{R}} \ten \ldots \ten \mathbf{1}_{\mathcal{R}}
\ten \underbrace{x}_{k-\mathrm{th}} \ten \mathbf{1}_{\mathcal{R}}
\ten \ldots \ten \mathbf{1}_{\mathcal{R}} \in
\mathcal{M}_{tens}^k.$$ The $\mathcal{M}_{tens}^k$'s form an
independent symmetric system of copies of $\mathcal{M}$ over
$\mathcal{N}$.

\vskip5pt

\noindent \textbf{Example 2.} Freely independent copies. Consider
$$\Big( \mathcal{A}_{free}, \mathcal{M}, \mathcal{M}_{free}^k,
\mathcal{N}; \mathsf{E}_\mathcal{N} \Big) = \Big( *_{\mathcal{N}
\! , \, k} (\mathcal{M}, \mathcal{E}_\mathcal{N}), \mathcal{M},
\pi_{free}^k(\mathcal{M}), \mathcal{N}; \mathsf{E}_{\mathcal{N}}
\Big),$$ the reduced $\mathcal{N}$-amalgamated free product of
$(\mathcal{M}, \mathcal{E}_\mathcal{N})$ with
$\mathcal{E}_\mathcal{N}: \mathcal{M} \to \mathcal{N}$ a normal
faithful conditional expectation, see e.g. \cite{JPX} for details
on the construction of reduced amalgamated free product von
Neumann algebras. The isomorphism $\pi_{free}^k$ is the canonical
embedding into the $k$-th component of the free product
$\mathcal{A}_{free}$. The $\mathcal{M}_{free}^k$'s form an
independent symmetric system of copies of $\mathcal{M}$ over
$\mathcal{N}$.

\vskip3pt

Our notion of noncommutative independent copies is quite general.
We refer to \cite{JX4} for more examples which arise naturally in
quantum probability. The first form of our transference principle
is the following.

\begin{theoA}
Let $1 \le p \le 2$ and let $(\mathcal{M}_k)_{k \ge 1}$ be an
increasingly independent family of top-subsymmetric copies of
$\mathcal{M}$ over $\mathcal{N}$. Then, there exists a positive
constant $c$ independent of $p$ and $n$ such that
$$\mathbb{E} \Big\| \sum_{k=1}^n \eps_k \hskip1pt \pi_k(x)
\Big\|_{L_p(\mathcal{A})} \sim_c \mathbb{E} \Big\| \sum_{k=1}^n
\eps_k \hskip1pt \pi_{free}^k(x)
\Big\|_{L_p(\mathcal{A}_{free})}.$$ If in addition the
$\mathcal{M}_k$'s are symmetric, then $$\Big\| \sum_{k=1}^n
\pi_k(x) \Big\|_{L_p(\mathcal{A})} \sim_c \Big\| \sum_{k=1}^n
\pi_{free}^k(x) \Big\|_{L_p(\mathcal{A}_{free})}.$$
\end{theoA}

\vskip3pt

\noindent The second form of transference is stated below.

\begin{theoB}
Let $1 \le p \le q \le \infty$ and let $(\mathcal{M}_k)_{k \ge 1}$
be an increasingly independent family of top-subsymmetric copies
of $\mathcal{M}$ over $\mathcal{N}$. Then, there exists a positive
constant $c$ independent of $p,q$ and $n$ such that
$$\Big\| \sum_{k=1}^n \pi_k(x) \ten \delta_k
\Big\|_{L_p(\mathcal{A}; \ell_q^n)} \sim_c \Big\| \sum_{k=1}^n
\pi_{free}^k(x) \ten \delta_k \Big\|_{L_p(\mathcal{A}_{free};
\ell_q^n)}.$$
\end{theoB}

We know from \cite{J3} that Theorem A holds for $p=1$. The strategy consists of applying our technique \cite{JP4,JP2} to provide a complete embedding $L_p \to L_1$ which preserves independence. This is done in Section \ref{Section1} and the rest of the paper will be essentially devoted to the proof of Theorem B, which is similar in nature but requires to adapt all the methods in \cite{J3,JP4,JP2}. As for Theorem A, our strategy is to prove the result in the extremal case $(p,q) = (1,\infty)$ and show that the general statement reduces to it. The extremal case is a consequence of Theorem \ref{4-term}, where we characterize the norm in $L_1(\mathcal{A}; \ell_\infty^n(\mathcal{R}))$ of increasingly independent top subsymmetric copies for any finite dimensional von Neumann algebra $\mathcal{R}$. The reduction argument is divided in two cb-embeddings $$L_p(\mathcal{A}; \ell_q^n) \to L_p(\widehat{\mathcal{A}}; \ell_\infty^{mn}) \to \prodd_{s,\mathcal{U}} L_1 \big(M_s(\widehat{\mathcal{A}})^{\otimes_{\mathrm{k}_s}}; \ell_\infty^{\mathrm{k}_smn} \big),$$ both preserving independence. Note that this map takes values in an ultraproduct of spaces of the form $L_1(\mathcal{A}'; \ell_\infty^n(\mathcal{R}))$, so that we are in position to apply Theorem \ref{4-term}. The second embedding might be of independent interest and will be proved in Theorem \ref{transfor}, while the first embedding is the content of Theorem \ref{vqq}. In both Theorems A and B, the main new difficulty relies on keeping track of independence in the construction of the embedding.

The drawback is that Theorems A and
B only hold for independent copies. This restriction goes back to
\cite{J3}. Any progress in the non identically distributed case
would be very desirable and thus we propose the following problem:

\vskip5pt

\noindent \textbf{Problem 1.} Do the scalar and mixed-norm
transference hold for non i.d. variables?

\renewcommand{\theequation}{$\Sigma_{p}$}
\addtocounter{equation}{-1}

\vskip5pt

Now we may revisit the singularity of certain constants mentioned
above. Given $2 \le p < \infty$, a probability space
$(\Omega,\mu)$ and $f_1, f_2, \ldots \in L_p(\Omega)$ a family of
mean-zero independent random variables, Rosenthal's classical
inequality gives
\begin{equation} \label{Ros>2}
{} \hskip24pt \Big( \int_\Omega \big| \sum_{k=1}^n f_k \big|^p
d\mu \Big)^{\frac1p} \sim_{c_p} \hskip5pt \max \left\{ \Big(
\sum_{k=1}^n \|f_k\|_p^p \Big)^{\frac1p}, \Big( \sum_{k=1}^n
\|f_k\|_2^2 \Big)^{\frac12} \right\}.
\end{equation}
As a byproduct, we obtain for $1 \le q \le p < \infty$
\renewcommand{\theequation}{$\Sigma_{pq}$}
\addtocounter{equation}{-1}
\begin{equation} \label{Mix>2}
\Big( \int_\Omega \big( \sum_{k=1}^n |f_k|^q \big)^{\frac{p}{q}}
d\mu \Big)^{\frac1p} \sim_{c_{p,q}} \max \left\{ \Big(
\sum_{k=1}^n \|f_k\|_p^p \Big)^{\frac1p}, \Big( \sum_{k=1}^n
\|f_k\|_q^q \Big)^{\frac1q} \right\}
\end{equation}
with the $f_k$'s not necessarily mean-zero. Indeed, the case $q=2$
easily follows from Khintchine and Rosenthal inequalities, while
the general case follows from an immediate renormalization
argument. Notice that in both cases we end up with the norm of an
intersection of Banach spaces, whose dual is the sum of the
corresponding dual spaces. This simple observation produces dual
inequalities for $1 < p \le 2$ and $1 < p \le q \le \infty$ as
follows
\begin{eqnarray*}
\Big( \int_\Omega \big| \sum_{k=1}^n f_k \big|^p d\mu
\Big)^{\frac1p} & \sim_{c_p} & \inf_{f_k = \phi_k + \psi_k}
\left\{ \Big( \sum_{k=1}^n \|\phi_k\|_p^p \Big)^{\frac1p}, \Big(
\sum_{k=1}^n \|\psi_k\|_2^2 \Big)^{\frac12} \right\}, \\ \Big(
\int_\Omega \big( \sum_{k=1}^n |f_k|^q \big)^{\frac{p}{q}} d\mu
\Big)^{\frac1p} & \sim_{c_{p,q}} & \inf_{f_k = \phi_k + \psi_k}
\left\{ \Big( \sum_{k=1}^n \|\phi_k\|_p^p \Big)^{\frac1p}, \Big(
\sum_{k=1}^n \|\psi_k\|_q^q \Big)^{\frac1q} \right\}.
\end{eqnarray*}
It is worth mentioning that the martingale version of Rosenthal
inequality \cite{B} was extended to $1 < p \le 2$ from \cite{JX} and
weak type estimates for $p=1$ were unknown until \cite{RMI}. On the
other hand, we know from \cite{JSZ} that the best constant $c_p$ in
Rosenthal's inequality behaves like $p / \log p$. In particular, we
find for \eqref{Ros>2} and \eqref{Mix>2} a non-removable singularity
at $\infty$ which is carried over to $1$ by our duality argument.
The problem is to decide whether this singularity is removable ---as
it happens for the Khintchine inequality--- or not, either in the
classical or in the quantum setting.

\vskip5pt

A direct argument to remove it not involving duality seems out of reach by now.
Note that precise decompositions $f_k = \phi_k + \psi_k$ have only been studied in
\cite{RMI,RanAnn} for martingales. However, free random variables
are fortunately at our disposal and we know from \cite{JP2,JPX} that
the free forms of \eqref{Ros>2} and \eqref{Mix>2} do not have a
singularity at $\infty$ and duality solves in that case our problem.
The validity of Khintchine and Rosenthal type inequalities in the
extremal case $p=\infty$ is a stamp of free probability, see
\cite{PP,RX} for related results. Our transference method in
Theorems A and B solves our problem for identically distributed
variables. It is also worth mentioning that the argument requires
freeness even in the classical case with commutative $f_k$'s!

\renewcommand{\theequation}{\arabic{equation}}

\vskip5pt

\noindent \textbf{Problem 2.} Are the dual forms of \eqref{Ros>2}
and \eqref{Mix>2} singular for non i.d. variables?

\vskip5pt

Let us comment some further applications of transference. We
recall Hiai's construction \cite{Hi} of the $q$-deformed analogue
of Shlyakhtenko's generalized circular variables. Consider a
complex Hilbert space $\mathcal{H}$ equipped with a distinguished
unit vector $\Omega$ and denote by $\mathcal{F}_{q}(\mathcal{H})$
the associated $q$-Fock space. If $q=\pm 1$, we find the
well-known Bosonic and Fermionic Fock spaces equipped with the
symmetric and antisymmetric structures. When $-1 < q < 1$ we
follow \cite{BKS} and equip it with the $q$-inner product induced
by $$\big\langle f_1\ten \cdots \ten f_n,\; g_1 \ten \cdots \ten
g_m \big\rangle_q =\delta_{nm} \sum_{\pi \in
\mathcal{S}_n}^{\null} q^{i(\pi)}\langle f_1,\; g_{\pi(1)} \rangle
\cdots \langle f_n,\; g_{\pi(n)}\rangle.$$ Let $\ell_q(e)$ and
$\ell_q^*(e)$ stand for the creation and annihilation operators
associated to a vector $e \in \mathcal{H}$, see \cite{BKS} for
precise definitions. Assume $\mathcal{H}$ is infinite dimensional
and separable, so that we can fix an orthonormal basis $(e_{\pm
k})_{k \ge 1}$. Given two sequences $(\lambda_k)_{k \ge 1}$ and
$(\mu_k)_{k \ge 1}$ of positive numbers, set $$gq_k = \lambda_k
\ell_q(e_k) + \mu_k \ell_q^*(e_{-k}) \quad \mbox{and} \quad
gq_{k,p} = d_{\phi_q}^{\frac{1}{2p}} gq_k
d_{\phi_q}^{\frac{1}{2p}}.$$ The von Neumann algebra generated by
the $gq_k$'s in the GNS-construction with respect to the vacuum
state $\phi_q(\cdot) = \langle \Omega, \cdot \, \Omega \rangle_q$
will be denoted by $\Gamma_q$ and represent the $q$-deformed
analogue of the corresponding Araki-Woods factor in the
antisymmetric case. Here $d_{\phi_q}$ denotes the density of
$\phi_q$.

\begin{coroAA}
Let $\mathcal{M}$ be a von Neumann algebra and $1 \le p \le 2$.
Let us consider a finite sequence $x_1, x_2, \ldots, x_n$ in
$L_p(\mathcal{M})$. Then, the following equivalences hold for any
$-1 \le q \le 1$ up to a constant $c$ independent of $p,q$ and $n$
\begin{eqnarray*}
\lefteqn{\hskip-5pt \Big\| \sum_{k=1}^n x_k \ten gq_{k,p}
\Big\|_{L_p(\mathcal{M} \bar\ten \Gamma_q)}} \\ & \sim_c &
\inf_{x_k = a_k + b_k} \Big\| \Big( \summ_k
\lambda_k^{\frac{2}{p}} \mu_k^{\frac{2}{p'}} a_ka_k^*
\Big)^{\frac12} \Big\|_{L_p(\mathcal{M})} + \Big\| \Big( \summ_k
\lambda_k^{\frac{2}{p'}} \mu_k^{\frac{2}{p}} b_k^*b_k
\Big)^{\frac12} \Big\|_{L_p(\mathcal{M})}.
\end{eqnarray*}
\end{coroAA}

The weighted Khintchine type inequalities considered above were
already proved in \cite{J3,JPX,X3}. The novelty of our result
relies on the nonsingularity of the constants involved. To be more
precise, we explain this point with a series of remarks:

\vskip5pt

\begin{itemize}
\item[\textbf{a)}] In the Fermionic case, Corollary A1 solves a
question by Xu in \cite{X3}. More concretely, Xu proved the
weighted Fermionic Khintchine inequality for $1 < p < \infty$ with
singularities at $1$ and $\infty$. The singularity at $\infty$ is
already predicted by the classical Khintchine inequality. However,
the first-named author proved in \cite{J3} the same inequality for
$p=1$ applying a central limit procedure to a Rosenthal type
inequality for independent copies in $L_1$. This motivated Xu to
ask whether the singularity at $1$ was removable.

\vskip5pt

\item[\textbf{b)}] On the other hand, Haagerup and Musat recently
used in \cite{HM2} another new argument to prove the weighted
Fermionic Khintchine inequality. Their method is not only simpler
than the one in \cite{J3}, but they also managed to improve the
constant up to $\sqrt{2}$. Unfortunately, their concrete approach
in $L_1$ seems not to work for $p > 1$ and our proof uses instead
the scalar-valued version of our transference method in Theorem A,
together with a central limit procedure as in \cite{J3}.

\vskip5pt

\item[\textbf{c)}] In \cite{JPX} Corollary A1 was proved with
$$\qquad c_{p,q} \le \Big( \frac{2}{\sqrt{1-|q|}} \Big)^{|1-\frac2p|}
\qquad \mbox{for} \qquad -1 < q <1.$$ Note that $q$ here has
nothing to do with the $q$ in \eqref{Mix>2}! It was also shown
that the same bound applies for the cb-complementation constant
$\gamma(p,q)$ of the subspace of $L_p(\Gamma_q)$ generated by the
generalized $q$-gaussians. Corollary A1 now provides a uniform
bound for $c_{p,q}$ as far as $1 \le p \le 2$. Moreover, by
Corollary A1 and the argument in \cite{JPX}, we have $$c_{p,q}
\sim \gamma(p,q) \quad \mbox{for} \quad p \ge 2.$$
\end{itemize}

\vskip5pt

\noindent The following question looks like the next step.

\vskip5pt

\noindent \textbf{Problem 3.} Find accurate estimates for
$\gamma(p,q)$ near $(p,q) = (\infty,\pm 1)$.

\vskip5pt

Our second application has to do with some recent results on
operator space $L_p$ embedding theory. Given $1 \le p < q \le 2$ and
a von Neumann algebra $\mathcal{M}$, the main result in
\cite{JP4,JP2} is the construction of a completely isomorphic
embedding of $L_q(\mathcal{M})$ into $L_p(\mathcal{A})$ for some
sufficiently large von Neumann algebra $\mathcal{A}$, where both
spaces are equipped with their natural operator space structure. We
refer to \cite{J2,P4,X3} for some prior results. The simplest
construction uses free probability techniques and does not produce
any singularity in the embedding constants. However, these
techniques are not the right ones to show the stability of
hyperfiniteness. In other words, whenever $\mathcal{M}$ is
hyperfinite we can also take $\mathcal{A}$ to be hyperfinite. Under
these conditions, all the known constructions produce a singularity
in the constant $\eta(p,q)$ of the cb-embedding $L_q \to L_p$ as $p
\to 1$. We solve this by transference.

\begin{coroAB}
Let $\mathcal{M}$ be hyperfinite and $1 \le p \le q \le 2$. Then,
there exists a completely isomorphic embedding of
$L_q(\mathcal{M})$ into $L_p(\mathcal{A})$ where both spaces are
equipped with their natural operator space structures and satisfy:
\begin{itemize}
\item[\textbf{i)}] $\mathcal{A}$ is hyperfinite.

\item[\textbf{ii)}] The constants are independent of $p,q$.
\end{itemize}
\end{coroAB}

So far we have provided applications of Theorem A. Our third
application now follows from Theorem B. Pisier proved in \cite{P5}
that there is no possible cb-embedding of $\mathrm{OH}$ into the
predual of a semifinite von Neumann algebra. This was generalized
by Xu, who proved in \cite{X3} that for $1 \le p < q \le 2$ we can
not cb-embed $C_q$ or $R_q$ into a semifinite $L_p$ space. In
particular, the same applies for $$S_q = C_q \ten_h R_q.$$

\vskip3pt

\noindent The following result completes the previous ones by
Pisier and Xu.

\begin{coroB}
If $1 \le p < q \le 2$, there is no cb-embedding of $\ell_q$ into
semifinite $L_p$.
\end{coroB}

In fact, as it was pointed out by Pisier our techniques go a bit
further, but that will be explained in the body of the paper.
Corollary B, together with the results by Pisier and Xu, justify
the relevance of type III von Neumann algebras in operator space
$L_p$ embedding theory. The main tools in our argument are:

\vskip5pt

\begin{itemize}
\item The noncommutative form of Rosenthal theorem from
\cite{JP3}.

\vskip3pt

\item Xu's nonembedding techniques from \cite{X3} {\`a} la
Grothendieck \cite{HM,PS,X2}.

\vskip3pt

\item A local cb-embedding $$x \in S_q^n \mapsto \frac{1}{n^{1/q}}
\sum_{k=1}^{n^2} \pi_{tens}^k(x) \ten \delta_k \in L_p \big(
M_{n^{n^2}}; \ell_q^{n^2} \big),$$ with constants independent of
$p,q$ and $n$. This result improves a previous one from \cite{JP},
where the most natural case $p=1$ was not proved and the constant
presented a singularity at $1$.
\end{itemize}

\section{A transference method}
\label{Section1}

\numberwithin{equation}{section}

Given a probability space $(\Omega, \mu)$, we denote the
expectation by $\mathbb{E}(f)=\int fd\mu$. We reserve the symbols
$(\eps_k)_{k\ge 1}$ for a sequence of independent Bernoulli random
variables, i.e. Prob$(\eps_k=\pm 1)=\frac12$. We shall write
$\mathcal{E_N}$ for $\mathsf{E}_\mathcal{N} \circ \pi_k$. We begin
by stating the main result in \cite{J3}. The second part is a
simple refinement from \cite{JP2}.

\begin{theorem} \label{Theorem-Araki-Woods} The
following inequalities hold for $x \in L_1(\mathcal{M}) \! :$

\vskip5pt

\begin{itemize}
\item[\textbf{i)}] If $(\mathcal{M}_k)_{k
\ge 1}$ are independent top-subsymmetric over $\mathcal{N}$, then
\begin{eqnarray*}
\quad \lefteqn{\hskip+15pt \mathbb{E} \Big\| \sum_{k=1}^n
\varepsilon_k \pi_k(x) \Big\|_{L_1(\mathcal{A})}} \\ & \sim &
\inf_{x = a+b+c} n \|a\|_{L_1(\mathcal{M})} + \sqrt{n} \big\|
\mathcal{E_N} (bb^*)^{\frac12} \big\|_{L_1(\mathcal{N})} +
\sqrt{n} \big\| \mathcal{E_N}(c^*c)^{\frac12}
\big\|_{L_1(\mathcal{N})}.
\end{eqnarray*}

\vskip5pt

\item[\textbf{ii)}] If moreover, $\mathcal{E_N}(x)=0$ and
$(\mathcal{M}_k)_{k \ge 1}$ are symmetric over $\mathcal{N}$, then
\begin{eqnarray*}
\quad \lefteqn{\qquad \Big\| \sum_{k=1}^n \pi_k(x)
\Big\|_{L_1(\mathcal{A})}} \\ & \sim & \inf_{x = a+b+c} n
\|a\|_{L_1(\mathsf{A})} + \sqrt{n} \big\| \mathcal{E_N}
(bb^*)^{\frac12} \big\|_{L_1(\mathcal{N})} + \sqrt{n} \big\|
\mathcal{E_N}(c^*c)^{\frac12} \big\|_{L_1(\mathcal{N})}.
\end{eqnarray*}
\end{itemize}
In both cases, the relevant constants are independent of $\mathcal{M}$ and $\mathcal{N}$.
\end{theorem}

We now give some elementary remarks on symmetric tensor products.
Given a positive integer $m$, the symmetric tensor product
$\ten_{sym}^m \mathcal{M}$ is defined as the von Neumann algebra
$$\ten_{sym}^m \mathcal{M} = \overline{\Big\{ \summ_k x_k^{\ten_m}
\big| \ x_k \in \mathcal{M} \Big\}}^{\mathrm{wot}} \subset
\mathcal{M}^{\bar\ten_m}.$$ Let us write $\mathcal{S}_m$ to denote
the symmetric group on $\{1,2,\ldots,m\}$. It is easily seen that
the symmetric tensor product $\ten_{sym}^m \mathcal{M}$ is exactly
the fix point algebra of the conditional expectation
$\mathcal{E}_{sym}(x_1 \ten \cdots \ten x_m) = \frac{1}{m!}
\sum_{\sigma \in \mathcal{S}_m} x_{\sigma(1)} \ten \cdots \ten
x_{\sigma(m)}.$ 

\begin{lemma} \label{Exchange}
Let $(\mathcal{M}_k)_{k \ge 1}$ be independent top-subsymmetric
copies over $\mathcal{N} \! :$
\begin{itemize}
\item[\textbf{i)}] The $\ten_{sym}^m \mathcal{M}_k$'s are
independent top-subsymmetric copies over $\ten_{sym}^m
\mathcal{N}$.

\item[\textbf{ii)}] Consider the $j$-th coordinate homomorphisms
$\pi_{tens}^j(x) = \mathbf{1} \ten \ldots \ten x \ten \ldots \ten
\mathbf{1}$ into the corresponding $m$-fold tensor product. Then,
given a von Neumann algebra $\mathcal{R}$ and $x \in
L_1(\mathcal{M} \bar\otimes \mathcal{R})$, the following estimate
holds with constants independent of $m,n$
$$\hskip38pt \mathbb{E} \Big\| \sum_{j=1}^m \pi_{tens}^j
\big( \sum_{i=1}^n \eps_i \big[ \pi_i \ten id \big] (x) \big)
\Big\|_1 \sim_c \mathbb{E} \Big\| \sum_{j=1}^m \pi_{tens}^j \big(
\sum_{i=1}^n \eps_i \big[ \pi^i_{free} \ten id \big] (x) \big)
\Big\|_1.$$

\item[\textbf{iii)}] Moreover, if $\mathcal{E_N} \ten id (x) = 0$
and the $\mathcal{M}_k$'s are symmetric $$\qquad \qquad \Big\|
\sum_{j=1}^m \pi_{tens}^j \big( \sum_{i=1}^n \big[ \pi_i \ten id
\big] (x) \big) \Big\|_1 \sim_c \Big\| \sum_{j=1}^m \pi_{tens}^j
\big( \sum_{i=1}^n \big[ \pi^i_{free} \ten id \big] (x) \big)
\Big\|_1.$$
\end{itemize}
\end{lemma}

\dem The first assertion trivially follows from the fact that
$$\E_{\ten^m_{sym}\mathcal{N}}(x \ten x \ten \cdots \ten
x) = \mathcal{E_N}(x) \ten \mathcal{E_N}(x) \ten \cdots \ten
\mathcal{E_N}(x)\, .$$ To prove ii), we set
\begin{eqnarray*}
\mathbb{M} & = & \ten_{sym}^m \big( \mathcal{M} \ten \mathcal{R}
\big), \\ \mathbb{M}_k & = & \ten_{sym}^m \big( \mathcal{M}_k \ten
\mathcal{R} \big).
\end{eqnarray*}
Let $\widehat{\pi}_k: \mathbb{M} \to \mathbb{M}_k$ given by
$\widehat{\pi}_k(x^{\ten_m})= \big[ \pi_k \ten id \big]
(x)^{\ten_m}$. Similarly, consider
$$\widehat{\pi}_{free}^k: x^{\ten_m} \in \mathbb{M} \mapsto \big[
\pi_{free}^k \ten id \big] (x)^{\ten_m} \in \mathbb{M}_{free}^k =
\ten_{sym}^m \big( \mathcal{M}_{free}^k \ten \mathcal{R} \big).$$

\vskip3pt

\noindent Let us note that
\begin{itemize}
\item $\sum_{j=1}^m \pi_{tens}^j(x) \in L_1(\mathbb{M})$.

\vskip3pt

\item According to i), if we set $$\mathbb{N} = \ten_{sym}^m \big(
\mathcal{N} \ten \mathcal{R} \big),$$ $(\mathbb{M}_k)_k$ and
$(\mathbb{M}_{free}^k)_k$ are independent top-subsymmetric copies
over $\mathbb{N}$.

\vskip3pt

\item Moreover, if we further assume that the $\mathcal{M}_k$'s
are symmetric, it also trivially follows that both the
$\mathbb{M}_k$'s and the $\mathbb{M}_{free}^k$'s are independent
symmetric copies of $\mathbb{M}$ over $\mathbb{N}$.
\end{itemize}
Under these conditions, we obtain from Theorem
\ref{Theorem-Araki-Woods} i)
\begin{eqnarray*}
\mathbb{E} \Big\| \sum_{j=1}^m \pi_{tens}^j \big( \sum_{i=1}^n
\eps_i \big[ \pi_i \ten id \big] (x) \big) \Big\|_1 & = &
\mathbb{E} \Big\| \sum_{i=1}^n \eps_i \widehat{\pi}_i
\big( \sum_{j=1}^m \pi_{tens}^j(x) \big) \Big\|_1 \\
& \sim & \mathbb{E} \Big\| \sum_{i=1}^n \eps_i
\widehat{\pi}_{free}^i \big( \sum_{j=1}^m \pi_{tens}^j(x) \big)
\Big\|_1 \\ & = & \mathbb{E} \Big\| \sum_{j=1}^m \pi_{tens}^j
\big( \sum_{i=1}^n \eps_i \big[ \pi_{free}^i \ten id \big] (x)
\big) \Big\|_1.
\end{eqnarray*}
For independent symmetric mean-zero copies we use Theorem
\ref{Theorem-Araki-Woods} ii) instead. \fin

\noindent {\bf Proof of Theorem A.} Given a complex Hilbert space
$\mathcal{H}$, let $\mathcal{F}_{-1}(\mathcal{H})$ denote its
antisymmetric Fock space. Write $c(e)$ and $a(e)$ for the creation
and annihilation operators associated to $e \in \mathcal{H}$.
Given an orthonormal basis $(e_{\pm k})_{k \ge 1}$ of
$\mathcal{H}$ and a family $(\mu_k)_{k \ge 1}$ of positive
numbers, the sequence $f_k = c(e_k) + \mu_k a(e_{-k})$ satisfies
the canonical anticommutation relations and we take $\mathcal{R}$
to be the von Neumann algebra generated by the $f_k$'s. Let
$\phi_{\mathcal{R}}$ be the quasi-free state on $\mathcal{R}$
determined by the vacuum. With this construction, $\mathcal{R}$ is
the Araki-Woods factor arising from the GNS construction applied
to the CAR algebra with respect to $\phi_{\mathcal{R}}$. In fact,
using a conditional expectation, we may replace the $\mu_k$'s by a
sequence $(\mu_k')_{k \ge 1}$ such that for every rational $0 <
\lambda < 1$ there are infinitely many $\mu_k' =
\lambda/(1+\lambda)$. According to Araki and Woods, we then obtain
the hyperfinite $\mathrm{III}_1$ factor $\mathcal{R}$.

\vskip5pt

Consider the amplification $\widehat{\mathcal{A}} = \mathcal{A}
\bar\ten \mathcal{B}(\ell_2)$ and assume for simplicity that
$\mathcal{A}$ is $\sigma$-finite. A normal strictly semifinite
faithful weight $\psi_{\widehat{\mathcal{A}}}$ is determined by a
sequence (a net in the general case) of pairs ($\psi_n, q_n$) such
that
\begin{itemize}
\item The $q_n$'s are increasing projections in
$\widehat{\mathcal{A}}$ with $\mathrm{SOT}-\lim_n q_n = 1$.

\item The $\psi_n$'s are normal positive functionals on
$\widehat{\mathcal{A}}$ with $\mathrm{supp} \hskip1pt \psi_n =
q_n$.

\vskip1pt

\item The $(\psi_n,q_n)$'s satisfy the compatibility condition
$\psi_{n+1} (q_n x q_n) = \psi_n(x)$.
\end{itemize}
Let us write $\mathrm{k}_n$ for the nondecreasing numbers
$\psi_n(q_n) \in (0,\infty)$. We refer to Propositions 8.10, 8.19
and 8.22 of \cite{JP2} for the fact that there is a normal
strictly semifinite faithful weight $\psi_{\widehat{\mathcal{A}}}$
with $\mathrm{k}_n \in \mathbb{N}$ and $\rho_n \in
L_1(\mathcal{B}(\ell_2) \bar\otimes \mathcal{R})$ such that
$$\|a\|_p \sim_c \lim_{n
\to \infty} \big\| a \ten \rho_n \big\|_{\mathcal{K}_{r \!
c_1}^1(\psi_n \otimes \phi_\mathcal{R})}$$ holds for every $a \in
L_p(\mathcal{A})$ up to an absolute constant $c$. The coefficients
$\rho_n$ are universal, i.e. do not depend on $\mathcal{A}$. The
specific form of $\psi_{\widehat{\mathcal{A}}}$ depends on the
spectrum of the operator $f_{\mid_{\partial_0}} \mapsto
f_{\mid_{\partial_1}}$ which takes the left boundary value of an
analytic function $f$ on the strip $0 < \mathrm{Re}z < 1$ to its
right boundary value, see Chapter 8 of \cite{JP2} for further
details. Writing $\psi_{\mathcal{R},n} = \psi_n \otimes
\phi_\mathcal{R}$, we know that $\psi_{\mathcal{R},n} =
\mathrm{k}_n \phi_{\mathcal{R},n}$ for some non-faithful state
$\phi_{\mathcal{R},n}$ with density $d_{\mathcal{R},n}$, so that
$\phi_{\mathcal{R},n}(x) = \mathrm{tr}(d_{\mathcal{R},n} x)$ and
the norm on the right takes the form
 $$\|z\|_{\mathcal{K}_{r \! c_1}^1(\psi_n \ten \phi_\mathcal{R})} =
 \inf_{z = z_1 + d_{\mathcal{R},n}^\frac12 z_r + z_c
 d_{\mathcal{R},n}^\frac12} \mathrm{k}_n \hskip1pt \|z_1\|_1 +
 \mathrm{k}_n^\frac12 \hskip2pt \|z_r\|_2
  + \mathrm{k}_n^\frac12 \hskip2pt
 \|z_c\|_2. $$
On the other hand, we have $\rho_n = \xi_n \ten \gamma$ where
$\gamma$ is a mean-zero element of $L_1(\mathcal{R},
\phi_{\mathcal{R}})$ and $\xi_n = \sum_{i,j \le n} e_{ij} \in
M_n$. Indeed, see \cite[Page 136 and proof of Proposition
8.10]{JP2} for the claim on $\gamma$ and the proofs of \cite[Lemma
8.21 and Proposition 8.22]{JP2} to get an idea of how to derive
the form of $\xi_n$, see also \cite{JP5} for more on this. Taking
this into account Theorem \ref{Theorem-Araki-Woods} ii) applies
for $\mathcal{N} = \langle \mathbf{1}_\mathcal{A} \rangle$ and
yields
\begin{eqnarray} \label{Eq1}
\quad \mathbb{E} \Big\| \summ_k \eps_k \hskip1pt \pi_k(x) \Big\|_p
& \sim & \lim_{n \to \infty} \mathbb{E} \Big\| \big( \summ_k
\eps_k \hskip1pt \pi_k(x) \big) \ten \rho_n \Big\|_{\mathcal{K}_{r
\! c_1}^1(\psi_n \ten \phi_\mathcal{R})} \\ \nonumber & \sim &
\lim_{n \to \infty} \mathbb{E} \Big\| \sum_{j=1}^n \pi_{tens}^j
\big( \summ_k \eps_k \hskip1pt \big[ \pi_k \ten id \big] (x \ten
\rho_n) \big) \Big\|_1.
\end{eqnarray}
Indeed, in the last step we use
\begin{itemize}
\item The mean-zero condition $$\phi_{\mathcal{R},n} \Big( \summ_k
\eps_k \hskip1pt \big[ \pi_k \ten id \big] (x \ten \rho_n) \Big) =
0.$$

\item The $\pi_{tens}^j$'s from Lemma \ref{Exchange} provide
symmetric independent copies.
\end{itemize}
The mean-zero condition follows from the fact that $\gamma$ is
mean-zero in $\mathcal{R}$. Now according to Lemma \ref{Exchange}
ii), we may replace $\pi_k$ by $\pi_{free}^k$ at the right hand
side of \eqref{Eq1}. This implies the first assertion calculating
backwards. For the second assertion, we first note that the
argument above is also valid without random signs whenever the
$\mathcal{M}_k$'s are symmetric and $\mathcal{E_N}(x) = 0$. In
other words, we have
\begin{equation} \label{Eq2}
\Big\| \summ_k \pi_k \big( x - \mathcal{E_N}(x) \big) \Big\|_p
\sim_c \Big\| \summ_k \pi_{free}^k \big( x - \mathcal{E_N}(x)\big)
\Big\|_p.
\end{equation}
Indeed, we just need to argue as we did above and use Lemma
\ref{Exchange} iii) for mean-zero symmetric copies. In the general
case, we may assume by approximation that only finitely many
$\pi_k$'s are considered. Then we have
\begin{eqnarray*}
\Big\| \sum_{k=1}^n \pi_k (x) \Big\|_p & \le & \Big\| \sum_{k=1}^n
\pi_k \big( x - \mathcal{E_N}(x) \big) \Big\|_p + \Big\|
\sum_{k=1}^n \pi_k \big( \mathcal{E_N}(x) \big) \Big\|_p \\ & = &
\Big\| \sum_{k=1}^n \pi_k \big( x - \mathcal{E_N}(x) \big)
\Big\|_p + \Big\| \sum_{k=1}^n \pi_{free}^k \big( \mathcal{E_N}(x)
\big) \Big\|_p \\ & \sim_c & \Big\| \sum_{k=1}^n \pi_{free}^k
\big( x - \mathcal{E_N}(x) \big) \Big\|_p + \Big\| \sum_{k=1}^n
\pi_{free}^k \big( \mathcal{E_N}(x) \big) \Big\|_p.
\end{eqnarray*}
We have used \eqref{Eq2} and the invariances
$${\pi_k}_{|_\mathcal{N}} = {\pi_{free}^k}_{|_\mathcal{N}} = id.$$
Moreover, using $\mathcal{E_N} = \mathsf{E}_{\mathcal{N}} \circ
\pi_{free}^k$ we find $\pi_{free}^k \circ \mathcal{E_N} =
\mathsf{E}_\mathcal{N} \circ \pi_{free}^k$ and $$\Big\|
\sum_{k=1}^n \pi_k (x) \Big\|_p \lesssim_c \Big\| \sum_{k=1}^n
\pi_{free}^k (x) \Big\|_p + 2 \Big\| \mathsf{E}_\mathcal{N} \big(
\sum_{k=1}^n \pi_{free}^k (x) \big) \Big\|_p \le 3 \Big\|
\sum_{k=1}^n \pi_{free}^k (x) \Big\|_p.$$ The reverse inequality
follows similarly.  \fin

\begin{remark}
\emph{If we allow the constant to be dependent on $p$, Theorem A
holds for $1 < p < \infty$ and not only for copies. Indeed, this
follows directly from the noncommutative Rosenthal inequality
\cite{JX4}. As in the classical case, we don't have a uniform
constant $c_p$ as $p \to \infty$. On the other hand, we know from
Claus K\"ostler \cite{Kost} that the symmetry assumption in the
second assertion of Theorem A can be replaced by the weaker notion
of subsymmetric copies defined in the Introduction.}
\end{remark}

\begin{remark}
\emph{Let $\mathsf{M}$ be a finite von Neumann algebra equipped
with a faithful normal trace $\tau$ and let $u_1, u_2, \ldots, u_n
\in \mathsf{M}$ be unitaries. We claim that the inequality below
holds for $1 \le p \le 2$, independent top-subsymmetric copies and
constants independent of $p$
\begin{eqnarray} \label{Unitaries}
\lefteqn{\hskip5pt \mathbb{E} \Big\| \sum_{k=1}^n \varepsilon_k
\pi_k(x) \ten u_k \Big\|_{L_p(\mathcal{A} \bar\ten \mathsf{M})}}
\\ \nonumber \null \qquad \null & \sim_c & \inf_{x = a+b+c} n^{\frac1p}
\|a\|_{L_p(\mathcal{M})} + \sqrt{n} \big\| \mathcal{E_N}
(bb^*)^{\frac12} \big\|_{L_p(\mathcal{N})} + \sqrt{n} \big\|
\mathcal{E_N}(c^*c)^{\frac12} \big\|_{L_p(\mathcal{N})}.
\end{eqnarray}
Indeed, the case $p=1$ is a more general version of Theorem
\ref{Theorem-Araki-Woods} that was already considered in
\cite{J3}. Moreover, it can be easily checked that our arguments
in Lemma \ref{Exchange} and Theorem A are stable under taking
tensors with arbitrary unitaries, as far as we work with this
refinement of Theorem \ref{Theorem-Araki-Woods}. This gives
$$\mathbb{E} \Big\| \summ_k \eps_k \hskip1pt \pi_k(x) \ten u_k
\Big\|_{L_p(\mathcal{A} \bar\ten \mathsf{M})} \sim_c \mathbb{E}
\Big\| \summ_k \eps_k \hskip1pt \pi_{free}^k(x) \ten u_k
\Big\|_{L_p(\mathcal{A} \bar\ten \mathsf{M})}.$$ If we now combine
it with the free Rosenthal inequality from \cite{JPX}
\begin{eqnarray*}
\lefteqn{\mathbb{E} \Big\| \sum_{k=1}^n \varepsilon_k \pi_{free}^k
(x_k) \Big\|_{L_p(\mathcal{A})}} \\ & \sim_c & \inf_{x_k =
a_k+b_k+c_k} \Big( \sum_{k=1}^n \|a_k\|_p^p \Big)^{\frac1p} +
\Big\| \big( \sum_{k=1}^n \mathcal{E_N} (b_kb_k^*) \big)^{\frac12}
\Big\|_p + \Big\| \big( \sum_{k=1}^n \mathcal{E_N}(c_k^*c_k)
\big)^{\frac12} \Big\|_p,
\end{eqnarray*}
we easily end up with \eqref{Unitaries}. This is a form of the
noncommutative Rosenthal inequality \cite{JX4} for independent
top-subsymmetric copies with no singularity at $p=1$ and will be
instrumental in the following paragraph.}
\end{remark}

\noindent {\bf Proof of Corollary A1.} Inequality \eqref{Unitaries}
implies the assertion by a central limit procedure which allows us
to pass from the three terms in Rosenthal inequality to the two
terms in the assertion. The argument for the particular case $p=1$
was given in Section 8 of \cite{J3}. The case $1 < p \le 2$ just
requires simple modifications that we leave to the reader. \fin

\noindent {\bf Proof of Corollary A2.} The construction in
\cite{JP2} is sketched as follows $$L_q (\mathcal{M})
\stackrel{(\alpha)}{\longrightarrow} \Big(
\mathcal{H}_{2p',2}^r(\mathcal{M}, \theta) \ten_{\mathcal{M}, h}
\mathcal{H}_{2p',2}^c(\mathcal{M}, \theta) \Big)_*
\stackrel{(\beta)}{\longrightarrow} \mathcal{K}_{r\! c_p}^p(\phi
\ten \psi \ten \xi) \stackrel{(\gamma)}{\longrightarrow}
L_p(\mathcal{A}).$$ We now review the cb-embeddings $(\alpha)$,
$(\beta)$ and $(\gamma)$ in some detail to identify where the
singularity appears. It will happen once in $(\beta)$ and once in
$(\gamma)$. However, in both cases our version \eqref{Unitaries} of
the noncommutative Rosenthal inequality for independent copies ---no
need of unitaries this time--- will allow us to remove the
singularity. Namely, the embedding $(\alpha)$ generalizes the
so-called \lq Pisier's exercise\rq${}$ to embed the Schatten class
$S_q = C_q \ten_h R_q$ into an operator space of the form
$$\big( C_p \oplus \mathrm{OH} \big) / graph(\Lambda_1)^\perp
\otimes_h \big( R_p \oplus \mathrm{OH} \big) /
graph(\Lambda_2)^\perp,$$ see \cite{JP4,X3} for further details.
When working with general von Neumann algebras this requires to
encode complex interpolation in terms of certain spaces of
analytic functions. This follows from the factorization
$$L_q(\mathcal{M}) = L_{2q}^r(\mathcal{M}) \ten_{\mathcal{M},h}
L_{2q}^c(\mathcal{M}),$$ Proposition 8.19 in \cite{JP2} and
duality. It turns out that the constants are independent of $p,q$
at this step. Moreover, a more convenient way to write the space
between $(\alpha)$ and $(\beta)$ is by means of a 4-term sum
$\mathcal{K}_{p,2}(\psi,\xi)$, with $\psi$ and $\xi$ normal
strictly semifinite faithful weights on $\mathcal{M}$ and
$\mathcal{B}(\ell_2)$ respectively, see Proposition 2.22
of \cite{JP2} for further details. The space
$\mathcal{K}_{p,2}(\psi \ten \xi)$ arises from a direct limit
$$\mathcal{K}_{p,2}(\psi \ten \xi) = \overline{\bigcup_{n \ge 1}
\mathcal{K}_{p,2}(\psi_n \ten \xi_n)},$$ where $\psi_n, \xi_n$ are
the restrictions of $\psi, \phi$ to certain increasing sequences
of finite projections. If $\psi_n \otimes \xi_n = \mathrm{k}_n \phi_n$ with
$\phi_n$ a state supported by the projection $q_n$, it turns out that $$\mathcal{K}_{p,2}(\psi_n \otimes \xi_n) = \sum_{u,v \in \{2p,4\}} \mathrm{k}_n^{\frac{1}{u} + \frac{1}{v}} L_{(u,v)}(q_n (\mathcal{M} \bar\otimes \mathcal{B}(\ell_2)) q_n).$$ The definition of asymmetric $L_{(u,v)}$ spaces is recalled in Section \ref{Section2} below. On the other hand, the space $\mathcal{K}_{r \!
c_p}^p(\phi \ten \psi \ten \xi)$ is also the direct limit of an
increasing family of 3-term sums $$\mathcal{K}_{r \! c_p}^p(\psi
\ten \phi \ten \xi) = \overline{\bigcup_{n \ge 1} \mathcal{K}_{r
\! c_p}^p(\psi_n \ten \phi \ten \xi_n)},$$ where $\phi$ is the
quasi-free state $\phi_\mathcal{R}$ in the proof of Theorem A above. The complete embedding $\mathcal{K}_{p,2}(\psi_n \ten
\xi_n) \to \mathcal{K}_{r \! c_p}^p(\psi_n \ten \phi \ten \xi_n)$
holds up to constants independent on $n$, see Proposition 8.10 in
\cite{JP2}. The main ingredient in the argument is the
noncommutative Rosenthal inequality for independent copies in
$L_p$. Equipped with \eqref{Unitaries}, we may now provide a
universal constant valid for any $1 \le p \le 2$. Thus the
embedding $(\beta)$ holds up to constants independent of $p$. Let
us consider the embedding $(\gamma)$. Its construction is given in
\cite[Theorem 8.11]{JP2}. The idea is to construct cb-embeddings
of $\mathcal{K}_{r \! c_p}^p(\psi_n \ten \phi \ten \xi_n)$ for
each $n \ge 1$, so that we can take direct limits and
hyperfiniteness is preserved. The main tool is a noncommutative
Poisson process, an algebraic construction that has nothing to
do with constants. Therefore the problem reduces to cb-embed
$\mathcal{K}_{r \! c_p}^p(\psi_n \ten \phi \ten \xi_n)$ into
$L_p(\mathcal{A}_n)$. A quick look at \cite[Theorem 8.11]{JP2}
shows that the only point where the singularity appears is once
more from the use of noncommutative Rosenthal inequality for
independent copies. Hence, Theorem A applies and produces
\eqref{Unitaries}, implying the assertion. \fin

\begin{remark}
\emph{According to \cite{JP4,JP2}, the assertion in Corollary A2
also holds with $L_q(\mathcal{M})$ replaced by any operator space
of the form $\mathrm{X}_1 \ten_h \mathrm{X}_2$, with
$\mathrm{X}_1$ a quotient of a subspace of $R \oplus \mathrm{OH}$
and $\mathrm{X}_2$ a quotient of a subspace of $C \oplus
\mathrm{OH}$.}
\end{remark}

\section{Mixed-norms of free variables}
\label{Section2}

In this section we recall several results from \cite{JP2} for the
convenience of the reader. The main result is a variation of the
free Rosenthal inequality from \cite{JPX} which will be instrumental in
the course of our argument. The correct formulation involves certain
noncommutative function spaces.

\subsection{Conditional $L_p$ spaces}

Inspired by Pisier's theory \cite{P2} of noncommutative
vector-valued $L_p$ spaces, several noncommutative function spaces
have been recently introduced in quantum probability. The first
insight came from some of Pisier's fundamental equalities which we
briefly review. Let $\mathcal{N}_1,\mathcal{N}_2$ be hyperfinite
von Neumann algebras. Given $1 \le p,q \le \infty$, we define $1/r
= |1/p - 1/q|$. If $p \le q$, the norm of $x \in L_p
(\mathcal{N}_1; L_q(\mathcal{N}_2))$ is given by
$$\inf \Big\{ \|\alpha\|_{L_{2r}(\mathcal{N}_1)}
\|y\|_{L_q(\mathcal{N}_1 \bar\otimes \mathcal{N}_2)}
\|\beta\|_{L_{2r}(\mathcal{N}_1)} \, \big| \ x = \alpha y \beta
\Big\}.$$ If $p \ge q$, the norm of $x \in L_p (\mathcal{N}_1;
L_q(\mathcal{N}_2))$ is given by
$$\sup \Big\{ \|\alpha x \beta \|_{L_q(\mathcal{N}_1 \bar\otimes
\mathcal{N}_2)} \, \big| \ \alpha, \beta \in
\mathsf{B}_{L_{2r}(\mathcal{N}_1)} \Big\}.$$

The hyperfiniteness is an essential assumption in \cite{P2}.
However, when dealing with mixed $L_p(L_q)$ norms, Pisier's
identities remain true for general von Neumann algebras, see
\cite{J1,JX4}. On the other hand, given any von Neumann algebra
$\mathcal{M}$, the \emph{row} and \emph{column} subspaces of $L_p$
are defined as follows
$$R_p^n(L_p(\mathcal{M})) = \Big\{ \sum_{k=1}^n x_k \ten e_{1k} \,
\big| \ x_k \in L_p(\mathcal{M}) \Big\} \subset L_p \big(
\mathcal{M} \bar\ten \mathcal{B}(\ell_2) \big),$$
$$C_p^n(L_p(\mathcal{M})) = \Big\{ \sum_{k=1}^n x_k \ten
e_{k1} \, \big| \ x_k \in L_p(\mathcal{M}) \Big\} \subset L_p
\big( \mathcal{M} \bar\ten \mathcal{B}(\ell_2) \big),$$ where
$(e_{ij})$ denotes the unit vector basis of $\mathcal{B}(\ell_2)$.
These spaces are crucial in the noncommutative
Khintchine/Rosenthal type inequalities and in noncommutative
martingale inequalities, where the row and column spaces are
traditionally denoted by $L_p(\mathcal{M}; \ell_2^r)$ and
$L_p(\mathcal{M}; \ell_2^c)$. The norm in these spaces is given by
$$\begin{array}{c} \displaystyle \Big\| \sum_{k=1}^n x_k \otimes
e_{1k} \Big\|_{R_p^n(L_p(\mathcal{M}))} = \Big\| \Big(
\sum_{k=1}^n x_k x_k^* \Big)^{\frac12} \Big\|_{L_p(\mathcal{M})},
\\ [12pt] \displaystyle \Big\| \sum_{k=1}^n x_k \otimes e_{k1}
\Big\|_{C_p^n(L_p(\mathcal{M}))} = \Big\| \Big( \sum_{k=1}^n x_k^*
x_k \Big)^{\frac12} \Big\|_{L_p(\mathcal{M})}.
\end{array}$$ In what follows we write
$R_p^n(L_p(\mathcal{M})) = L_p(\mathcal{M}; R_p^n)$ and
$C_p^n(L_p(\mathcal{M})) = L_p(\mathcal{M}; C_p^n)$.

\vskip3pt

Now, let us assume that $\mathcal{N}$ is a von Neumann subalgebra
of $\mathcal{M}$ and that there exists a normal faithful
conditional expectation $\mathcal{E_N}: \mathcal{M} \to
\mathcal{N}$. Then we may define $L_p$ norms of the
\emph{conditional square functions}
$$\Big( \sum_{k=1}^n \mathcal{E_N}(x_k x_k^*) \Big)^{\frac12} \quad
\mbox{and} \quad \Big( \sum_{k=1}^n \mathcal{E_N}(x_k^* x_k)
\Big)^{\frac12}.$$ These norms must be properly defined for $1 \le
p \le 2$, see \cite{J1} or \cite[Chapter 1]{JP2}. The resulting
spaces coincide with the row/column spaces above if $\mathcal{N}$
is $\mathcal{M}$ itself. When $n=1$ we recover the spaces
$L_p^r(\mathcal{M}, \mathsf{E})$ and $L_p^c(\mathcal{M},
\mathsf{E})$ from \cite{J1}.

\vskip3pt

We have already introduced $L_p(L_q)$ spaces, row and column
subspaces of $L_p$ and some variations associated to a given
conditional expectation. All the norms considered so far fit into
more general noncommutative function spaces ---for not necessarily
hyperfinite von Neumann algebras--- which we now define. Consider
the solid $\mathsf{K}$ in $\R^3$ determined by
$$\mathsf{K} = \Big\{(1/u,1/v,1/q) \, \big| \ 2 \le u,v \le
\infty, \ 1 \le q \le \infty, \ 1/u + 1/q + 1/v \le 1 \Big\}.$$
Let $\mathcal{N}$ be a conditioned subalgebra of $\mathcal{M}$ and
take $1/p = 1/u + 1/q + 1/v$ for some $(1/u,1/v,1/q) \in
\mathsf{K}$. Then we define the \emph{amalgamated $L_p$ space}
associated to the indices $(u,q,v)$ as the subspace
$L_u(\mathcal{N}) L_q(\mathcal{M}) L_v(\mathcal{N})$ of
$L_p(\mathcal{M})$ equipped with the norm
$$\inf \Big\{ \|\alpha\|_{L_u(\mathcal{N})}
\|y\|_{L_q(\mathcal{M})} \|\beta\|_{L_v(\mathcal{N})} \, \big| \ x
= \alpha y \beta \Big\},$$ where the infimum runs over all
possible factorizations $x = \alpha y \beta$ with
$(\alpha,y,\beta)$ belonging to $L_u(\mathcal{N}) \times
L_q(\mathcal{M}) \times L_v(\mathcal{N})$. Let us now fix
$(1/u,1/v,1/p) \in \mathsf{K}$ and take $1/s = 1/u+1/p+1/v$. Then
we define the \emph{conditional $L_p$ space} associated to the
indices $(u,v)$ as the completion $L_u^{-1}(\mathcal{N})
L_s(\mathcal{M}) L_v^{-1}(\mathcal{N})$ of $L_p(\mathcal{M})$ with
respect to the norm
$$\sup \Big\{ \|\alpha x \beta\|_{L_s(\mathcal{M})} \, \big| \
\|\alpha\|_{L_u(\mathcal{N})}, \|\beta\|_{L_v(\mathcal{N})} \le 1
\Big\}.$$ Both, amalgamated and conditional $L_p$ spaces, where
introduced in \cite{JP2} and we refer to that paper for a more
detailed exposition. It should also be noticed that our
terminology $L_u^{-1}(\mathcal{N}) L_s(\mathcal{M})
L_v^{-1}(\mathcal{N})$ for conditional $L_p$ spaces is different
from the one used in \cite{JP2}. Now we collect the complex
interpolation and duality properties of amalgamated and
conditional $L_p$ spaces from \cite{JP2}. Our interpolation
identities generalize some previous results by Pisier \cite{P0}
and recently by Xu \cite{X}.

\vskip3pt

\noindent Let $\mathsf{K}_0$ denote the interior of $\mathsf{K}$.
Then we have:

\vskip3pt

\begin{itemize}
\item[\textbf{a)}] $L_u(\mathcal{N})
L_q(\mathcal{M}) L_v(\mathcal{N})$ is a Banach space.

\vskip3pt

\item[\textbf{b)}] $L_{u_{\theta}}(\mathcal{N})
L_{q_{\theta}}(\mathcal{M}) L_{v_{\theta}}(\mathcal{N})$ is
isometrically isomorphic to
$$\Big[L_{u_0}(\mathcal{N}) L_{q_0}(\mathcal{M})
L_{v_0}(\mathcal{N}), L_{u_1}(\mathcal{N}) L_{q_1}(\mathcal{M})
L_{v_1}(\mathcal{N}) \Big]_{\theta}^{\null},$$ with $(\frac{1}{u_\theta},
\frac{1}{q_\theta}, \frac{1}{v_\theta}) = (\frac{1-\theta}{u_0} +
\frac{\theta}{u_1}, \frac{1-\theta}{q_0} + \frac{\theta}{q_1},
\frac{1-\theta}{v_0} + \frac{\theta}{v_1})$.

\vskip3pt

\item[\textbf{c)}] If $(1/u,1/v,1/q) \in \mathsf{K}_0$ and $1 -
1/p = 1/u + 1/q + 1/v$ $$\big( L_u(\mathcal{N}) L_q(\mathcal{M})
L_v(\mathcal{N}) \big)^* = L_u^{-1}(\mathcal{N})
L_{q'}(\mathcal{M}) L_v^{-1}(\mathcal{N}),$$
$$\big( L_u^{-1}(\mathcal{N})
L_{q'}(\mathcal{M}) L_v^{-1}(\mathcal{N}) \big)^* =
L_u(\mathcal{N}) L_q(\mathcal{M}) L_v(\mathcal{N}),$$ with respect to the antilinear duality bracket $\langle x,y \rangle = \mathrm{tr}(x^*y)$. A natural way to read the first identity (the second one is its dual) is to say that the dual of the amalgamated $L_{p'}$ space associated to $(u,v)$ is the conditional $L_p$ space associated to $(u,v)$, since $$1/p' = 1/u + 1/q + 1/v \quad \mbox{and} \quad 1/p = 1/q' - 1/u - 1/v.$$
\end{itemize}

\vskip3pt

\noindent We refer the reader to Part I of \cite{JP2} for some
refinements of these results.

\subsection{A variant of free Rosenthal's inequality}

In this paragraph we formulate the free analogue of inequality
\eqref{Mix>2} in the Introduction and its dual. To be precise, we
shall work for convenience with i.d. variables. In that case, it
is easily checked that \eqref{Mix>2} provides a natural way to
realize the space $$\mathcal{J}_{p,q}^n(\Omega) = n^{\frac1p}
L_p(\Omega) \cap n^{\frac1q} L_q(\Omega)$$ as an isomorph of a
subspace of $L_p(\Omega;\ell_q^n)$. Quite surprisingly, replacing
in \eqref{Mix>2} independent variables by matrices of independent
variables requires to intersect four spaces using the so-called
\emph{asymmetric} $L_p$ spaces. This phenomenon was discovered for
the first time in \cite{JP} and is partly motivated by the
isometry $L_p = L_{2p} L_{2p}$ meaning that the $p$-norm of $f$ is
the infimum of $\|g\|_{2p} \|h\|_{2p}$ over all factorizations $f
= g h$. Namely, if $L_{2p}^r$ and $L_{2p}^c$ denote the row and
column quantizations of $L_{2p}$ determined by definition
\eqref{Eq-RCLp} below, the operator space analogue of this
isometry is given by the complete isometry $L_p = L_{2p}^r
L_{2p}^c$, see below for further details. This leads us to
redefine $\mathcal{J}_{p,q}^n$ as
$$\mathcal{J}_{p,q}^n = \Big( n^{\frac{1}{2p}} L_{2p}^r \cap
n^{\frac{1}{2q}} L_{2q}^r \Big) \Big( n^{\frac{1}{2p}} L_{2p}^c
\cap n^{\frac{1}{2q}} L_{2q}^c \Big).$$ According to \cite{JP2},
we find
\begin{equation} \label{Eq-Aspect}
\mathcal{J}_{p,q}^n = n^{\frac{1}{p}} L_{2p}^r L_{2p}^c \cap
n^{\frac{1}{2p}+\frac{1}{2q}} L_{2p}^r L_{2q}^c \cap
n^{\frac{1}{2q} + \frac{1}{2p}} L_{2q}^r L_{2p}^c \cap
n^{\frac{1}{q}} L_{2q}^rL_{2q}^c.
\end{equation}
These spaces will be rigorously defined below. Our only aim here is
to motivate the forthcoming results and definitions. Let us now see
how the space in \eqref{Eq-Aspect} generalizes our first definition
of $\mathcal{J}_{p,q}^n(\Omega)$. On the Banach space level we have
the isometries $L_{2p}^r L_{2q}^c = L_s = L_{2q}^r L_{2p}^c$ with
$1/s = 1/2p + 1/2q$. Moreover, again by H\"{o}lder inequality it is
clear that
$$n^{\frac1s} \|f\|_s \le \max \Big\{ n^{\frac1p}
\|f\|_p, n^{\frac1q} \|f\|_q \Big\}$$ and the two cross terms in
\eqref{Eq-Aspect} disappear. However, in the category of operator
spaces the four terms have a significant contribution. The
operator space/free version of \eqref{Mix>2} is the main result in
\cite{JP2}, and goes further than its commutative counterpart.
More precisely, in contrast with the classical case, we find a
right formulation for ($\Sigma_{\infty q}$). Indeed, as for
Khintchine and Rosenthal inequalities, the limit case as $p \to
\infty$ holds when replacing independence by freeness.

\vskip5pt

Now we give detailed definitions and results. Let us write
$L_2^r(\mathcal{M})$ and $L_2^c(\mathcal{M})$ for the row/column
quantizations of $L_2(\mathcal{M})$ and let $2 \le q \le \infty$.
Then, the row/column structures on $L_q(\mathcal{M})$ are defined
as follows
\begin{equation} \label{Eq-RCLp}
\begin{array}{rcl}
L_q^r(\mathcal{M}) & = & \big[ \mathcal{M},
L_2^r(\mathcal{M})\big]_{\frac{2}{q}}, \\ [5pt] L_q^c(\mathcal{M})
& = & \big[ \mathcal{M}, L_2^c(\mathcal{M})\big]_{\frac{2}{q}}.
\end{array}
\end{equation}
In fact, a rigorous definition should take Kosaki's embeddings into account as done in \cite[Identity (1.3)]{JP2}, but we shall ignore such formalities. Now, if $2 \le u,v \le \infty$ and $1/p = 1/u + 1/v$ for some $1
\le p \le \infty$, we define the \emph{asymmetric $L_p$ space}
\label{LpAsimetrico} associated to the pair $(u,v)$ as the
$\mathcal{M}$-amalgamated Haagerup tensor product
\begin{equation} \label{Eq-Asymmetric}
L_{(u,v)}(\mathcal{M}) = L_u^r(\mathcal{M})
\otimes_{\mathcal{M},h} L_v^c(\mathcal{M}).
\end{equation}
That is, we consider the quotient of $L_u^r(\mathcal{M}) \otimes_h
L_v^c(\mathcal{M})$ by the closed subspace $\mathcal{I}$ generated
by the differences $x_1 \gamma \otimes x_2 - x_1 \otimes \gamma
x_2$ with $\gamma \in \mathcal{M}$. By a well known factorization
argument the norm of an element
$x$ in $L_{(u,v)}(\mathcal{M})$ is given by $\|x\|_{(u,v)}^{\null}
= \inf_{x = \alpha \beta} \|\alpha\|_{L_u(\mathcal{M})}
\|\beta\|_{L_v(\mathcal{M})}$, see Lemma 1.9 in \cite{JP2}.

\begin{itemize}
\item If $\mathcal{M} = M_m$, the space in \eqref{Eq-Asymmetric}
reduces to $S_{(u,v)}^m = C_{u/2}^m \otimes_h R_{v/2}^m$.

\item We have a cb-isometry $L_p(\mathcal{M}) =
L_{(2p,2p)}(\mathcal{M})$, see \cite[Remark 7.5]{JP2}.
\end{itemize}

Let $1 \le q \le p \le \infty$. According to the discussion which
led to \eqref{Eq-Aspect}, we know how the general aspect of
$\mathcal{J}_{p,q}^n(\mathcal{M})$ should be. Now, equipped with
asymmetric $L_p$ spaces we obtain a factorization of  noncommutative
$L_p$ spaces in the right way:
$$\mathcal{J}_{p,q}^n(\mathcal{M}) = \bigcap_{u,v \in \{2p,2q\}}
n^{\frac{1}{u} + \frac{1}{v}} \, L_{(u,v)}(\mathcal{M}).$$ If we
take $$\mathcal{M}_m = M_m(\mathcal{M}) \quad \mbox{and} \quad
\mathsf{E}_m = id_{M_m} \otimes \varphi: \mathcal{M}_m \to M_m$$
for $m \ge 1$ and define $$\frac1r = \frac1q - \frac1p \quad
\mbox{and} \quad \frac{1}{\gamma(u,v)} = \frac{1}{u} + \frac{1}{p}
+ \frac{1}{v},$$ we have an isometry
\begin{equation} \label{Lem-oss-Jpq}
S_p^m \big( \mathcal{J}_{p,q}^n(\mathcal{M}) \big) = \bigcap_{u,v
\in \{2r,\infty\}}^{\null} n^{\frac{1}{\gamma(u,v)}} \,
L_u^{-1}(M_m) L_{\gamma (u,v)}(\mathcal{M}_m) L_v^{-1}(M_m).
\end{equation}
The proof can be found in \cite{JP2}. Let $\mathcal{N}$ be a
conditioned subalgebra of $\mathcal{M}$ with corresponding
conditional expectation $\mathcal{E_N}: \mathcal{M} \to
\mathcal{N}$. According to \eqref{Lem-oss-Jpq}, we define the
$\mathcal{J}$-spaces
$$\mathcal{J}_{p,q}^n(\mathcal{M}, \mathcal{E_N}) = \bigcap_{u,v \in
\{2r,\infty\}}^{\null} n^{\frac{1}{\gamma(u,v)}} \,
L_u^{-1}(\mathcal{N}) L_{\gamma(u,v)}(\mathcal{M})
L_v^{-1}(\mathcal{N}).$$ The isometry \eqref{Lem-oss-Jpq} shows us
the way to follow. The philosophy is that complete boundedness
arises from amalgamation, see \cite{JP4,JP2}. Indeed, instead of
working with the o.s.s. of the spaces
$\mathcal{J}_{p,q}^n(\mathcal{M})$, it suffices to argue with the
Banach space structure of the more general spaces
$\mathcal{J}_{p,q}^n(\mathcal{M}, \mathcal{E_N})$. In this spirit,
for $1 \le q \le p \le \infty$ we set $1/r = 1/q - 1/p$ and
introduce the spaces
\begin{eqnarray*}
\mathcal{R}_{2p,q}^n(\mathcal{M}, \mathcal{E_N}) & = &
n^{\frac{1}{2p}} \, L_{2p}(\mathcal{M}) \, \cap \,
n^{\frac{1}{2q}} \, L_{2r}^{-1}(\mathcal{N}) L_{2q}(\mathcal{M})
L_\infty^{-1}(\mathcal{N}), \\ \mathcal{C}_{2p,q}^n \,
(\mathcal{M}, \mathcal{E_N}) & = & n^{\frac{1}{2p}} \,
L_{2p}(\mathcal{M}) \, \cap \, n^{\frac{1}{2q}} \, L_{\infty}^{-1}
(\mathcal{N}) L_{2q}(\mathcal{M}) L_{2r}^{-1}(\mathcal{N}).
\end{eqnarray*}

\begin{remark} \label{Re-oM}
\emph{The notion of $\mathcal{M}$-amalgamated Haagerup tensor product $\mathrm{X}_1 \otimes_{\mathcal{M},h} \mathrm{X}_2$ extends naturally to any pair $(\mathrm{X}_1, \mathrm{X}_2)$ of operator spaces such that $\mathrm{X}_1$ contains $\mathcal{M}$ as a right ideal and $\mathrm{X}_2$ does it as a left ideal. We shall write $\mathrm{X}_1 \oM \mathrm{X}_2$ to denote the underlying Banach space structure of $\mathrm{X}_1 \otimes_{\mathcal{M},h} \mathrm{X}_2$. According to the definition of the Haagerup tensor product and recalling the isometric embeddings $\mathrm{X}_j \subset \mathcal{B}(\mathcal{H}_j)$, we have $$\|x\|_{\mathrm{X}_1 \oM \mathrm{X}_2} = \inf \left\{ \Big\| \big( \summ_k x_{1k} x_{1k}^* \big)^{1/2} \Big\|_{\mathcal{B}(\mathcal{H}_1)} \Big\| \big( \summ_k x_{2k}^* x_{2k} \big)^{1/2} \Big\|_{\mathcal{B}(\mathcal{H}_2)} \right\},$$ where the infimum runs over all possible decompositions of $x + \mathcal{I}$ into a finite sum $$x = \summ_k x_{1k} \otimes x_{2k}
+ \mathcal{I}.$$ This uses the o.s.s. of $\mathrm{X}_j$ since row/column square functions live in
$\mathcal{B}(\mathcal{H}_j)$ but not necessarily in $\mathrm{X}_j$. Therefore, if $\mathrm{X}_1$ is stable under row square functions (finite sums) and $\mathrm{X}_2$ is stable under column square functions, no o.s.s. on the $\mathrm{X}_j$'s is needed to define $\mathrm{X}_1 \otimes_{\mathcal{M}} \mathrm{X}_2$. In particular, we may consider the Banach space $$\mathcal{R}_{2p,q}^n(\mathcal{M}, \mathcal{E_N}) \oM \mathcal{C}_{2p,q}^n(\mathcal{M}, \mathcal{E_N}).$$ In what follows, we might abuse of the terminology by writing $\mathrm{X}_1 \otimes_{\mathcal{M},h} \mathrm{X}_2$ for the space $\mathrm{X}_1 \otimes_\mathcal{M} \mathrm{X}_2$. If that happens, it will be clear what space do we mean from the context. The reader can find some of the basic properties of this construction in \cite[Chapters 6 and 7]{JP2}.}
\end{remark}

\begin{remark} \label{AmalgTrick}
\emph{The cb-isometry $$\mathrm{X}_1 \otimes_{\mathcal{M},h} \mathrm{X}_2 \, = \, \big( \mathrm{X}_1 \otimes_h R_m \big) \otimes_{M_m(\mathcal{M}),h} \big( C_m \otimes_h \mathrm{X}_2 \big)$$ reflects the behavior of row/column operator spaces with respect to amalgamated tensor products and it is a key property that we will be using along the paper. Indeed, it suffices to understand that $$\dim R_m \otimes_{M_m,h} C_m = 1,$$ which follows from the fact that $e_{1i} \otimes e_{j1} \sim \delta_{ij} e_{11} \otimes e_{11}$ where $\sim$ refers to the equivalence relation imposed by the quotient of amalgamation. These equivalences can be easily justified by the reader.}
\end{remark}

\noindent The isomorphisms below are the key results in
\cite{JP2}:
\begin{itemize}
\item[\textbf{a)}] If $1 \le q \le p \le \infty$, we have
$$\mathcal{J}_{p,q}^n(\mathcal{M}, \mathcal{E_N}) \simeq
\mathcal{R}_{2p,q}^n(\mathcal{M}, \mathcal{E_N}) \oM
\mathcal{C}_{2p,q}^n(\mathcal{M}, \mathcal{E_N}).$$

\item[\textbf{b)}] If $1 \le p \le \infty$ and $1/q = 1-\theta +
\theta/p$, we have
$$\mathcal{J}_{p,q}^n(\mathcal{M}, \mathcal{E_N}) \simeq \big[
\mathcal{J}_{p,1}^n(\mathcal{M}, \mathcal{E_N}),
\mathcal{J}_{p,p}^n(\mathcal{M}, \mathcal{E_N}) \big]_{\theta}.$$
\end{itemize}
Moreover, in all cases the involved relevant constants are
independent of $n$.

\begin{remark} \label{Singpq}
\emph{It is worth mentioning that the constants for $$\mathcal{J}_{p,q}^n(\mathcal{M}, \mathcal{E_N}) \simeq \big[ \mathcal{J}_{p,1}^n(\mathcal{M}, \mathcal{E_N}), \mathcal{J}_{p,p}^n(\mathcal{M}, \mathcal{E_N}) \big]_{\theta}$$ obtained in \cite[Theorem 7.2]{JP2} have a singularity when $(p,q) \sim (\infty,1)$. To be more explicit, for $q$ small and $p$ large we obtain a constant $c_{p,q} \sim (p-q)/(pq+q-p)$ and the same singularity appears in Theorem \ref{Theorem-Sigma-pq} below. However, the assertion in such theorem holds for the extremal values $(p,q)=(\infty,1)$ and the singularity seems to be removable. It appears as a byproduct of the noncommutative Burkholder inequality, which is used in the argument, see \cite[Remark 7.11]{JP2}. Fortunately, our construction of operator space $L_p$ embeddings in \cite{JP4,JP2} only uses Theorem \ref{Theorem-Sigma-pq} for $q=2$ and no singularity occurs in that case.}
\end{remark}

\newcommand{\M}{\mathcal{M}}
\newcommand{\NN}{\mathcal{N}}

In what follows, we shall also need to work with the predual space
of $\mathcal{J}_{p,q}^n(\mathcal{M}, \mathcal{E_N})$ which is
defined as follows. Given $1 \le p \le q \le \infty$, let us
consider the index $\frac1r = \frac1p - \frac1q$ and the
coefficient $\rho(u,v)$ determined by $\frac{1}{\rho(u,v)} =
\frac1p - \frac1u - \frac1v$. Then, we define the space
$$\mathcal{K}_{p,q}^n(\mathcal{M}, \mathcal{E_N}) =
\sum_{u,v \in \{2r,\infty\}} n^{\frac{1}{\rho(u,v)}}
L_u(\mathcal{N}) L_{\rho(u,v)}(\mathcal{M}) L_v(\mathcal{N}).$$
The antilinear bracket $\langle x,y \rangle_n = n \,
\mathrm{tr}(xy^*)$ gives $\mathcal{K}_{p,q}^n(\mathcal{M},
\mathcal{E_N})^* = \mathcal{J}_{p',q'}^n(\mathcal{M},
\mathcal{E_N})$, see \cite[Remark 7.4]{JP2}. As for the
$\mathcal{J}$-spaces, we shall write
$\mathcal{K}_{p,q}^n(\mathcal{M})$ to denote the space
$\mathcal{K}_{p,q}^n(\mathcal{M}, \mathcal{E_N})$ with
$(\mathcal{N}, \mathcal{E_N}) = (\langle \mathbf{1}_\mathcal{M}
\rangle, \varphi)$. Having defined the noncommutative forms of our
$\mathcal{J}$ and $\mathcal{K}$ spaces, we can now
state the free analog of $(\Sigma_{pq})$ and its dual  from
\cite{JP2}. It is fortunate that in the free case a more general
statement holds not requiring copies.

\begin{theorem} \label{Theorem-Sigma-pq}
If $1 \le p \le q \le \infty$, the maps
$$\begin{array}{rcl}
\displaystyle \sum_{k=1}^n x_k \ten \delta_k \in \hskip2.3pt
\mathcal{K}_{p,q}^1 \Big( \ell_{\infty}^n(\mathcal{M}),
\frac{1}{n} \sum_{k=1}^n \mathcal{E_N} \Big) \hskip1.3pt & \!\!
\mapsto \!\! & \displaystyle \sum_{k=1}^n \pi_{free}^k(x_k) \ten
\delta_k \in \hskip2.3pt L_p(\mathcal{A}_{free}; \ell_q^n) \\
[12pt] \displaystyle \sum_{k=1}^n x_k \ten \delta_k \in
\mathcal{J}_{p',q'}^1 \Big( \ell_{\infty}^n(\mathcal{M}),
\frac{1}{n} \sum_{k=1}^n \mathcal{E_N} \Big) & \!\! \mapsto \!\! &
\displaystyle \sum_{k=1}^n \pi_{free}^k(x_k) \ten \delta_k \in
L_{p'}(\mathcal{A}_{free}; \ell_{q'}^n) \end{array}$$ are
isomorphisms with complemented range and constants independent of
$n$. In particular, considering the restriction to the diagonal
subspaces $x_1 = x_2 = \ldots = x_n$ we obtain isomorphisms
$$\begin{array}{rcl} x \in \hskip2.3pt
\mathcal{K}_{p,q}^n(\mathcal{M}, \mathcal{E_N}) \hskip2.3pt &
\mapsto & \displaystyle \sum_{k=1}^n \pi_{free}^k(x) \ten \delta_k
\in L_p(\mathcal{A}_{free}; \ell_q^n), \\ [12pt] \displaystyle x
\in \mathcal{J}_{p',q'}^n(\mathcal{M}, \mathcal{E_N}) & \mapsto &
\displaystyle \sum_{k=1}^n \pi_{free}^k(x) \ten \delta_k \in
L_p(\mathcal{A}_{free}; \ell_q^n), \end{array}$$ with complemented
range and constants independent of $n$. Moreover, replacing
$(\mathcal{M}, \mathcal{N}, \mathcal{E_N})$ by $(M_m(\mathcal{M}),
M_m, id_{M_m} \otimes \varphi)$, we obtain complete isomorphisms
with completely complemented ranges and constants independent of
$n$ $$\begin{array}{rcl} x \in \hskip2.3pt
\mathcal{K}_{p,q}^n(\mathcal{M}) \hskip2.3pt & \mapsto &
\displaystyle \sum_{k=1}^n \pi_{free}^k(x) \ten \delta_k \in
L_p(\mathcal{R}_{free}; \ell_q^n), \\ [12pt] \displaystyle x \in
\mathcal{J}_{p',q'}^n(\mathcal{M}) & \mapsto & \displaystyle
\sum_{k=1}^n \pi_{free}^k(x) \ten \delta_k \in
L_p(\mathcal{R}_{free}; \ell_q^n). \end{array}$$
$\mathcal{R}_{free}$ stands for the non-amalgamated free product
$(\mathcal{M},\varphi) * (\mathcal{M},\varphi) * \cdots
* (\mathcal{M},\varphi)$.
\end{theorem}

\noindent \textbf{Sketch of the proof.} It clearly suffices to prove
the first assertion. Since the $\mathcal{J}_{pq}^n$-spaces form an
interpolation scale, as observed in Remark \ref{Singpq}, the same holds for $\mathcal{K}_{pq}^n$ by
duality. On the other hand, the spaces $L_p(\mathcal{A}_{free};
\ell_q^n)$ are particular examples of amalgamated or conditional
$L_p$ spaces (according to the value of $p$ and $q$) and hence also
form an interpolation scale. This is well-known for hyperfinite
algebras through Pisier's work, but note that the free product are
in general not hyperfinite. It suffices to show the extremal cases.
The case $p=q$ is trivial, while the argument for the map
$\mathcal{J}_{p',1}^1 \to L_{p'}(\ell_1^n)$ is essentially contained
in \cite{JP}. Indeed, using
$$\Big\| \sum_{k=1}^n x_k \ten \delta_k
\Big\|_{L_p(\mathcal{A}_{free}; \ell_1^n)} = \inf_{x_k = a_k b_k}
\Big\| \big( \sum_{k=1}^n a_k a_k^* \big)^{\frac12} \Big\|_p \,
\Big\| \big( \sum_{k=1}^n b_k^* b_k \big)^{\frac12} \Big\|_p,$$
see e.g. \cite{J1}, the norm estimate for
$\mathcal{J}_{p',1}^1$ boils down to
\begin{align} \label{opopop}
\Big\| \big( \sum_k \pi_{free}^k(x_k) \pi_{free}^k(x_k)^*
\big)^{\frac12} \Big\|_p \lesssim \Big( \sum_{k=1}^n
\|x_k\|_{2p}^{2p} \Big)^{\frac{1}{2p}} + \Big\| \big( \sum_{k=1}^n
\mathcal{E_N}(x_k x_k^*) \big)^{\frac12} \Big\|_{2p}
\end{align}
and its column version. Equation \eqref{opopop}  follows
from the free Rosenthal inequality \cite{JPX}, see \cite[Corollary
5.3]{JP2}. It is important to recall that, in contrast to the free
Rosenthal inequality, we do not require the $x_k$'s to be mean-zero
in \eqref{opopop}. The combination of \eqref{opopop} with
$\mathcal{J}_{p,q}^n(\mathcal{M}, \mathcal{E_N}) \simeq
\mathcal{R}_{2p,q}^n(\mathcal{M}, \mathcal{E_N}) \oM
\mathcal{C}_{2p,q}^n(\mathcal{M}, \mathcal{E_N})$ gives rise to the
boundedness of $\mathcal{J}_{p',1}^1 \to L_{p'}(\ell_1^n)$. By
interpolation, we deduce that $\mathcal{J}_{p',q'}^1 \to L_{p'}(\ell_{q'}^n)$ is bounded. By
duality, see the argument in \cite{JP}, the boundedness of the
inverse and the fact that the range is complemented follows from
the boundedness of the map $\mathcal{K}_{pq}^1 \to L_p(\ell_q^n)$.
This is even easier and follows from algebraic considerations that
can be found in \cite[Proposition 3.5]{JP}. The proof is complete.
\fin

\begin{remark} \label{RemKfree2}
\emph{Several remarks are in order:}
\begin{itemize}
\item[i)] \emph{The first result of this kind appeared in
\cite{JP}, where tensor independence played the role of freeness
and amalgamation was not considered. According to the classical
theory, we found in this case a non-removable singularity as $p
\to \infty$. Then, the duality argument produces a singularity as
$p \to 1$. As a byproduct of our methods, we shall show in this
paper that this singularity is removable. A key point for it is to
observe that, in the free setting considered in Theorem
\ref{Theorem-Sigma-pq}, we only find a (apparently removable)
singularity when $(p,q) \to (1,\infty)$ simultaneously.}

\vskip3pt

\item[ii)] \emph{The variables $\pi_{free}^k(x_k)$ are replaced by
$\pi_{free}^k(x_k,-x_k)$ in the formulation of Theorem
\ref{Theorem-Sigma-pq} given in \cite{JP2}. This was done to create
mean-zero random variables, a necessary condition for the free
Rosenthal inequality \cite{JPX}. In \eqref{opopop} mean-zero random
variables are not required and this simplifies our embedding.}

\vskip3pt

\item[iii)] \emph{The simpler formulation in \cite{JP} allowed us
to take values in an arbitrary operator space $\mathrm{X}$. In the
framework of Theorem \ref{Theorem-Sigma-pq}, this requires some
additional insight that will be analyzed in a forthcoming paper.}
\end{itemize}
\end{remark}

\section{Sums of independent copies in $L_1(\ell_\infty)$}
\label{Section3}

Let $(\mathcal{M}_k)_{k \ge 1}$ be an increasingly independent
family of top-subsymmetric copies of the von Neumann algebra
$\mathcal{M}$ over $\mathcal{N}$. Let us also consider a
$\sigma$-finite von Neumann algebra $\mathcal{R}$ equipped with a
normal faithful state $\phi$. Given $x \in L_1(\mathcal{M}) \ten
\mathcal{R}$, we want to study sums of independent copies in
$L_1(\mathcal{A}; \ell_\infty(\mathcal{R}))$. According to Section
\ref{Section2}, we know that $$\Big\| \sum_{k=1}^n \pi_k(x) \otimes
\delta_k \Big\|_{L_1(\mathcal{A}; \ell_\infty^n(\mathcal{R}))} =
\inf_{\pi_k(x) = \alpha y_k \beta} \|\alpha\|_{L_2(\mathcal{A})}
\Big( \sup_{1 \le k \le n} \|y_k\|_{\mathcal{A} \bar\ten
\mathcal{R}} \Big) \|\beta\|_{L_2(\mathcal{A})}.$$ On the other
hand, we have the following formula for any $z \in \mathcal{M}
\bar\ten L_1(\mathcal{R})$
\begin{equation} \label{amalgform}
\|z\|_{L_{\infty}(\mathcal{M}; L_1(\mathcal{R}))} = \inf \left\{
\Big\| \Big( \summ_k \mathsf{E}_\mathcal{M}(a_k a_k^*)
\Big)^{\frac12} \Big\|_\mathcal{M} \Big\| \Big( \summ_k
\mathsf{E}_\mathcal{M}(b_k^* b_k) \Big)^{\frac12}
\Big\|_\mathcal{M} \right\}
\end{equation}
where the infimum runs over all decompositions of $z$ into a
finite sum $\sum_k a_k b_k$ and $\mathsf{E}_\mathcal{M}:
\mathcal{M} \bar\ten \mathcal{R} \to \mathcal{M}$ is the
conditional expectation $\mathsf{E}_\mathcal{M} = id \ten \phi$.
Indeed, recall that the term on the right is the norm of $z$ in $$L_\infty^r(\mathcal{M} \bar\otimes \mathcal{R}; \mathsf{E}_\mathcal{M}) \otimes_{\mathcal{M} \bar\otimes \mathcal{R}} L_\infty^c(\mathcal{M} \bar\otimes \mathcal{R}; \mathsf{E}_\mathcal{M}).$$ Then, since the following identifications clearly hold
\begin{eqnarray*}
L_\infty^r(\mathcal{M} \bar\otimes \mathcal{R}; \mathsf{E}_\mathcal{M}) & = & L_2^{-1}(\mathcal{M}) L_2(\mathcal{M} \bar\otimes \mathcal{R}) L_\infty^{-1}(\mathcal{M}), \\ L_\infty^c(\mathcal{M} \bar\otimes \mathcal{R}; \mathsf{E}_\mathcal{M}) & = & L_\infty^{-1}(\mathcal{M}) L_2(\mathcal{M} \bar\otimes \mathcal{R}) L_2^{-1}(\mathcal{M}),
\end{eqnarray*}
we deduce that we are talking about the norm of $z$ in $$L_2^{-1}(\mathcal{M}) L_2(\mathcal{M} \bar\otimes \mathcal{R}) L_\infty^{-1}(\mathcal{M}) \otimes_{\mathcal{M} \bar\otimes \mathcal{R}} L_\infty^{-1}(\mathcal{M}) L_2(\mathcal{M} \bar\otimes \mathcal{R}) L_2^{-1}(\mathcal{M}).$$ In these terms, it is clear that \eqref{amalgform} follows from \cite[Proposition 6.9]{JP2} after taking $(\mathcal{M}, \mathcal{N})$ there to be our pair $(\mathcal{M} \bar\otimes \mathcal{R}, \mathcal{M})$, recall one more time that our terminology for conditional $L_p$ spaces is different from the one used in \cite{JP2}. Given $\lambda > 0$, we also define the spaces
\begin{eqnarray*}
\mathcal{R}_{\infty,1}^\lambda (\mathcal{M}, \mathcal{E_N}) & = &
\mathcal{M} \, \cap \, \sqrt{\lambda} \, L_{2}^{-1}(\mathcal{N})
L_{2}(\mathcal{M}) L_\infty^{-1}(\mathcal{N}), \\
\mathcal{C}_{\infty,1}^\lambda (\mathcal{M}, \mathcal{E_N}) & = &
\mathcal{M} \, \cap \, \sqrt{\lambda} \, L_{\infty}^{-1}
(\mathcal{N}) L_{2}(\mathcal{M}) L_{2}^{-1}(\mathcal{N}).
\end{eqnarray*}
According to \cite{JP2}, the norm on these spaces has the form
\begin{eqnarray*}
\|\alpha\|_{\mathcal{R}_{\infty,1}^\lambda(\mathcal{M},
\mathcal{E_N})} & = & \max \Big\{ \|\alpha\|_{\mathcal{M}},
\sqrt{\lambda} \, \big\|
\mathcal{E_N}(\alpha \alpha^*)^{\frac12} \big\|_{\mathcal{N}} \Big\}, \\
[5pt] \|\beta\|_{\hskip2pt
\mathcal{C}_{\infty,1}^\lambda(\mathcal{M}, \mathcal{E_N})}
\hskip1pt  & = & \max \Big\{ \|\beta\|_{\mathcal{M}},
\sqrt{\lambda} \, \big\|\mathcal{E_N}(\beta^* \beta)^{\frac12}
\big\|_{\mathcal{N}} \Big\}.
\end{eqnarray*}
Of course, when $\lambda = n$ we recover the spaces
$\mathcal{R}_{\infty,1}^n (\mathcal{M}, \mathcal{E_N})$ and
$\mathcal{C}_{\infty,1}^n (\mathcal{M}, \mathcal{E_N})$ from
the previous section. Now we want to consider the corresponding
$\mathcal{J}$-space with values in $L_1(\mathcal{R})$
$$\mathcal{J}_{\infty,1}^\lambda \big( \mathcal{M}, \mathcal{E_N};
L_1(\mathcal{R}) \big) = \mathcal{R}_{\infty,1}^\lambda(\mathcal{M}, \mathcal{E_N}) \ten_\mathcal{M} L_\infty \big( \mathcal{M};L_1(\mathcal{R}) \big) \ten_\mathcal{M} \mathcal{C}_{\infty,1}^\lambda(\mathcal{M}, \mathcal{E_N}).$$ Note that both $\mathcal{R}_{\infty,1}^\lambda(\mathcal{M}, \mathcal{E_N})$ and $\mathcal{C}_{\infty,1}^\lambda(\mathcal{M}, \mathcal{E_N})$ coincide algebraically with $\mathcal{M}$ itself. Thus, the norm in $\mathcal{J}_{\infty,1}^\lambda \big( \mathcal{M}, \mathcal{E_N}; L_1(\mathcal{R}) \big)$ of an element $z$ in the dense subspace $\mathcal{M} \otimes L_1(\mathcal{R})$ is given by $$\|z\|_{\mathcal{J}_{\infty,1}^\lambda(\mathcal{M}, \mathcal{E_N}; L_1(\mathcal{R}))} = \inf_{z = \alpha y \beta} \|\alpha\|_{\mathcal{R}_{\infty,1}^\lambda(\mathcal{M}, \mathcal{E_N})} \|y\|_{L_\infty(\mathcal{M};L_1(\mathcal{R}))} \|\beta\|_{\mathcal{C}_{\infty,1}^\lambda(\mathcal{M}, \mathcal{E_N})}.$$ When no confusion can arise, we shall write $\mathcal{J}_{\infty,1}^\lambda(L_1(\mathcal{R}))$.

\begin{lemma} \label{schu}
We have
\begin{eqnarray*}
\|z\|_{\mathcal{J}_{\infty,1}^\lambda(L_1(\mathcal{R}))} & \sim &
\max \Big\{ \|z\|_{L_{\infty}(\mathcal{M}; L_1(\mathcal{R}))}, \\
& & \sqrt{\lambda} \inf_{z = \alpha y} \big\| \mathcal{E_N}(\alpha
\alpha^*)^{\frac12} \big\|_{\mathcal{N}}
\|y\|_{L_{\infty}(\mathcal{M}; L_1(\mathcal{R}))}, \\
&  & \sqrt{\lambda} \inf_{z = y \beta} \hskip1pt
\|y\|_{L_{\infty}(\mathcal{M}; L_1(\mathcal{R}))} \big \|
\mathcal{E_N}(\beta^* \beta)^{\frac12} \big\|_{\mathcal{N}}, \\
& & \hskip5pt \lambda \hskip1pt \inf_{z = \alpha y \beta} \big\|
\mathcal{E_N}(\alpha \alpha^*)^{\frac12} \big\|_{\mathcal{N}}
\|y\|_{L_{\infty}(\mathcal{M}; L_1(\mathcal{R}))} \big\|
\mathcal{E_N}(\beta^*\beta)^{\frac12} \big\|_{\mathcal{N}} \Big\}.
\end{eqnarray*}
Moreover, the relevant constants are independent of
$\lambda$.
\end{lemma}

\dem The argument can be found in Chapter 6 of \cite{JP2}, see
Lemma 6.3. \fin

Once we have introduced the key spaces, we are ready for some
preliminary estimates. In what follows, we shall assume that $z$
is an element in
$\mathcal{J}_{\infty,1}^\lambda(L_1(\mathcal{R}))$ and we shall
work with a factorization $z= \alpha y \beta$ with
$$\max \Big\{ \|\alpha\|_{\mathcal{M}}, \sqrt{\lambda} \, \big\|
\mathcal{E_N}(\alpha \alpha^*)^{\frac12} \big\|_{\mathcal{N}}
\Big\} \|y\|_{L_\infty(\mathcal{M}; L_1(\mathcal{R}))} \max \Big\{
\|\beta\|_{\mathcal{M}}, \sqrt{\lambda} \, \big\|
\mathcal{E_N}(\beta^* \beta)^{\frac12} \big\|_{\mathcal{N}}
\Big\}$$ being $\le
\|z\|_{\mathcal{J}_{\infty,1}^\lambda(L_1(\mathcal{R}))} +
\varepsilon$ for small $\varepsilon$. Let $a = \displaystyle
\sqrt{\mathbf{1} - \alpha \alpha^*}, b = \sqrt{\mathbf{1} -
\beta^* \beta}$ and $$A_k = \pi_1(a) \pi_2(a) \cdots \pi_{k-1}(a)
\quad , \quad B_k = \pi_{k-1}(b) \cdots \pi_2(b) \pi_1(b).$$

\begin{lemma} \label{Estimacion1}
We have $$\Big\| \sum_{k=1}^n A_k \pi_k(z) B_k \ten \delta_k
\Big\|_{L_{\infty}(\mathcal{A}; \ell_1^n(L_1(\mathcal{R})))} \le
\|y\|_{L_{\infty}(\mathcal{M}; L_1(\mathcal{R}))}.$$
\end{lemma}

\dem We claim that
\begin{eqnarray*}
\lefteqn{\hskip-25pt \Big\| \sum_{k=1}^n \alpha_k y_k \beta_k \ten
\delta_k \Big\|_{L_{\infty}(\mathcal{A};
\ell_1^n(L_1(\mathcal{R})))}} \\ & \le & \Big\| \big( \sum_{k=1}^n
\alpha_k \alpha_k^* \big)^{\frac12} \Big\|_{\mathcal{A}} \Big(
\sup_{1 \le k \le n} \|y_k\|_{L_{\infty}(\mathcal{A};
L_1(\mathcal{R}))} \Big) \Big\| \big( \sum_{k=1}^n \beta_k^*
\beta_k \big)^{\frac12} \Big\|_{\mathcal{A}}.
\end{eqnarray*}
Indeed, if $(\alpha, \beta) = \big( \sum_k \alpha_k \alpha_k^*,
\sum_k \beta_k^* \beta_k \big)$, we find $v_k, w_k \in \mathcal{A}$ with
\begin{itemize}
\item $\alpha_k = \alpha^{\frac12} v_k$ and $\beta_k = w_k
\beta^{\frac12}$,

\vskip5pt

\item $\displaystyle \max \Big\{ \big\| \summ_k v_k v_k^*
\big\|_\mathcal{A}, \big\| \summ_k w_k^*w_k \big\|_{\mathcal{A}}
\Big\} \le 1$.
\end{itemize}
Since $L_{\infty}(\mathcal{A}; \ell_1^n(L_1(\mathcal{R})))$ is an
$\mathcal{A}$-bimodule, it suffices to show $$\Big\| \sum_{k=1}^n
v_k y_k w_k \ten \delta_k \Big\|_{L_{\infty}(\mathcal{A};
\ell_1^n(L_1(\mathcal{R})))} \le \sup_{1 \le k \le n}
\|y_k\|_{L_{\infty}(\mathcal{A}; L_1(\mathcal{R}))}.$$ Factorize
$y_k = \sum_s a_{ks} b_{ks}$ in such a way that $$\max \left\{
\Big\| \summ_s \mathsf{E}_\mathcal{A} (a_{ks} a_{ks}^*)
\Big\|_\mathcal{A}, \Big\| \summ_s \mathsf{E}_\mathcal{A}
(b_{ks}^* b_{ks}) \Big\|_\mathcal{A}^{\null} \right\} \le
\|y_k\|_{L_\infty(\mathcal{A}; L_1(\mathcal{R}))} + \varepsilon.$$ In
particular, we may factorize $v_k y_k w_k = \sum_s v_k a_{ks}
b_{ks} w_k$ to deduce the estimate
\begin{eqnarray*}
\lefteqn{\hskip-15pt \Big\| \sum_{k=1}^n v_k y_k w_k \ten \delta_k
\Big\|_{L_{\infty}(\mathcal{A}; \ell_1^n(L_1(\mathcal{R})))}} \\ &
\le & \Big\| \summ_{k,s} \mathsf{E}_\mathcal{A} \big( v_k a_{ks}
a_{ks}^* v_k^* \big) \Big\|_{\mathcal{A}}^{\frac12} \Big\|
\summ_{k,s} \mathsf{E}_\mathcal{A} \big( w_k^* b_{ks}^* b_{ks} w_k
\big) \Big\|_{\mathcal{A}}^{\frac12} \\ & \le & \Big\| \summ_k v_k
v_k^* \Big\|_\mathcal{A}^{\frac12} \hskip2pt \Big( \sup_{1 \le k
\le n} \|y_k\|_{L_\infty(\mathcal{A}; L_1(\mathcal{R}))} + \varepsilon \Big)
\hskip2pt \Big\| \sum_k w_k^* w_k \Big\|_\mathcal{A}^{\frac12} \\
& \le & \sup_{1 \le k \le n} \|y_k\|_{L_\infty(\mathcal{A};
L_1(\mathcal{R}))} + \varepsilon.
\end{eqnarray*}
This proves our claim if we let $\varepsilon \to 0^+$. Applying it for $$(\alpha_k, y_k, \beta_k)
= \big( A_k \pi_k(\alpha), \pi_k(y), \pi_k(\beta) B_k \big),$$
gives rise to
\begin{eqnarray*}
\lefteqn{\hskip-10pt \Big\| \sum_{k=1}^n A_k \pi_k(z) B_k \ten
\delta_k \Big\|_{L_{\infty}(\mathcal{A};
\ell_1^n(L_1(\mathcal{R})))}} \\ & \le & \Big\| \big( \sum_{k=1}^n
A_k \pi_k(\alpha \alpha^*) A_k^* \big)^{\frac12}
\Big\|_{\mathcal{A}} \|y\|_{L_{\infty}(\mathcal{M};
L_1(\mathcal{R}))} \Big\| \big( \sum_{k=1}^n B_k^* \pi_k(\beta^*
\beta) B_k \big)^{\frac12} \Big\|_{\mathcal{A}}.
\end{eqnarray*}
Therefore, the assertion follows from $$\Big\| \summ_k A_k
\pi_k(\alpha \alpha^*) A_k^* \Big\|_\mathcal{A} \le 1 \quad
\mbox{and} \quad \Big\| \sum_k B_k^* \pi_k(\beta^* \beta) B_k
\Big\|_\mathcal{A} \le 1.$$ Indeed, these estimates are implicit
in \cite[Lemma 7.6]{J3}. The proof is complete. \fin

The next lemma requires a further property of the space
$L_{\infty}(\mathcal{M}; L_1(\mathcal{R}))$ used in the definition
of the amalgamated tensor product
$\mathcal{J}_{\infty,1}^\lambda(L_1(\mathcal{R}))$. Indeed, we
require that for a conditioned subalgebra $\mathcal{N}$ of
$\mathcal{M}$ and a normal $*$-homomorphism $\rho: \mathcal{M} \to
\mathcal{N} \bar{\ten} \mathcal{B}(\ell_2)$, we have
\begin{equation} \label{techasspt}
\big\| \rho \ten id: L_{\infty} \big( \mathcal{M};
L_1(\mathcal{R}) \big) \to L_\infty \big( \mathcal{N} \bar{\ten}
\mathcal{B}(\ell_2); L_1(\mathcal{R}) \big) \big\| \le 1.
\end{equation}
To prove it we use an isometric inclusion
\begin{equation} \label{NDECEq}
L_\infty(\mathcal{M}; L_1(\mathcal{R})) \subset
\mathsf{NDEC}(\mathcal{R}^{op}, \mathcal{M}).
\end{equation}
In other words, by slicing $T_x(r)=(id\ten r)(x)$, we can view
$T_x: \mathcal{R}^{op}\to \mathcal{M}$ as a normal decomposable
map. The norm in the space of decomposable maps (see \cite{JR} for
details) is given by
$$\|T_x\|_{dec}= \inf \left\{ \Big\| \! \left(
\begin{array}{c} S_1 \ \ T_x \\ T_x^* \ S_2 \end{array}
\right) \! \Big\|_{M_2(\mathcal{R}) \to M_2(\mathcal{M})}
\vspace{0.3cm} \mbox{s.t.} \left( \begin{array}{c} S_1 \ \ T_x \\
T_x^* \ S_2 \end{array} \right) \mbox{completely positive}
\right\}.$$ To prove \eqref{NDECEq} let us take $a,b\in L_2(\M)$
and $M_{ab}(y)=ayb$. For $x\in \M\ten L_1(\mathcal{R})$ we see
that $M_{ab}T_x\in \mathsf{NDEC}(\mathcal{R}^{op},L_1(\M))$ and is
of finite rank. However, every completely positive map $T':
\mathcal{R}^{op} \to \M^{op}_*$ defines an element in
$(\mathcal{R}^{op}\ten_{\max}\M^{op})^*$. Given a finite tensor
$z=\sum_j r_j\ten m_j$ of norm less than one, we may lift it to an
element $\hat{z}$ of norm less than one in the unit ball of
$\mathcal{R}^{op}\ten_{\max}\M^{op}$. This implies
 \[ |\langle \hat{z},M_{ab}T_x\rangle| \le
 \|T_x\|_{dec}\|a\|_2\|b\|_2 \, .\]
Since $M_{ab}T_x$ is of finite rank we know that
 \[ |\langle \hat{z},M_{ab}T_x\rangle| =  |\langle z,axb\rangle|
 \, .\]
Then \cite[Proposition 6.9]{JP2} implies that
 \[ \|x\|_{L_{\infty}(\M;L_1(\mathcal{R}))} =
 \sup_{\|a\|_{L_2(\M)}, \, \|b\|_{L_2(\M)} \le 1} \|axb\|_1 \le
 \|T_x\|_{dec}.\]
According to \eqref{amalgform}, the converse follows easily by
factoring $x=x_1x_2$ and using that $x_1x_1^*$ and $x_2^*x_2$
correspond to completely positive maps. Then \eqref{techasspt}
follows from the fact that $\|\rho T_x\|_{dec}\le \|T_x\|_{dec}$.


\vskip3pt

\noindent In the following result, we shall use the conditional
expectations $\mathsf{E}_k: \mathcal{A} \to \mathcal{M}_k$.

\begin{lemma} \label{Estimacion2}
If $\delta \le 1/e$ and $\lambda = \delta^{-1} n$, we have
$$\Big\| \pi_n^{-1} \Big[ \mathsf{E}_n \Big( \sum_{k=1}^n
(\mathbf{1}-A_k)\pi_n(z)(\mathbf{1}-B_k) \Big) \Big]
\Big\|_{\mathcal{J}_{\infty,1}^\lambda(L_1(\mathcal{R}))} \le 2 n
e \delta \,
\|z\|_{\mathcal{J}_{\infty,1}^\lambda(L_1(\mathcal{R}))}.$$
\end{lemma}

\dem By the normal version of Kasparov's dilation theorem
\cite{Pas}, we may assume $\mathsf{E}_n(x) = e_{11} \rho(x)
e_{11}$, where $\rho: \mathcal{A} \to \mathcal{M}_n \bar\ten
\mathcal{B}(\ell_2)$ is a normal $*$-homomorphism. Let us
factorize $z= \alpha y \beta$ as we did before Lemma
\ref{Estimacion1}, with $\lambda = \delta^{-1} n$. Then we get
$$\mathsf{E}_n \Big( \sum_{k=1}^n
(\mathbf{1}-A_k)\pi_n(z)(\mathbf{1}-B_k) \Big) = \mathsf{RDC},$$
where $\mathsf{R,D,C}$ are given by
\begin{eqnarray*}
\mathsf{R} & = & \sum_{k=1}^n e_{11} \rho \big(
(\mathbf{1}-A_k)\pi_n(\alpha) \big) \ten e_{1k}, \\ \mathsf{D} & =
& \sum_{k=1}^n \rho(\pi_n(y)) \ten e_{kk}, \\ \mathsf{C} & = &
\sum_{k=1}^n \rho \big( \pi_n(\beta)(\mathbf{1}-B_k) \big) e_{11}
\ten e_{k1}.
\end{eqnarray*}
It is easily checked that $$\mathsf{RR^*} = \sum_{k=1}^n
\mathsf{E}_n \big( (\mathbf{1}-A_k) \pi_n(\alpha \alpha^*)
(\mathbf{1}-A_k)^* \big) \ten e_{11}.$$ According to the proof of
\cite[Lemma 7.8]{J3} we have
\begin{equation} \label{RR*1}
\big\|\pi_n^{-1} \big( (\mathsf{R} \mathsf{R}^*)^{\frac12} \big)
\big\|_{\mathcal{M}} = \big\| \pi_n^{-1} (\mathsf{R} \mathsf{R}^*)
\big\|_{\mathcal{M}}^{\frac12} \le \sqrt{2 n e \delta} \,
\|\alpha\|_{\mathcal{M}}.
\end{equation}
To estimate $\pi_n^{-1} \big( (\mathsf{R} \mathsf{R}^*)^{\frac12}
\big)$ in $\mathcal{R}_{\infty,1}^\lambda(\mathcal{M},
\mathcal{E_N})$ it remains to control the term $$\sqrt{\delta^{-1}
n} \, \big\| \mathcal{E_N} \big( \pi_n^{-1} (\mathsf{R}
\mathsf{R}^*) \big) \big\|_{\mathcal{N}}^{\frac12} =
\sqrt{\delta^{-1} n} \, \Big\| \sum_{k=1}^n \mathsf{E}_\mathcal{N}
\big( (\mathbf{1}-A_k) \pi_n(\alpha \alpha^*) (\mathbf{1}-A_k)^*
\big) \Big\|_\mathcal{N}^{\frac12}.$$ Finally, applying Lemma 7.1
(ii) and Lemma 7.7 (iv) from \cite{J3}, we obtain
\begin{eqnarray} \label{RR*2}
\lefteqn{\hskip5pt \big\| \mathcal{E_N} \big( \pi_n^{-1}
(\mathsf{R} \mathsf{R}^*) \big) \big\|_{\mathcal{N}}^{\frac12}  =
\Big\| \sum_{k=1}^n \mathsf{E}_\mathcal{N} \big( (\mathbf{1}-A_k)
\mathcal{E_N} (\alpha \alpha^*) (\mathbf{1}-A_k)^* \big)
\Big\|_\mathcal{N}^{\frac12}} \\ \nonumber \null \hskip10pt & \le
& \Big\| \sum_{k=1}^n \mathsf{E}_\mathcal{N} \big(
(\mathbf{1}-A_k) (\mathbf{1}-A_k)^* \big)
\Big\|_\mathcal{N}^{\frac12} \big\| \mathcal{E_N} (\alpha
\alpha^*) \big\|_\mathcal{M}^{\frac12} \le \sqrt{2 n e \delta} \,
\big\| \mathcal{E_N} (\alpha \alpha^*)
\big\|_\mathcal{M}^{\frac12}.
\end{eqnarray}
The combination of \eqref{RR*1} and \eqref{RR*2} (as well as a
symmetric argument) produces $$\big\| \pi_n^{-1} \big(
(\mathsf{RR^*})^{\frac12} \big)
\big\|_{\mathcal{R}_{\infty,1}^\lambda} \le \sqrt{2ne \delta} \,
\|\alpha\|_{\mathcal{R}_{\infty,1}^\lambda} \ \ \mbox{,} \ \ \big\|
\pi_n^{-1} \big( (\mathsf{C^*C})^{\frac12} \big)
\big\|_{\mathcal{C}_{\infty,1}^\lambda} \hskip3pt \le \sqrt{2ne
\delta} \, \|\beta\|_{\mathcal{C}_{\infty,1}^\lambda}.$$
Since we
have the factorization
$$\pi_n^{-1} \Big[ \mathsf{E}_n \Big( \sum_{k=1}^n
(\mathbf{1}-A_k)\pi_n(z)(\mathbf{1}-B_k) \Big) \Big] = \pi_n^{-1}
\big( (\mathsf{RR^*})^{\frac12} \big) \pi_n^{-1} \big( u
\mathsf{D} v \big) \pi_n^{-1} \big( (\mathsf{C^*C})^{\frac12}
\big)$$ for certain contractions $u,v$,  and
$L_\infty(M_n(\mathcal{M}); L_1(\mathcal{R}))$ is an
$\mathcal{M}$-bimodule, we deduce
\begin{eqnarray*}
\lefteqn{\hskip-25pt \Big\| \pi_n^{-1} \Big[ \mathsf{E}_n \Big(
\sum_{k=1}^n (\mathbf{1}-A_k)\pi_n(z)(\mathbf{1}-B_k) \Big) \Big]
\Big\|_{\mathcal{J}_{\infty,1}^\lambda(L_1(\mathcal{R}))}} \\ &
\le & 2ne \delta \, \|\alpha\|_{\mathcal{R}_{\infty,1}^\lambda}
\|\beta\|_{\mathcal{C}_{\infty,1}^\lambda} \big\|
\pi_n^{-1}(\mathsf{D}) \big\|_{L_\infty(M_n(\mathcal{M});
L_1(\mathcal{R}))} \\ & \le & 2ne \delta \,
\|\alpha\|_{\mathcal{R}_{\infty,1}^\lambda}
\|\beta\|_{\mathcal{C}_{\infty,1}^\lambda} \big\| \rho(\pi_n(y))
\big\|_{L_\infty(\mathcal{M}; L_1(\mathcal{R}))} \\ & \le & 2ne
\delta \, \Big( \|\alpha\|_{\mathcal{R}_{\infty,1}^\lambda}
\|y\|_{L_\infty(\mathcal{M}; L_1(\mathcal{R}))}
\|\beta\|_{\mathcal{C}_{\infty,1}^\lambda} \Big) \\ & \le & 2ne
\delta \, \Big(
\|z\|_{\mathcal{J}_{\infty,1}^\lambda(L_1(\mathcal{R}))} +
\varepsilon \Big).
\end{eqnarray*}
We have applied inequality \eqref{techasspt}. The assertion
follows by letting $\varepsilon \to 0$ above. \fin

\begin{lemma} \label{key}
If $\delta \le 1/e$ and $\lambda = \delta^{-1} n$, we have
$$\Big\| \sum_{k=1}^n \pi_k(x) \ten
\delta_k \Big\|_{L_1(\mathcal{A}; \ell_{\infty}^n(\mathcal{R}))}
\ge n (1 - 4e \delta) \,
\|x\|_{\mathcal{J}_{\infty,1}^\lambda(L_1(\mathcal{R}))^*}.$$
\end{lemma}

\dem By homogeneity, we will assume that
$$\|x\|_{\mathcal{J}_{\infty,1}^\lambda(L_1(\mathcal{R}))^*} =
1.$$ Let $z \in \mathcal{A} \otimes L_1(\mathcal{R})$ be a norm
$1$ element of $\mathcal{J}_{\infty,1}^\lambda(L_1(\mathcal{R}))$
such that $|\langle x , z \rangle| = 1-\gamma$ and factorize $z = \alpha y \beta$ with
$\|\alpha\|_{\mathcal{R}_{\infty,1}^\lambda} =
\|y\|_{L_\infty(\mathcal{M};L_1(\mathcal{R}))} =
\|\beta\|_{\mathcal{C}_{\infty,1}^\lambda} \le 1+\gamma.$ First we observe
from Lemma \ref{Estimacion1} that $$\Big| \sum_{k=1}^n \big\langle
\pi_k(x), A_k \pi_k(z) B_k \big\rangle \Big| \le (1+\gamma) \Big\|
\sum_{k=1}^n \pi_k(x) \ten \delta_k \Big\|_{L_1(\mathcal{A};
\ell_{\infty}^n(\mathcal{R}))}.$$ Now, to work through the error
estimate, we use $$z-azb =
z(\mathbf{1}-b)+(\mathbf{1}-a)z-(\mathbf{1}-a)z(\mathbf{1}-b).$$
Hence
\begin{eqnarray*}
n(1-\gamma) & \le & \Big| \sum_{k=1}^n \langle
\pi_k(x), \pi_k(z)\rangle \Big| \\
& \le & \Big| \sum_{k=1}^n \big\langle \pi_k(x), A_k \pi_k(z) B_k
\big\rangle \Big| + \Big| \sum_{k=1}^n \big\langle \pi_k(x),
\pi_k(z) - A_k \pi_k(z) B_k \big\rangle \Big| \\ & \le & (1+\gamma) \Big\|
\sum_{k=1}^n \pi_k(x) \ten \delta_k \Big\|_{L_1(\mathcal{A};
\ell_{\infty}^n(\mathcal{R}))}
\\ & + & \Big| \sum_{k=1} \big\langle \pi_k(x),
(\mathbf{1}-A_k) \pi_k(z) (\mathbf{1}-B_k) \big\rangle \Big| \\ &
+ & \Big| \sum_{k=1}^n \big\langle \pi_k(x), (\mathbf{1}-A_k)
\pi_k(z) \big\rangle \Big| + \Big| \sum_{k=1} \big\langle
\pi_k(x), \pi_k(z)(\mathbf{1}-B_k) \big\rangle \Big|.
\end{eqnarray*}
By top-subsymmetry and \cite[Lemma 7.1]{J3}, we deduce
\begin{align*}
&\Big| \sum_{k=1}^n \big\langle \pi_k(x), (\mathbf{1}-A_k) \,
\pi_k(z) \big\rangle \Big| =  \Big| \big\langle
\pi_n(x), \sum_{k=1}^n (\mathbf{1}-A_k) \, \pi_n(z) \big\rangle \Big| \\
& =  \Big| \big\langle \pi_n(x), \mathsf{E}_n \Big( \sum_{k=1}^n
(\mathbf{1}-A_k) \Big) \, \pi_n(z) \big\rangle \Big| = \Big|
\big\langle x, \mathsf{E}_\mathcal{N} \Big( \sum_{k=1}^n
(\mathbf{1}-A_k) \Big) z \big\rangle \Big| \\
& \le  \Big\| \mathsf{E}_\mathcal{N} \Big( \sum_{k=1}^n
(\mathbf{1}-A_k) \Big) z
\Big\|_{\mathcal{J}_{\infty,1}^\lambda(L_1(\mathcal{R}))}
\|x\|_{\mathcal{J}_{\infty,1}^\lambda(L_1(\mathcal{R}))^*} \le
\Big\| \mathsf{E}_\mathcal{N} \Big( \sum_{k=1}^n (\mathbf{1}-A_k)
\Big) \Big\|_\mathcal{N}.
\end{align*}
Similarly, we find $$\Big| \sum_{k=1}^n \big\langle \pi_k(x),
\pi_k(z) \, (\mathbf{1}-B_k) \big\rangle \Big| \le \Big\|
\mathsf{E}_\mathcal{N} \Big( \sum_{k=1}^n (\mathbf{1}-B_k) \Big)
\Big\|_\mathcal{N}.$$ We refer to \cite[Lemma 7.7 (iii)]{J3} for
$$\max \left\{ \Big\| \mathsf{E}_\mathcal{N} \Big( \sum_{k=1}^n
(\mathbf{1}-A_k) \Big) \Big\|_\mathcal{N}, \Big\|
\mathsf{E}_\mathcal{N} \Big( \sum_{k=1}^n (\mathbf{1}-B_k) \Big)
\Big\|_\mathcal{N}\right\} \le ne \delta.$$ Here we are using
implicitly that we have $$\max \left\{ \big\| \mathcal{E_N}(\alpha
\alpha^*)^{\frac12} \big\|_\mathcal{N}, \big\| \mathcal{E_N}(\beta^*
\beta)^{\frac12} \big\|_\mathcal{N} \right\} \le
\frac{1}{\sqrt{\delta^{-1} n}}.$$ Our argument for the symmetric
term uses Lemma \ref{Estimacion2} instead
\begin{eqnarray*}
\lefteqn{\Big| \sum_{k=1}^n \big\langle \pi_k(x), (\mathbf{1}-A_k)
\, \pi_k(z) \, (\mathbf{1}-B_k) \big\rangle \Big|} \\ & = & \Big|
\big\langle x, \pi_n^{-1} \Big[ \mathsf{E}_n \Big( \sum_{k=1}^n
(\mathbf{1}-A_k) \, \pi_n(z) \, (\mathbf{1}-B_k) \Big) \Big]
\big\rangle \Big| \le 2ne \delta \,
\|x\|_{\mathcal{J}_{\infty,1}^\lambda(L_1(\mathcal{R}))^*}.
 \end{eqnarray*}
This yields $$n(1-\gamma) \le (1+\gamma) \Big\| \sum_{k=1}^n \pi_k(x) \ten
\delta_k \Big\|_{L_1(\mathcal{A}; \ell_{\infty}^n(\mathcal{R}))} +
4ne \delta.$$ Taking $\gamma \to 0^+$, we deduce the
assertion and the proof is complete. \fin

\begin{remark}
\emph{Apart from our references to \cite{J3} in Section \ref{Section1}, the estimation of the error terms above is the only place in this paper where top-subsymmetry really takes place.}
\end{remark}

\noindent Let us consider the following norms
\begin{eqnarray*}
\|x\|_{L_1^r(\mathcal{M}, \mathcal{E_N}; \mathcal{R})} & = &
\inf_{x=ayb} \|a\|_{L_2(\mathcal{N})} \hskip1.5pt
\|y\|_{\mathcal{M} \bar\ten \mathcal{R}} \|b\|_{L_2(\mathcal{M})},
\\ \|x\|_{L_1^c(\mathcal{M}, \mathcal{E_N}; \mathcal{R})} & = &
\inf_{x=ayb} \|a\|_{L_2(\mathcal{M})} \|y\|_{\mathcal{M} \bar\ten
\mathcal{R}} \|b\|_{L_2(\mathcal{N})} \hskip1.5pt , \\
\|x\|_{L_1^s(\mathcal{M}, \mathcal{E_N}; \mathcal{R})} & = &
\inf_{x=ayb} \|a\|_{L_2(\mathcal{N})} \hskip1.5pt \|y\|_{\mathcal{M}
\bar\ten \mathcal{R}} \|b\|_{L_2(\mathcal{N})} \hskip1.5pt.
\end{eqnarray*}

\begin{lemma} \label{dual}
If $z \in \mathcal{M} \ten L_1(\mathcal{R})$, we have
\begin{eqnarray*}
\|z\|_{L_1^r(\mathcal{M}, \mathcal{E_N}; \mathcal{R})^*} & = &
\hskip3pt \inf_{z = \alpha y} \hskip2pt \big\|
\mathcal{E_N}(\alpha \alpha^*)^{\frac12} \big\|_\mathcal{N}
\|y\|_{L_{\infty}(\mathcal{M}; L_1(\mathcal{R}))}, \\
\|z\|_{L_1^c(\mathcal{M}, \mathcal{E_N}; \mathcal{R})^*} & = &
\hskip3pt \inf_{z = y \beta} \hskip2.5pt
\|y\|_{L_{\infty}(\mathcal{M}; L_1(\mathcal{R}))}
\big\| \mathcal{E_N}(\beta^* \beta)^{\frac12} \big\|_\mathcal{N}, \\
\|z\|_{L_1^s(\mathcal{M}, \mathcal{E_N}; \mathcal{R})^*} & = &
\inf_{z = \alpha y \beta} \big\| \mathcal{E_N}(\alpha
\alpha^*)^{\frac12} \big\|_\mathcal{N}
\|y\|_{L_{\infty}(\mathcal{M}; L_1(\mathcal{R}))} \big\|
\mathcal{E_N}(\beta^* \beta)^{\frac12} \big\|_\mathcal{N}.
\end{eqnarray*}
\end{lemma}

\dem Given $x \in \mathcal{M} \ten \mathcal{R}$, let $x=a_jy_jb_j$
with
\begin{eqnarray*}
\|x\|_{L_1^r(\mathcal{M}, \mathcal{E_N}; \mathcal{R})} & \sim &
\|a_1\|_{L_2(\mathcal{N})} \hskip1.6pt \|y_1\|_{\mathcal{M}
\bar\ten
\mathcal{R}} \hskip1.5pt \|b_1\|_{L_2(\mathcal{M})}, \\
\|x\|_{L_1^c(\mathcal{M}, \mathcal{E_N}; \mathcal{R})} & \sim &
\|a_2\|_{L_2(\mathcal{M})} \|y_2\|_{\mathcal{M} \bar\ten
\mathcal{R}} \hskip1.5pt \|b_2\|_{L_2(\mathcal{N})} \hskip1.5pt , \\
\|x\|_{L_1^s(\mathcal{M}, \mathcal{E_N}; \mathcal{R})} & \sim &
\|a_3\|_{L_2(\mathcal{N})} \hskip1.6pt \|y_3\|_{\mathcal{M} \bar\ten
\mathcal{R}} \hskip1.5pt \|b_3\|_{L_2(\mathcal{N})} \hskip1.5pt .
\end{eqnarray*}
Here $\sim$ means up to $(1+\delta)$ for an arbitrary $\delta>0$.
Then, we have
\begin{eqnarray*}
\langle x,z \rangle & = & \mathrm{tr}_{\mathcal{M} \bar\ten
\mathcal{R}} \big( y_jb_jz^*a_j \big) \le \|y_j\|_{\mathcal{M}
\bar\ten \mathcal{R}} \|a_j^* z b_j^*\|_{L_1(\mathcal{M} \bar\ten
\mathcal{R})} \\ & \le & \|a_j\|_{L_2(\mathsf{A}_j)}
\|y_j\|_{\mathcal{M} \bar\ten \mathcal{R}}
\|b_j\|_{L_2(\mathsf{B}_j)} \, \sup_{\begin{subarray}{c}
\|\alpha_j\|_{L_2(\mathsf{A}_j)} \le 1
\\ \|\beta_j\|_{L_2(\mathsf{B}_j)} \le 1 \end{subarray}} \|\alpha_j z
\beta_j\|_{L_1(\mathcal{M} \bar\ten \mathcal{R})},
\end{eqnarray*}
with respect to anti-linear duality and where $$(\mathsf{A}_1,
\mathsf{B}_1, \mathsf{A}_2, \mathsf{B}_2, \mathsf{A}_3,
\mathsf{B}_3) = (\mathcal{N}, \mathcal{M}, \mathcal{M},
\mathcal{N}, \mathcal{N}, \mathcal{N}).$$ According to this, it is
easily seen that the closure of $\mathcal{M} \ten
L_1(\mathcal{R})$ with respect to the norm (for each $j=1,2,3$)
given by the supremum above embeds isometrically into the dual of
$L_1^\bullet(\mathcal{M}, \mathcal{E_N}; \mathcal{R})$ with
$(\bullet,j) = (r,1), (c,2), (s,3)$. Therefore, it suffices to see
that
\begin{eqnarray}
\label{1} \sup_{\alpha_1, \beta_1} \|\alpha_1 z
\beta_1\|_{L_1(\mathcal{M} \bar\ten \mathcal{R})} & = & \inf_{z =
\alpha y} \hskip2pt \big\| \mathcal{E_N}(\alpha
\alpha^*)^{\frac12} \big\|_\mathcal{N}
\|y\|_{L_{\infty}(\mathcal{M}; L_1(\mathcal{R}))}, \\ \label{2}
\sup_{\alpha_2, \beta_2} \|\alpha_2 z \beta_2\|_{L_1(\mathcal{M}
\bar\ten \mathcal{R})} & = & \inf_{z = y \beta} \hskip2.5pt
\|y\|_{L_{\infty}(\mathcal{M}; L_1(\mathcal{R}))} \big\|
\mathcal{E_N}(\beta^* \beta)^{\frac12} \big\|_\mathcal{N},
\end{eqnarray}
where $\alpha_1, \beta_2 \in \mathsf{B}_{L_2(\mathcal{N})}$ and
$\alpha_2, \beta_1 \in \mathsf{B}_{L_2(\mathcal{M})}$; as well as
\begin{eqnarray}\label{3}
\lefteqn{\hskip-15pt \sup_{\alpha_3, \beta_3} \|\alpha_3 z
\beta_3\|_{L_1(\mathcal{M} \bar\ten \mathcal{R})}} \\ \nonumber &
= & \inf_{z = \alpha y \beta} \big\| \mathcal{E_N}(\alpha
\alpha^*)^{\frac12} \big\|_\mathcal{N}
\|y\|_{L_{\infty}(\mathcal{M}; L_1(\mathcal{R}))} \big\|
\mathcal{E_N}(\beta^* \beta)^{\frac12} \big\|_\mathcal{N}
\end{eqnarray}
with $\alpha_3, \beta_3 \in \mathsf{B}_{L_2(\mathcal{N})}$. Since
the proof of \eqref{1} and \eqref{2} is quite similar to that of
\eqref{3}, we shall only give a detailed argument for the last
one. Given a factorization $z= \alpha y \beta$ and $\alpha_0,
\beta_0 \in L_2(\mathcal{N})$, the upper estimate follows from
\begin{eqnarray*}
\lefteqn{\hskip-10pt \|\alpha_0 z \beta_0\|_{L_1(\mathcal{M}
\bar\ten \mathcal{R})}} \\ & \le & \|\alpha_0
\alpha\|_{L_2(\mathcal{M})} \|y\|_{L_\infty(\mathcal{M};
L_1(\mathcal{R}))} \|\beta \beta_0\|_{L_2(\mathcal{M})} \\ & \le &
\|\alpha_0\|_{L_2(\mathcal{N})} \big\|\mathcal{E_N}(\alpha
\alpha^*)^{\frac12} \big\|_{\mathcal{N}}
\|y\|_{L_\infty(\mathcal{M}; L_1(\mathcal{R}))} \big\|
\mathcal{E_N}(\beta^* \beta)^{\frac12} \big\|_\mathcal{N}
\|\beta_0\|_{L_2(\mathcal{N})}.
\end{eqnarray*}
For the lower estimate, we set $$||| z ||| = \inf_{z=\alpha y
\beta} \big\|\mathcal{E_N}(\alpha \alpha^*)^{\frac12}
\big\|_\mathcal{N} \|y\|_{L_{\infty}(\mathcal{M};
L_1(\mathcal{R}))} \big\|\mathcal{E_N} (\beta^* \beta)^{\frac12}
\big\|_\mathcal{N}.$$ This expression defines a norm. Indeed, the
positive definiteness follows from $$|||z||| \ge
\|z\|_{L_\infty^r(\mathcal{M}; \mathcal{E}_\mathcal{N})
\otimes_\mathcal{M} L_\infty(\mathcal{M};L_1(\mathcal{R}))
\otimes_\mathcal{M} L_\infty^c(\mathcal{M};
\mathcal{E}_\mathcal{N})},$$ while the triangle inequality can be
proved following Pisier's factorization argument in \cite[Lemma
3.5]{P2}. Given $z_0 \in \mathcal{M} \ten L_1(\mathcal{R})$, let
us consider a norm $1$ linear functional
$$\phi_{z_0}: \big( \mathcal{M} \ten L_1(\mathcal{R}), ||| \ ||| \big) \to
\C \quad \mbox{such that} \quad \phi_{z_0}(z_0) = |||z_0|||.$$
Note that we have $$\big| \phi_{z_0}(\alpha y \beta) \big| \le
\big\| \mathcal{E_N}(\alpha \alpha^*) \big\|_\mathcal{N}^{\frac12}
\|y\|_{L_\infty(\mathcal{M}; L_1(\mathcal{R}))} \big\|
\mathcal{E_N}(\beta^* \beta) \big\|_\mathcal{N}^{\frac12}.$$ In
particular, we may apply ---as in \cite[Theorem 3.16 + Proposition
6.9]{JP2}--- a standard Grothendieck-Pietsch separation argument
to find states $\varphi_1$ and $\varphi_2$ in $\mathcal{N}^*$ with
associated densities $d_1, d_2$ in $L_1(\mathcal{N}^{**})$, so
that
\begin{eqnarray} \label{bidualphi}
\big| \phi_{z_0}(\alpha y \beta) \big| & \le &
\varphi_1(\mathcal{E_N}(\alpha \alpha^*))^{\frac12}
\|y\|_{L_{\infty}(\mathcal{M};
L_1(\mathcal{R}))} \varphi_2(\mathcal{E_N}(\beta^* \beta))^{\frac12} \\
\nonumber & = & \|d_1^{\frac12} \alpha \|_{L_2(\mathcal{M}^{**})}
\|y\|_{L_{\infty}(\mathcal{M}; L_1(\mathcal{R}))} \|\beta
d_2^{\frac12}\|_{L_2(\mathcal{M}^{**})}.
\end{eqnarray}
We want to construct a norm one functional $\psi: L_1(\mathcal{M}^{**}
\bar{\otimes} \mathcal{R}) \to \C$ with $$\phi_{z_0}(\alpha y
\beta) = \psi \big( d_1^\frac12 \alpha y \beta d_2^\frac12
\big).$$ Let $e_j = \mathrm{supp} \, d_j$ be the support of $d_j$
for $j=1,2$. We know that the space $L_1(\mathcal{M}^{**}
\bar{\ten}\mathcal{R})=L_1(\mathcal{M}^{**})
\hat{\ten}L_1(\mathcal{R})$ is given by the operator space tensor
product. Therefore elements of the form  $$\xi = \sum_{i,j,k,l=1}^n a_{ik}b_{lj} \,  x_{kl}
 \ten y_{ij}$$ with $$\|a\|_2 \big\| (x_{kl})
\big\|_{M_n(L_1(\mathcal{M}^{**}))} \big\| (y_{kl})
\big\|_{M_n(L_1(\mathcal{R}))} \|b\|_2 \le 1$$ are dense in the
unit ball of $L_1(\mathcal{M}^{**}) \hat{\ten}L_1(\mathcal{R})$.
Note also that $\eta = (\sum_{kl} a_{ik} x_{kl} b_{lj})_{ij}$ is
of norm $\le 1$ in $S_1^n(L_1(\mathcal{M}^{**})) = Dec \,
(M_n,L_1(\mathcal{M}^{**}))$ where decomposable refers to linear
combination of positive elements. Thus we can find $h_1, h_2 \in
L_2(\mathcal{M}^{**})$ and $u: M_n\to \mathcal{M}^{**}$ in the
unit ball of $Dec \, (M_n,\mathcal{M}^{**})$ such that
$$\summ_{kl} a_{ik}x_{kl}b_{lj} = h_1u(e_{ij})h_2.$$ Recall that
$Dec \, (M_n,\mathcal{M}^{**}) = Dec \, (M_n,\mathcal{M})^{**}$
and therefore  we can find a net of maps $u_{s}$ in the unit ball
of $Dec \, (M_n,\mathcal{M})$ such that $h_1 u(e_{ij}) h_2 =
\lim_s h_1 u_{s}(e_{ij}) h_2$. Passing to convex combinations, we
may assume that $u_{s}(e_{ij})$ converges in the strong and
strong$^*$ topologies, so that $h_1 u_{s} h_2$ converges to $\eta$
in norm. If we assume additionally that $\xi = e_1 \xi e_2$, we
may replace $h_1$ and $h_2$ by $e_1h_1$ and $h_2e_2$. According to
Kaplansky's density theorem and the norm density of $\sqrt{d_1}
\mathcal{M}^{**}$ in $e_1L_2(\mathcal{M}^{**})$, we see that
$\sqrt{d_1} \mathcal{M}$ is norm dense in $e_1
L_2(\mathcal{M}^{**})$. Similarly, $\mathcal{M} \sqrt{d_2}$ is
norm dense in $L_2(\M^{**}) e_2$. Thus we can find $m_{t_1}$,
$\tilde{m}_{t_2} \in \M$ such that $$e_1 h_1 = \limm_{t_1}
d_1^\frac12 m_{t_1} \quad \mbox{and} \quad h_2 e_2 = \limm_{t_2}
\tilde{m}_{t_2} d_2^\frac12.$$ This shows that
$$\xi =  \limm_{s} \limm_{t_1,t_2} \sum_{i,j=1}^n d_1^\frac12
m_{t_1} u_{s}(e_{ij}) \tilde{m}_{t_2} d_2^\frac12 \ten y_{ij}.$$
We deduce from \eqref{bidualphi} that for fixed $s,t_1,t_2$
\begin{eqnarray*}
\lefteqn{\hskip-15pt \Big| \phi_{z_0} \big( \sum_{i,j=1}^n m_{t_1} u_{s}(e_{ij}) \tilde{m}_{t_2} \ten
y_{ij} \big) \Big|} \\ & \le & \big\| d_1^\frac12 m_{t_1} \big\|_2 \Big\| \sum_{i,j=1}^n u_{s}(e_{ij}) \ten y_{ij} \Big\|_{L_{\infty}(\mathcal{\M}; L_1(\mathcal{R}))} \big\| \tilde{m}_{t_2} d_2^\frac12 \big\|_2
\end{eqnarray*}
Recall that $(y_{ij})\in M_n(L_1(\mathcal{R}))$ has norm $\le 1$.
Since $u_{s}$ is decomposable we see that $u_{s} \ten id: M_n(L_1(\mathcal{R})) \to L_{\infty}(\mathcal{M};L_1(\mathcal{R}))$ is a contraction, which is easy to check for completely positive $u_{s}$. Thus we get
$$\Big\| \sum_{i,j=1}^n u_{s}(e_{ij}) \ten y_{ij} \Big\|_{L_{\infty}(\mathcal{\M} \ten \mathcal{R})}
\le \|u_{s}\|_{dec} \big\| (y_{ij})
\big\|_{M_n(L_1(\mathcal{R}))},$$ and therefore $$\psi(\xi) =
\limm_{s} \limm_{t_1,t_2} \Big| \psi_{z_0} \big( \sum_{i,j=1}^n
m_{t_1}u_{s}(e_{ij}) \tilde{m}_{t_2} \ten y_{ij}) \Big| \le  1.$$
Let us resume what we have proved so far. For fixed $n\in
\mathbb{N}$ we have shown that $\psi(\sqrt{d_1} m \sqrt{d_2} \ten
y ) = \phi_{z_0}(m\ten y)$ extends to a continuous functional on
the Banach space projective tensor product
$e_1L_1(\mathcal{M}^{**})e_2\ten_{\pi}L_1(\mathcal{R})$ such that
 \[  \Big| \psi \big( \sum_{ij=1}^n x_{ij} \ten y_{ij} \big) \Big| \le
 \big\| (x_{ij}) \big\|_{S_1^n(L_1(\mathcal{M}^{**}))}
 \big\| (y_{ij}) \big\|_{M_n(L_1(\mathcal{R}))}  \, . \]
Since left and right multiplications with $e_1$, $e_2$ are
completely contractive, we may extend $\psi$ to
$L_1(\mathcal{M}^{**})\ten_{\pi}L_1(\mathcal{R})$ satisfying the
same inequality. This means $\psi$ induces a linear map
$T_{\psi}:L_1(\mathcal{R})\to
L_1(\mathcal{M}^{**})^*=\mathcal{M}^{op**}$ such that
 \[ \big\| id_{M_n}\ten T_{\psi}:M_n(\mathcal{R})\to
 M_n(\mathcal{M}^{op**}) \big\|\le 1  \, . \]
Since this is true for all $n\in \mathbb{N}$ we deduce that
$T_{\psi}$ is completely bounded. According to Effors/Ruan's
theorem \cite[Theorem 7.2.4]{ER}
$\mathcal{CB}(L_1(\mathcal{R}),\mathcal{M}^{op**}) =
\mathcal{R}^{op}\bar{\ten} \mathcal{M}^{op**} =
(L_1(\mathcal{R}\bar{\ten}\mathcal{M}^{**}))^*$. Therefore $\psi$
corresponds to a norm one functional on
$L_1(\mathcal{R}\bar{\ten}\mathcal{M}^{**})$ such that
 \[ \psi \big( d_1^\frac12 x d_2^\frac12 \ten y \big) =
 \phi_{z_0}(x\ten y) \, . \]
Now we have to replace $d_1, d_2 \in L_1(\mathcal{N}^{**})$ using
an ultraproduct procedure. We recall from \cite[Section 6.2]{JP2}
that we have a completely positive, completely isometric
$\M$-bimodule map $$\rho: L_1(\mathcal{M}^{**})\to \prodd_{\U}
L_1(\mathcal{M}),$$ such that $\rho^*: (\prod_{\U}
L_1(\mathcal{M}))^*\to (\mathcal{M}^{op})^{**}$ is a conditional
expectation. Thus $$\rho \ten id: L_1(\mathcal{M}^{**}) \hat{\ten}
L_1(\mathcal{R}) \to \prodd_\U L_1(\mathcal{M}) \hat{\ten}
L_1(\mathcal{R}).$$ Therefore we find a norm one functional
$\psi': \prod_\U L_1(\mathcal{M}) \hat{\ten}L_1(\mathcal{R})\to
\C$ such that $\psi' \circ (\rho \ten id) = \psi$. The map $\rho$
also induces a map $\rho_p:L_p(\mathcal{M}^{**})\to \prod_\U
L_p(\mathcal{M})$ which remains a $\M$-bimodule map. In
particular, we get $$\rho \big( d_1^\frac12 x d_2^\frac12 \big) =
\rho_2 \big( d_1^\frac12 \big) \hskip1pt x \hskip1pt \rho_2 \big(
d_2^\frac12 \big) \quad \mbox{for} \quad x \in \mathcal{M}.$$ Let
us recall that the inclusion $L_2(\mathcal{N}^{**}) \subset
L_2(\mathcal{M}^{**})$ is defined with the help of the conditional
expectation $\mathcal{E_N}: \mathcal{M} \to \mathcal{N}$, more
precisely $\mathcal{E_N}^{**}$ which is still a (maybe
non-faithful) conditional expectation. We recall from \cite[Lemma
6.2 ii)]{JP2} that $\rho(L_1(\mathcal{N}^{**})) \subset \prod_\U
L_1(\mathcal{N})$ and hence $\rho_2(\sqrt{d_j}) \in \prod_\U
L_2(\mathcal{N})$. Therefore we find
\begin{eqnarray*}
|||z_0||| & = & \big| \psi' \circ (\rho \otimes id) \big( d_1^\frac12 z_0 d_2^\frac12 \big) \big| \\
& \le & \big\| \rho_2 \big( d_1^\frac12 \big) z_0 \rho_2
\big(d_2^\frac12 \big) \big\|_{\prod_\U
L_1(\mathcal{M}\bar{\ten}\mathcal{R})} \\ & = & \limm_{i,\U}
\big\| \rho_2 \big( d_1^\frac12 \big)_i z_0 \rho_2
\big(d_2^\frac12 \big)_i
\big\|_{L_1(\mathcal{M}\bar{\ten}\mathcal{R})} \\ & \le &
\sup_{\|\alpha\|_{L_2(\NN)}, \|\beta\|_{L_2(\NN)} \le 1} \big\|
\alpha z_0 \beta \big\|_{L_1(\mathcal{M} \bar{\ten}\mathcal{R})}.
\end{eqnarray*}
This concludes the proof of \eqref{3}. The argument for \eqref{1} and \eqref{2} is similar. \fin

\noindent Let us define the space
\begin{eqnarray*}
\mathcal{K}_{1,\infty}^\lambda (\mathcal{M}, \mathcal{E_N};
\mathcal{R}) & = & \lambda \, L_1(\mathcal{M}; \mathcal{R}) \\ & + &
L_1^s(\mathcal{M}, \mathcal{E_N}; \mathcal{R}) \\ & + &
\sqrt{\lambda} \, L_1^r(\mathcal{M}, \mathcal{E_N}; \mathcal{R}) \\
& + & \sqrt{\lambda} \, L_1^c(\mathcal{M}, \mathcal{E_N};
\mathcal{R}),
\end{eqnarray*}
where $L_1(\mathcal{M}; \mathcal{R})$ is a shortened way of writing $L_1(\mathcal{M}; L_\infty(\mathcal{R}))$.
We shall often write $\mathcal{K}_{1,\infty}^\lambda(\mathcal{R})$. The norm of $x \in \mathcal{K}_{1,\infty}^\lambda(\mathcal{R})$ is given by $$\inf_{x = \sum_{1}^4 x_j} \lambda \hskip1pt
\|x_1\|_{L_1(\mathcal{M}; \mathcal{R})} + \sqrt{\lambda} \hskip1pt
\|x_2\|_{L_1^r(\mathcal{M}, \mathcal{E_N}; \mathcal{R})} +
\sqrt{\lambda} \hskip1pt \|x_3\|_{L_1^c(\mathcal{M}, \mathcal{E_N};
\mathcal{R})} + \|x_4\|_{L_1^s(\mathcal{M}, \mathcal{E_N};
\mathcal{R})}.$$ The following result probably holds in larger generality. However,
this requires additional fine tuning on the assumptions. For our
purpose, finite-dimensional $\mathcal{R}$'s are enough.

\begin{theorem} \label{4-term}
Let us consider a conditioned subalgebra $\mathcal{N}$ of
$\mathcal{M}$ and a finite dimensional von Neumann algebra
$\mathcal{R}$. Let $(\mathcal{M}_k)_{k \ge 1}$ be an increasingly
independent family of top-subsymmetric copies of $\mathcal{M}$
over $\mathcal{N}$. Then, the following estimate holds up to an
absolute constant for any $n \ge 1$ $$\Big\| \sum_{k=1}^n \pi_k(x)
\ten \delta_k \Big\|_{L_1(\mathcal{A};
\ell_{\infty}^n(\mathcal{R}))} \sim
\|x\|_{\mathcal{K}_{1,\infty}^n(\mathcal{R})}.$$
\end{theorem}

\dem By the triangle inequality
\begin{equation} \label{indest1}
\Big\| \sum_{k=1}^n \pi_k(x) \ten \delta_k
\Big\|_{L_1(\mathcal{A}; \ell_{\infty}^n(\mathcal{R}))} \le n
\|x\|_{L_1(\mathcal{M}; \mathcal{R})}.
\end{equation}
If $x=ayb$ with $a \in L_2(\mathcal{N})$ and $b \in
L_2(\mathcal{M})$, then $$\pi_k(x) = a \pi_k(y) \pi_k(b) \quad
\mbox{and} \quad \sum_{k=1}^n \pi_k(x) \ten \delta_k = a \Big(
\sum_{k=1}^n \pi_k(y) u_k \ten \delta_k \Big) \Big( \sum_{k=1}^n
\pi_k(b^*b) \Big)^{\frac12},$$ where the $u_k$'s are contractions
in $\mathcal{A}$. This immediately gives
\begin{equation} \label{indest2}
\Big\| \sum_{k=1}^n \pi_k(x) \ten \delta_k
\Big\|_{L_1(\mathcal{A}; \ell_{\infty}^n(\mathcal{R}))} \le
\sqrt{n} \, \|a\|_{L_2(\mathcal{N})} \|y\|_{\mathcal{M} \bar\ten
\mathcal{R}} \|b\|_{L_2(\mathcal{M})}.
\end{equation}
In fact, the same argument provides the remaining individual
estimates
\begin{equation}
\label{indest3} \Big\| \sum_{k=1}^n \pi_k(x) \ten \delta_k
\Big\|_{L_1(\mathcal{A}; \ell_{\infty}^n(\mathcal{R}))} \le \min
\Big\{ \sqrt{n} \, \|x\|_{L_1^c(\mathcal{M}, \mathcal{E_N};
\mathcal{R})}, \|x\|_{L_1^s(\mathcal{M}, \mathcal{E_N};
\mathcal{R})} \Big\}.
\end{equation}
The combination of \eqref{indest1}, \eqref{indest2} and
\eqref{indest3} shows that the upper estimate holds contractively.
Let us now prove the lower estimate. Since $\mathcal{R}$ is finite
dimensional we may characterize the dual space of
$\mathcal{K}_{1,\infty}^{\la}$. Indeed, it follows from  Lemma
\ref{schu} and Lemma  \ref{dual} that
\begin{eqnarray*}
\|z\|_{\mathcal{J}_{\infty,1}^\lambda} & \sim & \max \Big\{
\|z\|_{L_{\infty}(\mathcal{M}; L_1(\mathcal{R}))}, \sqrt{\lambda}
\inf_{z = \alpha y} \big\| \mathcal{E_N}(\alpha
\alpha^*)^{\frac12} \big\|_{\mathcal{N}}
\|y\|_{L_{\infty}(\mathcal{M}; L_1(\mathcal{R}))}, \\
&  & \quad  \quad \quad \sqrt{\lambda} \inf_{z = y \beta}
\hskip1pt \|y\|_{L_{\infty}(\mathcal{M}; L_1(\mathcal{R}))} \big
\|
\mathcal{E_N}(\beta^* \beta)^{\frac12} \big\|_{\mathcal{N}}, \\
& & \quad  \quad \quad\hskip5pt \lambda \hskip1pt \inf_{z = \alpha
y \beta} \big\| \mathcal{E_N}(\alpha \alpha^*)^{\frac12}
\big\|_{\mathcal{N}} \|y\|_{L_{\infty}(\mathcal{M};
L_1(\mathcal{R}))} \big\| \mathcal{E_N}(\beta^*\beta)^{\frac12}
\big\|_{\mathcal{N}} \Big\}
\\ & = & \la \hskip1pt \sup \Big\{ |\mbox{tr}(x^*z)|
\, \big| \ \|x\|_{\mathcal{K}_{1,\infty}^{\la}} \le 1 \Big\}.
\end{eqnarray*}
Since the embedding of $\mathcal{K}_{1,\infty}^{\la}$ in its
bidual is isometric, we deduce from Lemma \ref{key}
$$\|x\|_{\mathcal{K}_{1,\infty}^{\la}(\mathcal{R})} \lesssim \la
\hskip1pt \|x\|_{\mathcal{J}_{\infty,1}^{\la}(L_1(\mathcal{R}))^*}
\le 16 \hskip1pt e \hskip1pt \Big\| \sum_{k=1}^n \pi_k(x) \ten
\delta_k \Big\|_{L_1(\mathcal{A};
\ell_{\infty}^n(\mathcal{R}))}.$$ We used $\lambda = n/\delta$ and
$\delta=1/8e$ so that $1-4e\delta=\frac12$. The proof is complete.
\fin

\section{A vector-valued embedding result}
\label{Section4}

Given two von Neumann algebras $\mathcal{M}$ and $\mathcal{R}$ as
in the previous section, our aim now is to find a complete
embedding of $L_p(\mathcal{M}; \mathcal{R})$ for each $1 < p <
\infty$ into an ultraproduct of the form
$$\prodd_{n, \U} L_1 \big( \mathcal{A}_n;
\ell_\infty^{\mathrm{k}_n}(\mathcal{R}) \big).$$ Moreover, for our
applications we also need some additional information on how such
an embedding is constructed in order to maintain the notion of
independent copies. In the following, $X_\mathcal{M}$ will be an
operator space containing $\mathcal{M}$ as a two-sided ideal. Then we may define
$$L_{2p}(\mathcal{M}) X_\mathcal{M} L_{2p}(\mathcal{M}) = L_{2p}^r(\mathcal{M})
\ten_{\mathcal{M},h} X_\mathcal{M} \ten_{\mathcal{M},h}
L_{2p}^c(\mathcal{M}).$$
We will also work with subspaces and
quotients of $$\big( \mathcal{M} \oplus L_2^r(\mathcal{M}) \big)
\ten_{\mathcal{M},h} X_\mathcal{M} \ten_{\mathcal{M},h} \big(
\mathcal{M} \oplus L_2^c(\mathcal{M}) \big).$$ Our main tool is a
standard modification of the so-called Pisier's exercise, see
\cite{JP4,X3} and \cite[Exercise 7.9]{P3}. In other words, a way to
reformulate complex interpolation in this setting. We follow the
same approach as in \cite{JP4,JP2}. Indeed, let $\mathcal{S}$ be the
strip of complex numbers $z$ with $0 \le \mathrm{Re}(z) \le 1$ and
let $\partial_0 \cup \partial_1$ be the partition of its boundary
$\partial \mathcal{S}$ with $\partial_j$ the line of $z$'s with
$\mathrm{Re}(z)=j$. If $0 < \theta < 1$, let $\mu_\theta$ be the
harmonic measure of the point $z = \theta$. This is a probability
measure on the boundary $\partial \mathcal{S}$ (with density given
by the Poisson kernel in the strip) that can be written as
$\mu_\theta = (1 - \theta) \mu_0 + \theta \mu_1$, with $\mu_j$ being
probability measures supported by $\partial_j$ and such that
\begin{equation} \label{Harmonic}
f(\theta) = \int_{\partial \mathcal{S}} f d\mu_\theta
\end{equation}
for any bounded analytic $f$ extended non-tangentially to
$\partial \mathcal{S}$. Let $$\mathcal{S}^r_{\mathcal{M}} = \Big(
L_2^c(\partial_0) \bar\ten \mathcal{M} \Big) \oplus \Big( L_2^r(\partial_1)
\ten_h L_2^r(\mathcal{M}) \Big),$$ $$\mathcal{S}^c_{\mathcal{M}} =
\Big( \mathcal{M} \bar\ten L_2^r(\partial_0) \Big) \oplus \Big( L_2^c(\mathcal{M})
\ten_h L_2^c(\partial_1) \Big).$$ The von Neumann algebra tensor product used above is the weak closure of the minimal tensor product, which in this particular case coincides with the Haagerup tensor product since we have either a column space on the left or a row space on the right. In particular, the only difference is that we are taking the closure in the weak operator topology. The direct sums will be taken Hilbertian. Then, if $\mathcal{M}$ comes equipped with a normal strictly semifinite faithful weight $\psi$ and $d_{\psi}$ denotes the associated density, we define $\mathsf{H}_\theta^r(\mathcal{M})$ as the subspace of all pairs $(f_0,f_1)$ of functions in $\mathcal{S}^r_{\mathcal{M}}$ such that for every scalar-valued analytic function $g: \mathcal{S} \to \C$ (extended non-tangentially to the boundary) with $g(\theta)=0$, we have $$(1-\theta) \int_{\partial_0} g(z) \, d_\psi^{\frac12} f_0(z) \, d\mu_0(z) + \theta \int_{\partial_1} g(z) f_1(z) \, d\mu_1(z) = 0.$$ Similarly, the condition on $\mathsf{H}_\theta^c(\mathcal{M}) \subset \mathcal{S}^c_{\mathcal{M}}$ is $$(1-\theta) \int_{\partial_0} g(z) f_0(z) \, d_\psi^{\frac12} \, d\mu_0(z) + \theta \int_{\partial_1} g(z) f_1(z) \, d\mu_1(z) = 0.$$ We shall also need to consider the subspaces
\begin{eqnarray*}
\mathsf{H}_{r,0} & = & \Big\{ (f_0,f_1) \in
\mathsf{H}_\theta^r(\mathcal{M}) \, \big| \ (1-\theta)
\int_{\partial_0} d_{\psi}^{\frac12} f_0 d\mu_0 + \theta
\int_{\partial_1} f_1 d\mu_1 = 0 \Big\},
\\ \mathsf{H}_{c,0} & = & \Big\{ (f_0,f_1) \in
\mathsf{H}_\theta^c(\mathcal{M}) \, \big| \ (1-\theta)
\int_{\partial_0} f_0 d_{\psi}^{\frac12} d\mu_0 + \theta
\int_{\partial_1} f_1 d\mu_1 = 0 \Big\}.
\end{eqnarray*}
We define the $\mathcal{M}$-bimodules
$$\mathcal{H}_r(\mathcal{M},\theta) =
\mathsf{H}_\theta^r(\mathcal{M}) / \mathsf{H}_{r,0} \quad \mbox{and}
\quad \mathcal{H}_c(\mathcal{M},\theta) =
\mathsf{H}_\theta^c(\mathcal{M}) / \mathsf{H}_{c,0}.$$

\begin{remark}
\emph{We may think of $\mathcal{H}_r(\mathcal{M},\theta)$ as the
space of $\mathcal{M} + L_2^r(\mathcal{M})$-valued analytic
functions $f$ on the strip, with $f(\partial_0) \subset
\mathcal{M}$ and $f(\partial_1) \subset L_2^r(\mathcal{M})$
quotiented by the equivalence relation $f_1 \sim f_2$ iff both
take the same value at $\theta$. A similar observation holds for
$\mathcal{H}_c(\mathcal{M},\theta)$. It is somewhat encoded in the
proof of Proposition \ref{emb} that indeed
$$\mathcal{H}_r(\mathcal{M},\theta) = L_{2p}^r(\mathcal{M}) \qquad
\mbox{and} \qquad \mathcal{H}_c(\mathcal{M},\theta) =
L_{2p}^c(\mathcal{M}).$$}
\end{remark}
\noindent In the following we use the notation
$$\mathcal{H}_\theta(X_\mathcal{M}) =
\mathcal{H}_r(\mathcal{M},\theta) \ten_{\mathcal{M},h}
X_\mathcal{M} \ten_{\mathcal{M},h}
\mathcal{H}_c(\mathcal{M},\theta).$$

\begin{lemma} \label{step1}
Given $1 \le p \le \infty$, we have a contractive inclusion
$$S_{2p'}^m L_{2p}(M_m \ten \mathcal{M}) M_m(X_\mathcal{M})
L_{2p}(M_m \ten \mathcal{M}) S_{2p'}^m \subset R_m \ten_h
\mathcal{H}_{\frac1p}(X_\mathcal{M}) \ten_h C_m.$$
\end{lemma}

\dem We claim that the inclusion $$S_{2p'}^m L_{2p}(M_m \ten
\mathcal{M}) \subset R_m \ten_h \mathcal{H}_r(\mathcal{M},1/p)
\ten_h R_m$$ is contractive.  Let $x=ab$ be such that $a \in
S_{2p'}^m$ and $b \in L_{2p}(M_m \ten \mathcal{M})$ are norm $1$
elements. Using the fact that $R_m \ten_h \mathcal{H}_r(\mathcal{M},1/p)
\ten_h R_m$ is a right $M_m(\mathcal{M})$-module we may apply polar
decomposition and assume that $a$ and $b$ are positive. Indeed,
write $ab = |a^*| u_ab = |a^*| |b^* u_a^*| u_{ab}$ and use that
$$\big\| |a^*| \big\|_{2p'} = \|a\|_{2p'} \qquad \mbox{and} \qquad \big\|
|b^* u_a^*| \big\|_{2p} \le \|b\|_{2p}.$$ Define the analytic function
$f: z \in \mathcal{S} \mapsto a^{(1-z)p'}b^{pz}$ with $f(\frac1p)=x$.
If we set $f_j = f_{\mid_{\partial_j}}$, it is clear that
\begin{eqnarray*}
\lefteqn{\hskip-5pt \|x\|_{R_m \ten_h
\mathcal{H}_r(\mathcal{M},1/p) \ten_h R_m}} \\ [5pt] & \le &
\big\|(f_0,f_1) \big\|_{R_m \ten_h \mathcal{S}^r_\mathcal{M}
\ten_h R_m} \\ & = & \Big( (1-\mbox{$\frac1p$}) \hskip1pt
\|f_0\|_{R_m \ten_h (L_2^c(\partial_0) \bar\ten \mathcal{M})
\ten_h R_m}^2 + \mbox{$\frac1p$} \hskip1pt \|f_1\|_{R_m \ten_h
(L_2^r(\partial_1) \ten_h L_2^r(\mathcal{M})) \ten_h R_m}^2
\Big)^{\frac12}.
\end{eqnarray*}
The space $R_m \ten_h (L_2^c(\partial_0) \bar\ten \mathcal{M})
\ten_h R_m$ is completely isometric to $$(R_m \ten_h R_m)
\ten_{M_m,h} \big( C_m \ten_h (L_2^c(\partial_0) \bar\ten
\mathcal{M}) \ten_h R_m \big) \ten_{M_m,h} (C_m \ten_h R_m),$$
which in turn is isometric to $S_2^m L_\infty (M_m \ten
\mathcal{M}; L_2^c(\partial_0))$. On the other hand, given any $z
\in \partial_0$ we have that $f(z) = a^{p'} u_z$ with $u_z$ being
a unitary in $M_m(\mathcal{M})$ for each $z \in \partial_0$.
Hence, we get $$\|u\|_{L_\infty (M_m \ten \mathcal{M};
L_2^c(\partial_0))} = \Big\| \big( \int_{\partial_0} u_z^* u_z \,
d \mu_0(z) \big)^{\frac12} \Big\|_{M_m \otimes \mathcal{M}} = 1$$
and $\|f_0\|_{R_m \ten_h (L_2^c(\partial_0) \bar\ten \mathcal{M})
\ten_h R_m}^2 \le \|a\|_{2p'}^{2p'} \le 1$. Moreover, it is easy
to check that $$\|f_1\|_{R_m \ten_h (L_2^r(\partial_1) \ten_h
L_2^r(\mathcal{M})) \ten_h R_m}^2 = \int_{\partial_1} \|v_z
b^p\|_{L_2(M_m \ten \mathcal{M})}^2 \, d \mu_1(z) =
\|b\|_{2p}^{2p} \le 1.$$ Putting altogether, we deduce our claim.
Similarly, $$L_{2p}(M_m \ten N) S_{2p'}^m \subset C_m \ten_h
\mathcal{H}_c(\mathcal{M},1/p) \ten_h C_m$$ holds contractively.
Therefore, the assertion follows from
\begin{eqnarray*}
\lefteqn{R_m \ten_h \mathcal{H}_{\frac1p}(X_\mathcal{M}) \ten_h
C_m} \\ & = & \big( R_m \ten_h \mathcal{H}_r(\mathcal{M}, 1/p)
\ten_h R_m \big) \bullet M_m(X_\mathcal{M}) \bullet \big( C_m
\ten_h \mathcal{H}_c(\mathcal{M},1/p) \ten_h C_m \big).
\end{eqnarray*}
where the symbol $\bullet$ stands for the amalgamated tensor
product $\ten_{M_m(\mathcal{M}),h}$. \fin

\begin{proposition} \label{emb}
We have a complete isometry $$L_p(\mathcal{M}; \mathcal{R}) \to
\mathcal{H}_{\frac1p}(\mathcal{M} \bar\ten \mathcal{R}) \quad
\mbox{for each} \quad 1 < p < \infty.$$
\end{proposition}

\dem We have
\begin{eqnarray*}
S_1^m(L_p(\mathcal{M}; \mathcal{R})) & = & S_{2p'}^m L_p(M_m \ten
\mathcal{M}; \mathcal{R}) S_{2p'}^m \\
& = & S_{2p'}^m L_{2p}(M_m \ten \mathcal{M}) M_m(\mathcal{M}
\bar\ten \mathcal{R}) L_{2p}(M_m \ten \mathcal{M}) S_{2p'}^m.
\end{eqnarray*}
Hence Lemma \ref{step1} implies that $$S_1^m(L_p(\mathcal{M};
\mathcal{R})) \subset R_m \ten_h \mathcal{H}_{\frac1p}(\mathcal{M}
\bar\ten \mathcal{R}) \ten_h C_m =
S_1^m(\mathcal{H}_{\frac1p}(\mathcal{M} \bar\ten \mathcal{R}))$$ is
contractive for all $m$. Therefore, the inclusion of
$L_p(\mathcal{M}; \mathcal{R})$ into $\mathcal{H}_{1/p}(\mathcal{M}
\bar\ten \mathcal{R})$ is a complete contraction. To complete the
argument, we proceed by duality and analyze the inclusion
\begin{equation} \label{ContractDual}
S_1^m (L_{p'}(\mathcal{M}; L_1(\mathcal{R}))) \subset
M_m(\mathcal{H}_{\frac1p}(\mathcal{M} \bar\ten \mathcal{R}))^*.
\end{equation}
As in the proof of Lemma \ref{step1}, we may factorize $x= abscd$
with $a,d \in S_{2p}^m$ and $b,c \in L_{2p'}(M_m \ten
\mathcal{M})$ positive norm $1$ elements and $s$ being a not
necessarily positive norm $1$ element in
$L_\infty(M_m(\mathcal{M}); L_1(\mathcal{R}))$. Let us now
consider a norm $1$ element $y \in
M_m(\mathcal{H}_{1/p}(\mathcal{M} \bar\ten \mathcal{R}))$. Then we
may find an analytic function $\xi = \alpha w \beta$ in the
equivalence class determined by $y$ such that

\vskip3pt

\begin{itemize} \item We have $w \in M_k(\mathcal{M}
\bar\ten \mathcal{R})$ for some $k \ge 1$ and
$$\null \hskip27pt \alpha =(\alpha_0, \alpha_1) \in M_{m,k} \big(
L_2^c(\partial_0) \bar\ten \mathcal{M} \big) \oplus M_{m,k} \big(
L_2^r(\partial_1) \ten_h L_2^r(\mathcal{M}) \big),$$ $$\null
\hskip27pt \beta = (\beta_0, \beta_1) \hskip1.5pt \in M_{k,m}
\big( \mathcal{M} \bar\ten L_2^r(\partial_0) \big) \oplus M_{k,m}
\big( L_2^c(\mathcal{M}) \ten_h L_2^c(\partial_1) \big).$$

\vskip3pt

\item The estimate $\|w\|_{M_k(\mathcal{M} \bar\ten \mathcal{R})}
\le 1$ holds and $$\null \hskip25pt \Big( (1 - \mbox{$\frac1p$})
\hskip1pt \|\alpha_0\|_{M_{m,k}(L_2^c(\partial_0) \bar\ten
\mathcal{M})}^2 + \mbox{$\frac1p$} \hskip1pt
\|\alpha_1\|_{M_{m,k}(L_2^r(\partial_1) \ten_h
L_2^r(\mathcal{M}))}^2 \Big)^{\frac12} \le 1,$$ $$\null \hskip25pt
\Big( (1 - \mbox{$\frac1p$}) \hskip1pt
\|\beta_0\|_{M_{k,m}(\mathcal{M} \bar\ten L_2^r(\partial_0))}^2 +
\mbox{$\frac1p$} \hskip1pt \|\beta_1\|_{M_{k,m}(L_2^c(\mathcal{M})
\ten_h L_2^c(\partial_1))}^2 \hskip1pt \Big)^{\frac12} \le 1.$$
\end{itemize}

\vskip3pt

\noindent By adding zeros if necessary, we assume $m=k$ for
simplicity. As in Lemma \ref{step1}, we may define $g(z) =
a^{zp}b^{(1-z)p'}sc^{(1-z)p'}d^{zp}.$ Note that $g$ is also
analytic and hence the identity below holds $$\langle x,y \rangle
= \big\langle g(\mbox{$\frac1p$}), \xi(\mbox{$\frac1p$})
\big\rangle = \int_{\partial \mathcal{S}} \big\langle g(z), \xi(z)
\big\rangle \, d\mu_{\frac1p}(z).$$ Now we claim that
\begin{eqnarray} \label{claimholder}
\Big| \int_{\partial_0} \big\langle g(z), \xi(z) \big\rangle \, d
\mu_0(z) \Big| & \le & \Big\| \big( \int_{\partial_0}
\alpha_0(z)^* \alpha_0(z) \, d\mu_0(z) \big)^{\frac12}
\Big\|_{M_m(\mathcal{M})} \\ \nonumber & \times & \Big\| \big(
\int_{\partial_0} \beta_0(z) \hskip1.5pt \beta_0(z)^* \, d\mu_0(z)
\big)^{\frac12} \Big\|_{M_m(\mathcal{M})}
\end{eqnarray}
Indeed, since $g_{\mid_{\partial_0}}(z) = u_z b^{p'} s c^{p'} v_z$
with $u_z, v_z$ unitaries and $b,c,s,w$ are norm $1$
\begin{eqnarray*}
\lefteqn{\Big| \int_{\partial_0} \big\langle g(z), \xi(z)
\big\rangle \, d \mu_0(z) \Big|} \\ & = & \Big| \int_{\partial_0}
\mathrm{tr} \big( w \beta_0(z) v_z^* c^{p'} s^* b^{p'} u_z^* \alpha_0(z)
\big) \, d\mu_0(z) \Big| \\ & \le & \Lambda \, \Big\| \int_{\partial_0}
\alpha_0(z)^* \alpha_0(z) \, d\mu_0(z) \Big\|_{M_m(\mathcal{M})}^{\frac12} \Big\| \int_{\partial_0} \beta_0(z)
\hskip1.5pt \beta_0(z)^* \, d\mu_0(z) \Big\|_{M_m(\mathcal{M})}^{\frac12},
\end{eqnarray*}
where $$\Lambda \, = \, \|w\|_{M_m(\mathcal{M} \bar\ten \mathcal{R})} \, \sup_{z \in \partial_0} \big\| v_z^* c^{p'} s^* b^{p'} u_z^* \big\|_{L_1(M_m(\mathcal{M} \bar\ten \mathcal{R}))} \, = \, \Lambda_1 \Lambda_2.$$ We have $\Lambda_1 \le 1$ by hypothesis, while H\"older's inequality gives $$\Lambda_2 \, \le \, \|c^{p'}\|_{L_2(M_m(\mathcal{M}))} \|s^*\|_{L_\infty(M_m(\mathcal{M}); L_1(\mathcal{R}))} \|b^{p'}\|_{L_2(M_m(\mathcal{M}))} \, \le \, 1.$$ This proves \eqref{claimholder}. Similarly, we have
\begin{eqnarray*}
\lefteqn{\Big| \int_{\partial_1} \big\langle g(z), \xi(z)
\big\rangle \, d \mu_1(z) \Big| \ = \ \Big| \int_{\partial_1}
\big\langle \tilde{u}_z s \tilde{v}_z, a^p \alpha_1(z) w
\beta_1(z) d^p \, d\mu_1(z) \big\rangle \Big|} \\ & \le &
\int_{\partial_1} \big\| a^p \alpha_1(z) w \beta_1(z) d^p
\big\|_{L_1(M_m(\mathcal{M}); \mathcal{R})} \, d\mu_1(z) \\ & \le
& \Big( \int_{\partial_1} \|a^p \alpha_1(z)\|_2^2 \, d\mu_1(z)
\Big)^{\frac12} \|w\|_{M_m(\mathcal{M} \bar\ten \mathcal{R})}
\Big( \int_{\partial_1} \|\beta_1(z) d^p\|_2^2 \, d\mu_1(z)
\Big)^{\frac12}\\ & \le & \Big\| \int_{\partial_1}
\mathrm{tr}_\mathcal{M} \big( \alpha_1(z) \alpha_1(z)^* \big) \,
d\mu_1(z) \Big\|_{M_m}^{\frac12} \Big\| \int_{\partial_1}
\mathrm{tr}_\mathcal{M} \big( \beta_1(z)^* \beta_1(z) \big) \,
d\mu_1(z) \Big\|_{M_m}^{\frac12}.
\end{eqnarray*}
In the last inequality we use that $a,d$ are norm $1$ in
$S_{2p}^m$. Summarizing, we get
\begin{eqnarray*}
\big| \langle x,y \rangle \big| & \le & (1 - \mbox{$\frac1p$})
\hskip1pt \Big| \int_{\partial_0} \langle g,\xi \rangle \,
d\mu_{0} \Big| + \mbox{$\frac1p$} \hskip1pt \Big|
\int_{\partial_1} \langle g, \xi \rangle \, d\mu_{1} \Big|
\\ & \le & (1 - \mbox{$\frac1p$}) \,
\|\alpha_0\|_{M_m(L_2^c(\partial_0) \bar\ten \mathcal{M})}
\, \|\beta_0\|_{M_m(\mathcal{M} \bar\ten L_2^r(\partial_0))} \\
[5pt] & + & \mbox{$\frac1p$} \,
\|\alpha_1\|_{M_m(L_2^r(\partial_1) \ten_h L_2^r(\mathcal{M}))} \,
\|\beta_1\|_{M_m(L_2^c(\mathcal{M}) \ten_h L_2^c(\partial_1))} \\
& \le & \Big( (1 - \mbox{$\frac1p$}) \, \|\alpha_0\|^2 +
\mbox{$\frac1p$} \, \|\alpha_1\|^2 \Big)^{\frac12} \Big( (1 -
\mbox{$\frac1p$}) \, \|\beta_0\|^2 + \mbox{$\frac1p$} \,
\|\beta_1\|^2 \Big)^{\frac12} \, \le \, 1.
\end{eqnarray*}
Therefore, the inclusion \eqref{ContractDual} is contractive and
the assertion follows by duality. \fin

\begin{theorem} \label{transfor} Given $1 < p < \infty$ and
$\mathcal{M}, \mathcal{R}$ as above assuming in addition that
$\mathcal{R}$ is finite dimensional. Then there exist states $\phi_n$ on
$M_{n}$, positive integers $\mathrm{k}_n$ and elements $\xi_n \in L_1(M_n)$
such that we have a complete embedding $$x \in L_p(\mathcal{M};
\mathcal{R}) \mapsto \Big( \sum_{j=1}^{\mathrm{k}_n} \pi_{tens}^j(\xi_n \ten x) \ten
\delta_j \Big) \in \prodd_{n,\U} L_1 \big( M_n(\mathcal{M})^{\ten_{\mathrm{k}_n}};
\ell_{\infty}^{\mathrm{k}_n}(\mathcal{R}) \big).$$
\end{theorem}

\dem When restricted to analytic functions, the operator $\Lambda(f|_{\partial_0}) = f|_{\partial_1}$ is densely defined and injective. In combination with \eqref{Harmonic}, this allows us to see the subspace of analytic functions on $\mathcal{S}$ vanishing at $1/p$ as the annihilator of the graph of $\Lambda$,
conveniently regarded as a space of analytic functions. The reader
is referred to \cite[Remark 2.2]{JP4} for further details.
Moreover, it is also observed in \cite{JP4} that we can replace
$\Lambda$ by a strictly positive diagonal operator
$\mathsf{d}_\la^{-1}$ on $\ell_2$ without changing the operator
space structure. In other words, we have complete isomorphisms
$$u_r: \big( L_2^c(\partial_0) \oplus L_2^r(\partial_1) \big)
\Big/ \big\{ f \ \mbox{analytic s.t.} \ f(\mbox{$\frac1p$})=0
\big\} \to \big( C \oplus R \big) \big/
graph(\mathsf{d}_\lambda^{-1})^\perp,$$
$$u_c: \big( L_2^r(\partial_0) \oplus
L_2^c(\partial_1) \big) \Big/ \big\{ f \ \mbox{analytic s.t.} \
f(\mbox{$\frac1p$})=0 \big\} \to \big( R \oplus C \big) \big/
graph(\mathsf{d}_\lambda^{-1})^\perp.$$ For further reference, we
set $\xi_r = u_r(1)$ and $\xi_c = u_c(1)$, where $1$ denotes the
constant function $1$ on the strip. Here $\mathsf{d}_\lambda:
\ell_2 \to \ell_2$ is a diagonal operator
$\mathsf{d}_\lambda(\delta_k) = \lambda_k e_k$ with $0 < \lambda_k
< \infty$ and the fact that we may consider the same operator in
both cases is also justified in \cite{JP4}. The exact same argument
mentioned above shows that, if we tensorize with the identity map,
we also obtain complete isomorphisms $$u_r: \mathcal{S}^r_\mathcal{M}
/ \mathsf{H}_{r,0} \to \Big[ \big( C \bar\ten \mathcal{M} \big) \oplus \big( R \ten_h
L_2^r(\mathcal{M}) \big) \Big] \Big/ graph(\mathsf{d}_\lambda^{-1})^\perp,$$
$$u_c: \mathcal{S}^c_\mathcal{M} / \mathsf{H}_{c,0} \to \Big[ \big(
\mathcal{M} \bar\ten R \big) \oplus \big( L_2^c(\mathcal{M})
\ten_h C \big) \Big] \Big/ graph(\mathsf{d}_\lambda^{-1})^\perp.$$
Let us define
$$\widetilde{\mathcal{H}}_r(\mathcal{M},1/p) =
u_r(\mathcal{H}_r(\mathcal{M},1/p)) \quad \mbox{and} \quad
\widetilde{\mathcal{H}}_c(\mathcal{M},1/p) =
u_c(\mathcal{H}_c(\mathcal{M},1/p)).$$ As in \cite{JP4}, we
observe that $$\Big[ \big( C \bar\ten \mathcal{M} \big) \oplus
\big( R \ten_h L_2^r(\mathcal{M}) \big) \Big] \Big/
graph(\mathsf{d}_\lambda^{-1})^\perp,$$ $$\Big[ \big( \mathcal{M}
\bar\ten R \big) \oplus \big( L_2^c(\mathcal{M}) \ten_h C \big)
\Big] \Big/ graph(\mathsf{d}_\lambda^{-1})^\perp,$$ can also be
understood as $\mathcal{K}$-spaces. Indeed, if we use anti-linear
duality, we have $\mathsf{d}_\lambda^{-1}: x \in C \mapsto
(\mathsf{d}_\lambda^{-1} x)^\mathrm{t} \in R$ in the first case
and $\mathsf{d}_\lambda^{-1}: x \in R \mapsto (x
\mathsf{d}_\lambda^{-1})^\mathrm{t} \in C$ in the second one. This
means that $graph(\mathsf{d}_\lambda^{-1})^\perp$ is spanned by
elements of the form $(-\mathsf{d}_\lambda^{-1} x,x)$ in the first
quotient and by $(-x \mathsf{d}_\lambda^{-1},x)$ in the second
one. In conclusion this allows us to cb-embed
\begin{eqnarray*}
\big[ u_r \ten id \ten u_c \big] \big(
\mathcal{H}_{\frac1p}(\mathcal{M} \bar\ten \mathcal{R}) \big) & =
& \widetilde{\mathcal{H}}_{\frac1p}(\mathcal{M} \bar\ten
\mathcal{R}) \\ & = & \widetilde{\mathcal{H}}_r(\mathcal{M},1/p)
\ten_{\mathcal{M},h} \big( \mathcal{M} \bar\ten \mathcal{R} \big)
\ten_{\mathcal{M},h} \widetilde{\mathcal{H}}_c(\mathcal{M},1/p)
\end{eqnarray*}
into a four term $\mathcal{K}$-space $\mathcal{K}_\lambda$ with
norm given by
$$\|x\|_{\mathcal{K}_\lambda} = \inf_{x= x_1 +
\mathsf{d}_\lambda^{-1} x_2 \mathsf{d}_\lambda^{-1} +
\mathsf{d}_\lambda^{-1} x_3 + x_4 \mathsf{d}_\lambda^{-1}}
\sum_{j=1}^4 \|x_j\|_{\mathsf{E}_j},$$ where the spaces
$\mathsf{E}_1, \mathsf{E}_2, \mathsf{E}_3, \mathsf{E}_4$ are given
by
\begin{eqnarray*}
\mathsf{E}_1 & = & \mathcal{M} \bar\ten \mathcal{B}(\ell_2)
\bar\ten \mathcal{R}, \\ \mathsf{E}_2 & = & L_1(\mathcal{M}
\bar\ten \mathcal{B}(\ell_2); \mathcal{R}), \\ \mathsf{E}_3 & = &
L_2^r(\mathcal{M} \bar\ten \mathcal{B}(\ell_2)) \ten_{\mathcal{M}
\bar\ten \mathcal{B}(\ell_2),h} (\mathcal{M} \bar\ten
\mathcal{B}(\ell_2) \bar\ten \mathcal{R}), \\ \mathsf{E}_4 & = &
(\mathcal{M} \bar\ten \mathcal{B}(\ell_2) \bar\ten \mathcal{R})
\ten_{\mathcal{M} \bar\ten \mathcal{B}(\ell_2),h}
L_2^c(\mathcal{M} \bar\ten \mathcal{B}(\ell_2)).
\end{eqnarray*}
Indeed, all these spaces can be essentially obtained by applying
Remark \ref{AmalgTrick}. For instance, to obtain $\mathsf{E}_4$ we
have to show that the given space comes from the choice $C
\bar\ten \mathcal{M}$ and $L_2^c(\mathcal{M}) \ten_h C$ for the
left and right spaces. More concretely, we have a completely
isometric embedding of $$(C \bar\ten \mathcal{M})
\ten_{\mathcal{M},h} (\mathcal{M} \bar\ten \mathcal{R})
\ten_{\mathcal{M},h} (L_2^c(\mathcal{M}) \ten_h C)$$ into
$\mathsf{E}_4$. However, according to Remark \ref{AmalgTrick} we
may embed it cb-isometrically into $(C \bar\ten \mathcal{M} \ten_h
R) \ten_{\mathcal{M} \bar\ten \mathcal{B}(\ell_2),h} (C \ten_h
\mathcal{M} \bar\ten \mathcal{R} \ten_h R) \ten_{\mathcal{M}
\bar\ten \mathcal{B}(\ell_2),h} (C \ten_h L_2^c(\mathcal{M})
\ten_h C)$ which in turn embeds in
$$(\mathcal{M} \bar\ten \mathcal{B}(\ell_2)) \ten_{\mathcal{M}
\bar\ten \mathcal{B}(\ell_2),h} (\mathcal{M} \bar\ten
\mathcal{B}(\ell_2) \bar\ten \mathcal{R}) \ten_{\mathcal{M}
\bar\ten \mathcal{B}(\ell_2),h} L_2^c(\mathcal{M} \bar\ten
\mathcal{B}(\ell_2)).$$ This completes the argument since the
latter space is $\mathsf{E}_4$. On the other hand, since
$\mathcal{R}$ is of finite dimension ($m$ say) we know that
topologically we may write $\mathcal{K}_{\la}^*$ as follows
\begin{eqnarray*}
\mathcal{K}_{\la}^* & \simeq & \big( (\mathcal{M} \bar\ten
\mathcal{B}(\ell_2))^m \big)^* \cap \big( \mathsf{d}_{\la}^{-1} \mathcal{M}
\bar\ten \mathcal{B}(\ell_2) \mathsf{d}_{\la}^{-1} \big)^m \\
& \cap & \big( \mathsf{d}_{\la}^{-1} L_2(\mathcal{M} \bar{\ten}
\mathcal{B}(\ell_2)) \big)^m \cap \big( L_2(\mathcal{M} \bar{\ten}
\mathcal{B}(\ell_2) \mathsf{d}_{\la}^{-1} )^m \\ & = & L_1
\big(\mathcal{M} \bar\ten \mathcal{B}(\ell_2) \big)^m \cap
\big(\mathsf{d}_{\la}^{-1} \mathcal{M} \bar\ten
\mathcal{B}(\ell_2) \mathsf{d}_{\la}^{-1} \big)^m \\ & \cap &
\big( \mathsf{d}_{\la}^{-1} L_2(\mathcal{M}
\bar{\ten}\mathcal{B}(\ell_2)) \big)^m \cap \big( L_2(\mathcal{M}
\bar{\ten} \mathcal{B}(\ell_2)) \mathsf{d}_{\la}^{-1} \big)^m.
\end{eqnarray*}
Here we used the fact that a matrix $[x_{ij}]$ with $x_{ij}\in
\mathcal{M}\subset L_1(\mathcal{M})$ belonging to $(\mathcal{M}
\bar\ten \mathcal{B}(\ell_2))^{**}$ already belongs to
$L_1(\mathcal{M} \bar\ten \mathcal{B}(\ell_2))$. Let $(p_n)_{n \ge
1}$ be an increasing sequence of orthogonal projections commuting
with $\mathsf{d}_{\la}$ and converging strongly to $\mathbf{1}$.
Then we deduce that for $x\in \mathcal{K}_{\la}$ $\lim_{n\to
\infty} \big\langle p_nxp_n,y \big\rangle = \langle x,y \rangle$
because we can use convergence in the norm of $L_1$ or $L_2$ on at
least one side of the bracket. Using a weak$^*$-limit we obtain
\[ \|x\|_{\mathcal{K}_{\la}} = \|x\|_{\mathcal{K}_{\la}^{**}}
\le \lim_{n,\U} \big\| (1\ten p_n\ten 1) x (1\ten p_n\ten 1)
\big\|_{\mathcal{K}_{\la}} \le \|x\|_{\mathcal{K}_{\la}} \] for
any free ultrafilter $\U$ on the integers. Therefore, allowing to
take ultraproducts (as we do) in the final space, it suffices to
consider the finite-dimensional case where $\mathcal{B}(\ell_2)$ is replaced
by the matrix algebra $M_n$. Define on $M_n$
$$\psi_n \Big( \sum_{i,j=1}^n \alpha_{ij} e_{ij} \Big) =
\sum_{k=1}^n \la_k^2 \, \alpha_{kk} \quad \mbox{and} \quad \phi_n(x)
= \frac{\psi_n(x)}{\psi_n(\mathbf{1})}.$$ Since the original
$\mathsf{d}_{\la}$ is unbounded, we may assume that $\sum_k \la_k^2
> 1$. Then, by approximation we can indeed assume that
$\psi_n(\mathbf{1}) = \mathrm{k}_n$ is an integer, see \cite{JP4}
for more details. If $d_{\psi_n}$ and $d_{\phi_n}$ stand for the
corresponding densities, we clearly have $$\mathsf{d}_\lambda =
d_{\psi_n}^{\frac12} = \sqrt{\mathrm{k}_n} \,
d_{\phi_n}^{\frac12}.$$ In particular, we may replace
$\mathcal{K}_\lambda$ by $\mathcal{K}_n$ with
$$\|x\|_{\mathcal{K}_n} = \inf_{x=\sum_1^4 x_j}
\|x_1\|_{\mathsf{E}_{n1}} + \mathrm{k}_n \, \big\|
d_{\phi_n}^{\frac12} x_2 d_{\phi_n}^{\frac12}
\big\|_{\mathsf{E}_{n2}} + \sqrt{\mathrm{k}_n} \, \big\|
d_{\phi_n}^{\frac12} x_3 \big\|_{\mathsf{E}_{n3}} +
\sqrt{\mathrm{k}_n} \, \big\| x_4 d_{\phi_n}^{\frac12}
\big\|_{\mathsf{E}_{n4}}$$ and where $\mathsf{E}_{nj}$ is the
result of replacing in $\mathsf{E}_j$ the algebra
$\mathcal{B}(\ell_2)$ by $M_n$. We can identify this space in the
terminology of Theorem \ref{4-term}. Namely, if we fix a positive
integer $m$ and set $\mathcal{E}_m = id \ten \phi_n \ten \varphi:
M_m \ten M_n \ten \mathcal{M} \to M_m$, then it is easily checked
that we have the following isometric isomorphism
$$S_1^m(\mathcal{K}_n) =
\mathcal{K}_{1,\infty}^{\mathrm{k}_n} \big( M_{mn}(\mathcal{M}),
\mathcal{E}_m; \mathcal{R} \big).$$ For instance, according to Remark \ref{AmalgTrick},
$S_1^m(\mathsf{E}_{n3})$ can be written as
$$\Big[ R_m \ten_h L_2^r(M_n(\mathcal{M}))
\ten_h R_m \Big] \ten_{M_{mn}(\mathcal{M}),h} \Big[ C_m \ten_h
M_{n}(\mathcal{M} \bar\ten \mathcal{R}) \ten_h R_m \Big]
\ten_{M_m,h} C_{m^2},$$ which in turn is cb-isometric to $$L_2^r(M_{mn}(\mathcal{M})) \ten_{M_{mn}(\mathcal{M}),h} (M_{mn}(\mathcal{M}) \bar\ten \mathcal{R}) \ten_{M_m,h} L_2^c(M_m).$$ In other words, we obtain the space $L_1^c(M_{mn}(\mathcal{M}), \mathcal{E}_m; \mathcal{R})$. The term
$\mathsf{E}_{n4}$ is handled in the same way while the terms
$\mathsf{E}_{n1}$ and $\mathsf{E}_{n2}$ are even easier. Applying
Theorem \ref{4-term} for tensor copies, we have an embedding
$$x \in \mathcal{K}_{1,\infty}^{\mathrm{k}_n} \big(
M_{mn}(\mathcal{M}), \mathcal{E}_m;\mathcal{R} \big) \mapsto
\sum_{j=1}^{\mathrm{k}_n} \pi_{tens}^j(x) \ten \delta_j \in L_1
\big( \mathcal{A}_{m,n} ;\ell_\infty^{\mathrm{k}_n}(\mathcal{R})
\big),$$ with $\mathcal{A}_{m,n} = M_m \ten
M_n(\mathcal{M})^{\ten_{\mathrm{k}_n}}$. In particular, this
produces a cb-embedding $$w_n: x \in \mathcal{K}_n \mapsto
\sum_{j=1}^{\mathrm{k}_n} \pi_{tens}^j(x) \ten \delta_j \in L_1
\big( M_n(\mathcal{M})^{\ten_{\mathrm{k}_n}} ;
\ell_\infty^{\mathrm{k}_n}(\mathcal{R}) \big),$$ with relevant
constants independent of $n$. Finally, the construction of our
complete embedding is as follows. We first apply Proposition
\ref{emb} and then we proceed as above. Namely, if $u = u_r \ten
id \ten u_c$, we have
$$L_p(\mathcal{M}; \mathcal{R}) \stackrel{j}{\longrightarrow}
\mathcal{H}_{\frac1p}(\mathcal{M} \bar\ten \mathcal{R})
\stackrel{u}{\longrightarrow}
\widetilde{\mathcal{H}}_{\frac1p}(\mathcal{M} \bar\ten
\mathcal{R}) \stackrel{id}{\longrightarrow} \overline{\bigcup_{n
\ge 1} \mathcal{K}_n}$$ and we have constructed a complete
embedding
$$\prodd_{n, \mathcal{U}} w_n: \overline{\bigcup_{n \ge 1}
\mathcal{K}_n} \longrightarrow \prodd_{n, \mathcal{U}} L_1 \big(
M_n(\mathcal{M})^{\ten_{\mathrm{k}_n}} ;
\ell_\infty^{\mathrm{k}_n}(\mathcal{R}) \big).$$ Let $q_n:
\mathcal{B}(\ell_2) \to M_n$ be the projection into the upper left
corner and let us define the element $\xi_n = q_n \xi_r q_n \xi_c
q_n \in L_1(M_n)$. Since $u(j(x)) = u(1 \ten x) = \xi_r \xi_c \ten
x$, the form of the embedding $\prod_{n, \mathcal{U}} w_n \circ u
\circ j$ is the one given in the assertion. \fin

\section{Mixed-norm transference and applications}
\label{Section5}

Given a Hilbert space $\mathcal{H}$, we shall write $\mathcal{H}_r$
and $\mathcal{H}_c$ to denote the row and column operator space
structures on $\mathcal{H}$. Accordingly, $\mathcal{H}_{r_p}$ and
$\mathcal{H}_{c_p}$ stand for the complex interpolation spaces
$[\mathcal{H}_r, \mathcal{H}_c]_{1/p}$ and $[\mathcal{H}_c,
\mathcal{H}_r]_{1/p}$, respectively. Let us fix $1 \le p \le q \le
\infty$ and $n \ge 1$ a positive integer. By Pisier's exercise
\cite{P3} and some refinements \cite{J3,JP4,X3}, we may construct
complete embeddings $$\al_q: C_q^n \to
L_2^{c}(\Omega,\mu_q;\ell_2^n) + L_2^{c_p}(\Omega,\nu_q;
\ell_2^n),$$ $$\beta_q: R_q^n \hskip1pt \to L_2^{r}(\Omega,\mu_q;
\ell_2^n) + L_2^{r_p}(\Omega,\nu_q; \ell_2^n),$$ for suitable
measures $\mu_q$ and $\nu_q$ on a finite set $\Omega =
\{1,2,\ldots,m\}$ with $m$ depending on $n$. In fact, an elaborated version of this result was already used in the previous section. A much more concrete approach is available in \cite[Lemma 2.2]{JP5}. Let $\mu_i = \mu_q\{i\}$ and $\nu_i = \nu_q\{i\}
= \la_i \mu_i$ for some $\la_i > 0$. Let us write
$\mathsf{d}_\lambda$ for the diagonal operator on $\ell_2^m$
determined by the $\la_i$'s. That is, $\mathsf{d}_\la = \sum_k \la_k
e_{kk}$. The symbol $+$ above refers as in the previous section to
the quotient of the direct sums
$$L_2^{c}(\Omega,\mu_q;\ell_2^n) \oplus L_2^{c_p}(\Omega,\nu_q;
\ell_2^n) \quad \mbox{and} \quad L_2^{r}(\Omega,\mu_q; \ell_2^n)
\oplus L_2^{r_p}(\Omega,\nu_q; \ell_2^n)$$ by the subspace
$$S = \Big\{ \big( a_{ij}, - \la_i^{-\frac12}a_{ij} \big) \, \big|
\ 1 \le i \le m, \ 1 \le j \le n \Big\}.$$ More concretely, in the
first case $a = (a_{ij}) \in L_2^{c}(\Omega,\mu_q; \ell_2^n)$ is a
column with $m$ entries in $\ell_2^n$, while in the second case $a
= (a_{ij}) \in L_2^r(\Omega, \mu_q; \ell_2^n)$ is a row. In
particular, we may write $S$ in each case as follows $$S_\alpha =
\Big\{ \big( a, - \mathsf{d}_\la^{-\frac12} a \big) \, \big| \ a
\in L_2^{c}(\Omega,\mu_q; \ell_2^n) \Big\} = \Big\{ \big( -
\mathsf{d}_\la^{\frac12} a, a \big) \, \big| \ a \in
L_2^{c}(\Omega,\mu_q; \ell_2^n) \Big\},$$
$$S_\beta = \Big\{ \big( a, - a \mathsf{d}_\la^{-\frac12} \big) \,
\big| \ a \in L_2^{r}(\Omega,\mu_q; \ell_2^n) \Big\} = \Big\{
\big( - a \mathsf{d}_\la^{\frac12}, a \big) \, \big| \ a \in
L_2^{r}(\Omega,\mu_q; \ell_2^n) \Big\}.$$ The embedding is of the
form $$\al_q(a) = \mathbf{1}_\Omega \ten a + S_\al.$$ The formula
for $\beta_q$ is the same. Let us define $$S_{\alpha \beta} =
S_\al \ten_h \Big( L_2^{r}(\Omega,\mu_q; \ell_2^n) \oplus
L_2^{r_p}(\Omega,\nu_q; \ell_2^n) \Big) + \Big(
L_2^{c}(\Omega,\mu_q; \ell_2^n) \oplus L_2^{c_p}(\Omega,\nu_q;
\ell_2^n) \Big) \ten_h S_\beta.$$

\begin{lemma} \label{PEx}
If $1\le p\le q\le \infty$, we have a cb-embedding
$$v_q: \ell_q^n \to \Big( L_2^{c}(\Omega,\mu_q; \ell_2^n) +
L_2^{c_p}(\Omega,\nu_q; \ell_2^n) \Big) \ten_h \Big(
L_2^{r}(\Omega,\mu_q; \ell_2^n) + L_2^{r_p}(\Omega,\nu_q; \ell_2^n)
\Big),$$ $$v_q \big( \xi_1, \xi_2, \ldots, \xi_n \big) =
\sum_{k=1}^n \xi_k \, \al_q(e_{k1}) \ten \beta_q(e_{1k}) = \Big[
\sum_{k=1}^n \xi_k\,  \Big( \sum_{i,j=1}^m e_{ij} \Big) \ten e_{kk}
\Big] + S_{\alpha \beta}.$$
\end{lemma}

\dem It follows from our considerations above and $\ell_q^n
\subset C_q^n \ten_h R_q^n$. \fin

The embedding $v_q$ is special in the sense that its range is
contained in the subalgebra $M_m \ten \ell_{\infty}^n$, after a
suitable change of variables. To explain this we recall that
$\mu_q\{i\} = \mu_i$ and $\nu_q\{i\} = \la_i \mu_i$. Therefore, the
map $$j: (a_{ij}) \in L_2^c(\Omega,\mu_q; \ell_2^n) \mapsto \big(
\sqrt{\mu_i} \, a_{ij} \big) \in C_{mn}$$ is a complete isometry. To
respect the sum operation, we have to apply $j$ also on
$L_2^{c_p}(\Omega,\nu_q; \ell_2^n)$. If $\lambda$ stands for the
measure on $\Omega$ given by $\lambda\{i\} = \lambda_i$, we find
another complete isometry $j: L_2^{c_p}(\Omega,\mu_q;\ell_2^n) \to
L_2^{c_p}(\Omega,\lambda,\ell_2^n)$. Hence, applying the same
argument for the other side and using the terminology
$$L_2^{c_p}(\Omega,\lambda; \ell_2^n) =
\mathsf{d}_\lambda^{\frac12} C_p^{mn} \quad \mbox{and} \quad
L_2^{r_p}(\Omega,\lambda; \ell_2^n) = R_p^{mn}
\mathsf{d}_\lambda^{\frac12},$$ we find a complete isometry
\begin{eqnarray*}
\lefteqn{\hskip-130pt J: \big( L_2^c(\Omega,\mu_q; \ell_2^n)
\oplus L_2^{c_p}(\Omega,\nu_q; \ell_2^n) \big) \ten_h \big(
L_2^r(\Omega,\mu_q, \ell_2^n) \oplus L_2^{r_p}(\Omega,\nu_q;
\ell_2^n) \big)} \\ [3pt] \hskip130pt  \null & \to & \big( C_{mn}
\oplus \mathsf{d}_\lambda^{\frac12} C_p^{mn} \big) \ten_h \big(
R_{mn} \oplus R_p^{mn} \mathsf{d}_\lambda^{\frac12} \big)
\end{eqnarray*}
with $J = (j,j) \ten (j,j)$. Passing to quotients, we may replace
$\oplus$ by $+$. The key observation here is that algebraically we
have $J \, (v_q(\ell_q^n)) \subset M_m\otimes \ell_{\infty}^n$.
Indeed, note
\begin{equation} \label{jjeelaee}
J \Big[ \Big(\sum_{i,j=1}^m e_{ij} \Big) \ten e_{kk} \Big] = \Big(
\sum_{i,j=1}^m \mu_i^{\frac12} \mu_j^{\frac12} e_{ij} \Big) \ten
\delta_{k}.
\end{equation}
Before we proceed with our next result, we review the
$\mathcal{K}_{p,\infty}^\mathrm{k}(M_m,\phi)$ space given by a
state $\phi$. In this case, given any pair $u,v \in\{2p,\infty\}$,
the inclusion map $S_{(u,v)}^m \subset S_p^m$ depends on $\phi$.
Indeed, we may and will assume that $\phi(x) = \sum_k \phi_k
x_{kk}$. Then, the density $d_{\phi}$ is indeed a diagonal
operator with coefficients $\phi_k$ and we have Kosaki's embedding
$$x \in S_{(u,v)}^m \mapsto d_\phi^{\frac{1}{2p}-\frac1u} x
\, d_\phi^{\frac{1}{2p}-\frac1v} \in S_p^m.$$ Therefore, we find
$$\|x\|_{\mathcal{K}_{p,\infty}^\mathrm{k}(M_m,\phi)} = \inf
\Big\{\mathrm{k}^{\frac1p} \|x^1\|_{S_p^m} +
\mathrm{k}^{\frac{1}{2p}} \|x^2\|_{S_{2p}^m} +
\mathrm{k}^{\frac{1}{2p}} \|x^3\|_{S_{2p}^m} + \|x^4\|_{M_m}
\Big\},$$ where the infimum runs over $x = x^1 + x^2
d_{\phi}^{\frac{1}{2p}} + d_{\phi}^{\frac{1}{2p}} x^3 +
d_{\phi}^{\frac{1}{2p}} x^4 d_{\phi}^{\frac{1}{2p}}$. This gives
$$\|x\|_{\mathcal{K}_{p,\infty}^\mathrm{k}(M_m,\phi)} = \inf
\Big\{ \big\| d_{\mathrm{k}\phi}^{\frac{1}{2p}} x^1
d_{\mathrm{k}\phi}^{\frac{1}{2p}} \big\|_{S_p^m} + \big\|
d_{\mathrm{k}\phi}^{\frac{1}{2p}} x^2 \big\|_{S_{2p}^m} + \big\|
x^3 d_{\mathrm{k}\phi}^{\frac{1}{2p}} \big\|_{S_{2p}^m} +
\|x^4\|_{M_m} \Big\},$$ where this time the infimum is taken over
$d_{\phi}^{-\frac{1}{2p}} x \, d_{\phi}^{-\frac{1}{2p}} = x^1 +
x^2 + x^3 + x^4$. A similar calculation applies in the
operator-valued setting. In the following result, we shall use the
notation $$\|x\|_{d^\alpha X d^\beta} = \|d^\alpha x \hskip1pt
d^\beta\|_{X}.$$

\begin{lemma} \label{concreteK}
Let us consider a von Neumann algebra $\mathcal{M}$, positive
integers $m,n,\mathrm{k}$ and a state $\phi$ on $M_m$. Let
$\mathcal{E_M}: M_m(\mathcal{M}) \ten \ell_\infty^{\mathrm{k} n}
\to \mathcal{M}$ stand for the conditional expectation
$$\mathcal{E_M} \Big(\sum_{i=1}^\mathrm{k} \sum_{j=1}^n x_{ij}
\ten \delta_{ij} \Big) = \frac{1}{n \mathrm{k}}
\sum_{i=1}^\mathrm{k} \sum_{j=1}^n \phi \ten id (x_{ij}).$$ Let
$\mathcal{R} = M_m(\mathcal{M}) \ten \ell_\infty^n$ and consider
the space $X_{\phi,n}^p(\mathcal{M})$ defined by
$$L_p(\mathcal{M};\ell_{\infty}^{n}(M_m)) +
d_{\mathrm{k}\phi}^{\frac{1}{2p}} L_{2p}(\mathcal{R})
L_{2p}(\mathcal{M}) + L_{2p}(\mathcal{M}) L_{2p}(\mathcal{R})
d_{\mathrm{k} \phi}^{\frac{1}{2p}} +
d_{\mathrm{k}\phi}^{\frac{1}{2p}} L_p(\mathcal{R})
d_{\mathrm{k}\phi}^{\frac{1}{2p}}.$$ Then the following identity
holds
$$\Big\|\sum_{i=1}^\mathrm{k} \sum_{j=1}^n d_{\phi}^{\frac{1}{2p}}
x_j \, d_{\phi}^{\frac{1}{2p}} \otimes \delta_{ij}
\Big\|_{\mathcal{K}_{p,\infty}^{\mathrm{k}n}(M_m(\M) \ten
\ell_\infty^{\mathrm{k}n},\E_{\M})} = \Big\| \sum_{j=1}^n x_j \ten
\delta_j \Big\|_{X_{\phi,n}^p(\M)}.$$
\end{lemma}

\dem The subspace of sequences
$$\big(x_1, \ldots, x_1, x_2, \ldots, x_2, \ldots, x_n, \ldots, x_n
\big),$$ with every $x_j$ repeated $\mathrm{k}$ times is
complemented in $\mathcal{K}_{p,\infty}^{\mathrm{k}n}(M_m(\M) \ten
\ell_\infty^{\mathrm{k}n},\E_{\M})$ since it is complemented in
the four spaces composing it. Let us assume for simplicity that
$\M$ has a normalized trace $\tau$. Our reference state and trace
in the construction of the Haagerup $L_p$ spaces are
$$\psi(x_{ij}) = \frac{1}{n \mathrm{k}} \sum_{i=1}^\mathrm{k}
\sum_{j=1}^n \phi \ten \tau(x_{ij}) \quad \mbox{and} \quad
\mathrm{tr}(x_{ij}) = \frac{1}{n \mathrm{k}} \sum_{i=1}^\mathrm{k}
\sum_{j=1}^n \mathrm{tr}_{M_m} \ten \tau(x_{ij}).$$ The density is
given by $d_\psi = d_\phi \ten \mathbf{1}_\M$ and letting
$\widehat{\M}_{\mathrm{k}mn} = M_m(\M) \ten
\ell_\infty^{\mathrm{k}n}$ we have
\begin{align*}
&\mathcal{K}_{p,\infty}^{\mathrm{k}n} \big(
\widehat{\M}_{\mathrm{k}mn}, \mathcal{E_M} \big) = (n
\mathrm{k})^{\frac1p} \ d_\psi^{\frac{1}{2p}}
L_p(\widehat{\M}_{\mathrm{k}mn}) \, d_\psi^{\frac{1}{2p}} \\
& + (n \mathrm{k})^{\frac{1}{2p}} \, d_\psi^{\frac{1}{2p}}
L_{2p}(\widehat{\M}_{\mathrm{k}mn}) L_{2p}(\M)  +  (n
\mathrm{k})^{\frac{1}{2p}} \, L_{2p}(\M)
L_{2p}(\widehat{\M}_{\mathrm{k}mn}) d_\psi^{\frac{1}{2p}} \\
& + L_{2p}(\M) \, L_{\infty}(\widehat{\M}_{\mathrm{k}mn}) \,
L_{2p}(\M) \ = \ Z_1 + Z_2 + Z_3 + Z_4.
\end{align*}
This means that $$\Big\|\sum_{i=1}^\mathrm{k} \sum_{j=1}^n
d_{\phi}^{\frac{1}{2p}} x_j \, d_{\phi}^{\frac{1}{2p}} \otimes
\delta_{ij} \Big\|_{\mathcal{K}_{p,\infty}^{\mathrm{k}n}
(\widehat{\M}_{\mathrm{k}mn},\E_{\M})} = \inf_{x_j = x_j^1 + x_j^2
+ x_j^3 + x_j^4} \sum_{s=1}^4 \Big\| \summ_{i,j} x_j^s \otimes
\delta_{ij} \Big\|_{Z_s}.$$ Let us compute the four norms
\begin{eqnarray*}
\Big\|\summ_{i,j} x_j \otimes \delta_{ij} \Big\|_{Z_1} & = & (n
\mathrm{k})^{\frac1p} \Big( \frac{1}{n \mathrm{k}} \summ_{i,j}
\big\| d_\phi^{\frac{1}{2p}} x_j d_\phi^{\frac{1}{2p}}
\big\|_{L_p(M_m(\M))}^p \Big)^{\frac1p} \\ & = &
\mathrm{k}^{\frac1p} \Big( \sum_{j=1}^n \big\| d_\phi^{\frac{1}{2p}}
x_j d_\phi^{\frac{1}{2p}} \big\|_p^p \Big)^{\frac1p}\!\! \ = \
\mathrm{k}^{\frac1p} \Big\| \sum_{j=1}^n x_j \ten \delta_j
\Big\|_{d_\phi^{\frac{1}{2p}} L_p(\mathcal{R})
d_\phi^{\frac{1}{2p}}}, \\ \Big\| \summ_{i,j} x_j \otimes
\delta_{ij} \Big\|_{Z_2} & = & \mathrm{k}^{\frac{1}{2p}}
\inf_{d_\phi^{1/2p} x_j = z_j b} \Big( \sum_{j=1}^n
\|z_j\|_{L_{2p}(M_m(\M))}^{2p} \Big)^{\frac{1}{2p}}
\|b\|_{L_{2p}(\M)} \\ & = & \mathrm{k}^{\frac{1}{2p}} \Big\|
\sum_{j=1}^n x_j \ten \delta_j \Big\|_{d_\phi^{\frac{1}{2p}}
L_{2p}(\mathcal{R}) L_{2p}(\M)}, \\ \Big\| \summ_{i,j} x_j \otimes
\delta_{ij} \Big\|_{Z_4} & = & \inf_{x_j = a z_j b}
\|a\|_{L_{2p}(\M)} \sup_{1 \le j \le n} \|z_j\|_{M_m(\M)}
\|b\|_{L_{2p}(\M)} \\ & = & \Big\| \sum_{j=1}^n x_j \ten \delta_j
\Big\|_{L_p(\M; \ell_\infty^n(M_m))}.
\end{eqnarray*}
The $Z_3$-term is calculated as $Z_2$. The proof is complete. \fin

\begin{theorem} \label{vqq}
Let $1 \le p \le q \le \infty$ and set $m = |\Omega|$ as above.
Then, there exists a state $\phi_m$ on $M_m$ and a positive
integer $\mathrm{k}_m$ such that we have a complete embedding
$$L_p(\mathcal{M}; \ell_q^n) \to L_p \Big( \mathcal{M}
\bar\ten \Big[*_{\mathrm{k}_m n} \big(M_m, \phi_m \big) \Big] ;
\ell_\infty^{\mathrm{k}_m n} \Big)$$ given by the relation
$$\sum_{j=1}^n x_j \ten \delta_j \mapsto
\sum_{i=1}^{\mathrm{k}_m} \sum_{j=1}^n x_j \ten
\pi_{free}^{ij}(a_m) \ten \delta_{ij},$$ where $a_m = a_m(p,q) \in
M_m$. The relevant constants are independent of $m$ and $n$.
\end{theorem}

\dem By enlarging $m$ if necessary, we may assume that $\sum_i
\la_i^p = \mathrm{k}_m$ is an integer. Then we define the normalized
state $\phi_m(x) = \mathrm{k}_m^{-1} \sum_{i} \la_i^p x_{ii}$. We
observe that for arbitrary elements $x_j \in M_m \ten L_p(\M)$, the
right hand side of Lemma \ref{concreteK} coincides with the norm of
diagonal sequences (i.e. $mn \times mn$ matrices which are diagonal
on its $n$-component) in the space $$L_p \Big(\M; \big( C_{mn} +
\mathsf{d}_\lambda^\frac12 C_p^{mn} \big) \ten_h \big( R_{mn} +
R_p^{mn} \mathsf{d}_\lambda^\frac12 \big) \Big).$$ Indeed, we use
the properties of the Haagerup tensor product and the fact that
projection onto the diagonal $M_n \to \ell_{\infty}^n$ is completely
contractive with respect to all the four interpolation norms. The
embedding is given by $u \circ \big( id_{L_p(\M)} \ten (J \circ v_q)
\big)$ where $u$ is the first map from Theorem
\ref{Theorem-Sigma-pq} for $q=\infty$ and $\mathrm{k}_m n$ instead
of $n$. Note that $u$ is well-defined because of $J(v_q(\ell_q^n))
\subset M_m \otimes \ell_{\infty}^n$. Moreover, identity
\eqref{jjeelaee} tells us that $$a_m = d_\phi^{\frac{1}{2p}} \Big(
\sum_{i,j \le \mathrm{k}_m'} \sqrt{\mu_i \mu_j} e_{ij} \Big)
d_\phi^{\frac{1}{2p}}$$ with $\mathrm{k}_m' \le \mathrm{k}_m$. The
proof is complete. \fin

\begin{remark}
\emph{The constants in Theorem \ref{vqq} are also independent of
$p,q$ as far as we do not have $p \to 1$ and $q \to \infty$
simultaneously. The use of Theorem \ref{Theorem-Sigma-pq} produces
such singularity, see Remarks \ref{Singpq} and \ref{RemKfree2} for
further details.}
\end{remark}

\noindent {\bf Proof of Theorem B.} Let us consider
$$\sum_{k=1}^n \pi_k(x) \ten \delta_k \in L_p(\mathcal{A};
\ell_q^n).$$ According to Theorem \ref{vqq}, the following
equivalence holds
$$\Big\| \sum_{k=1}^n \pi_k(x) \ten \delta_k
\Big\|_{L_p(\mathcal{A}; \ell_q^n)} \sim_c \Big\|
\sum_{i=1}^{\mathrm{k}_m} \sum_{j=1}^n \pi_j(x) \ten
\pi_{free}^{ij} (a_m) \ten \delta_{ij} \Big\|_{L_p(\mathcal{A}
\ten \mathcal{B}; \ell_\infty^{\mathrm{k}_m n})},$$ with
$\mathcal{B} = *_{\mathrm{k}_m n} (M_m, \phi_m)$. Applying the
complete embedding that we constructed in Theorem \ref{transfor}
to the term on the right hand side (note that both Theorems \ref{transfor}
and \ref{vqq} provide constants independent of $m$ and $n$) we obtain a new term in the
ultraproduct $$\prodd_{s,\U} L_1 \Big( \big[ M_s \ten \mathcal{A} \ten
\mathcal{B} \big]^{\ten_{\mathrm{k}_s}}; \ell_\infty^{\mathrm{k}_m
\mathrm{k}_s n} \Big)$$ of the following form
$$\left( \sum_{i=1}^{\mathrm{k}_m} \sum_{j=1}^n
\sum_{w=1}^{\mathrm{k}_s} \pi_{tens}^w \Big[ \xi_s \ten \pi_j(x)
\ten \pi_{free}^{ij} (a_m) \Big] \ten \delta_{ijw} \right)_s$$ for
a fixed family of matrices $\xi_s \in M_s$. On the other hand, we
have $$\pi_{free}^{ij} = \pi_{free}^j \circ \pi_{free}^i$$ by the
transitivity of free products, see e.g. Proposition 2.5.5. in
\cite{VDN}. Therefore, if we let $\alpha_j = \pi_j \ten
\pi_{free}^j$ and amalgamate over $M_s$, we may rewrite the term
above as follows
$$\left( \sum_{i=1}^{\mathrm{k}_m} \sum_{j=1}^n
\sum_{w=1}^{\mathrm{k}_s} \pi_{tens}^w \Big[ \alpha_j \big( \xi_s
\ten x \ten \pi_{free}^{i} (a_m) \big) \Big] \ten \delta_{ijw}
\right)_s.$$ Then, arguing as in the proof of Lemma \ref{Exchange}
(in particular part i), we obtain $$\left( \sum_{j=1}^n
\widehat{\alpha}_j \Big[ \sum_{i=1}^{\mathrm{k}_m}
\sum_{w=1}^{\mathrm{k}_s} \pi_{tens}^w \big( \xi_s \ten x \ten
\pi_{free}^{i} (a_m) \big) \ten \delta_{iw} \Big] \ten \delta_j
\right)_s$$ with $\widehat{\alpha}_j$ a tensor amplification of
$\alpha_j$. However, according to Lemma \ref{Exchange} i) the
$\widehat{\alpha}_j$'s provide an increasingly independent family
of top-subsymmetric copies over the symmetric tensor product of
$M_s$. Hence, we are in position to apply Theorem \ref{4-term}
with $\pi_j$ replaced by $\widehat{\alpha}_j$ and $\mathcal{R} =
\ell_\infty^{\mathrm{k}_s \mathrm{k}_m}$. Again, the constants are
independent of the involved parameters. This gives us a new term
which does not depend on the choice of the morphisms
$\widehat{\alpha}_j$, so that we may use
$$\widehat{\alpha}_j = \pi_{free}^j \ten \pi_{free}^j$$ instead.
The assertion is then obtained by calculating backwards. \fin

\begin{remark} \emph{If we do not require the
constant to be $(p,q)$-independent, Theorem B also holds for $p>q$
by a simple duality argument. The singularity arises in this case
from the complementation constant of the subspace of independent
copies in $L_p(\mathcal{A}; \ell_q^n)$. As for $(\Sigma_{pq})$,
this singularity is not removable.}
\end{remark}

\begin{corollary} \label{CorAPThB}
Let $1 \le p \le q \le \infty$ and let $(\mathcal{M}_k)_{k \ge 1}$
be an increasingly independent family of top-subsymmetric copies
of $\mathcal{M}$ over $\mathcal{N}$. Then, we have an isomorphic
embedding $$x \in \mathcal{K}_{p,q}^n(\mathcal{M}, \mathcal{E_N})
\mapsto \sum_{k=1}^n \pi_k(x) \otimes \delta_k \in L_p \big(
\mathcal{A}; \ell_q^n \big)$$ with complemented range and
constants independent of $n$. In particular, replacing
$(\mathcal{M}, \mathcal{N}, \mathcal{E_N})$ by $(M_m(\mathcal{M}),
M_m, id_{M_m} \otimes \varphi)$ so that $\mathcal{A} =
M_m(\mathcal{R})$ for some $\mathcal{R}$, we obtain a complete
isomorphism with completely complemented range and constants
independent of $n$ $$x \in \mathcal{K}_{p,q}^n(\mathcal{M})
\mapsto \sum_{k=1}^n \pi_k(x) \otimes \delta_k \in L_p \big(
\mathcal{R}; \ell_q^n \big).$$
\end{corollary}

\dem It follows immediately from Theorem \ref{Theorem-Sigma-pq}
and Theorem B. \fin

\begin{remark} \label{Obsnota}
\emph{According to Remarks \ref{Singpq} and \ref{RemKfree2} we
know that, except for the case $(p,q) \sim (1,\infty)$, the
constants in Corollary \ref{CorAPThB} are also independent of
$p,q$. On the other hand, since we are using transference, we need
to work with independent copies. In the free case we can also work
with non i.d. variables, see Theorem \ref{Theorem-Sigma-pq}.}
\end{remark}

The rest of the paper is devoted to the proof of Corollary B. We
begin by stating a refinement of \cite[Theorem 4.2]{JP} which
follows easily from our previous results in this paper. We shall
write $L_p(M_n)$ for the Schatten class $S_p^n$ equipped with the
normalized trace $\frac1n \mathrm{tr}_n$.

\begin{lemma} \label{discrteembedding}
Let $1 \le p \le q \le 2$ and a positive integer $n \ge 1$. Then,
the following mapping is a complete isomorphism onto a completely
complemented subspace with constants independent of $p,q$ and $n$
$$\Psi_{pq}: x \in L_q(M_n) \mapsto \frac{1}{n^{2/q}}
\sum_{k=1}^{n^2} \pi_{tens}^k(x) \ten \delta_k \in
L_p(M_{n^{n^2}}; \ell_q^{n^2}).$$ As before, if $1 \le p \le q \le
\infty$, the same holds with a singularity when $(p,q) \sim
(1,\infty)$.
\end{lemma}

\dem According to Theorem 4.2 and Remark 4.3 in \cite{JP}, the
assertion holds for $1 < p \le q \le \infty$ with a constant $c_p$
majorized by $p/p-1$. The fact that it also holds for $p=1$ now
follows from Theorem \ref{4-term} and the argument in \cite{JP}.
The universality of the constants follows by Corollary
\ref{CorAPThB} + Remark \ref{Obsnota} followed by the original
argument \cite{JP} again. \fin

\begin{remark}
\emph{The choice $m=n^2$ in $L_p(M_{n^m}; \ell_q^m)$ is optimal,
see \cite{JP} for details.}
\end{remark}

In what follows, we will need some preparation on ultraproducts of
semifinite von Neumann algebras. Let $(\mathcal{M}_n)$ be a family
of semifinite von Neumann algebras with normal semifinite faithful
traces $(\tau_n)$. We may define $\tau_\U(x_n) = \limm_{n,\U}
\tau_n(x_n)$ on the ultraproduct von Neumann algebra
$$\mathcal{M}_\U = \big( \prodd_{n,\U} L_1(\mathcal{M}_n) \big)^*.$$
Let us set $$\mathcal{M}_{\U,sf} = \overline{\Big\{ (q_n x_n
q_n)^\bullet \, \big| \ \limm_{n,\U} \tau_n(q_n) < \infty, \,
(x_n)^\bullet \in \mathcal{M}_\U \Big\}}^{\mathrm{wot}}.$$ Then it
turns out that $\mathcal{M}_{\U,sf}$ is a semifinite von Neumann
subalgebra of $\mathcal{M}_\U$ and $\tau_\U((x_n)^{\bullet})) =
\lim_{n,\U} \tau_n(x_n)$ defines an trace on
$\mathcal{M}_{\U,sf}$. An appropriate way to check this consist in
checking the axioms of a (tracial) Hilbert algebra $$\mathcal{A} =
\Big\{ (x_n)^{\bullet} \, \big| \ \limm_{n,\U} \|x_n\|<\infty \
\mbox{and} \ \limm_{n,\U} \tau(x_n^*x_n) < \infty \Big\},$$ where
$(x_n)^{\bullet}$ corresponds to the equivalence class of a
bounded sequence of positive elements in $(\prod_{\U}
L_1(\mathcal{M}_n))^*$. Then $\tau_\U$ can be extended to a normal
semifinite trace on $\mathcal{M}_{\U,sf}$ viewed as the closure of
$\mathcal{A}$ in the GNS-representation of the Hilbert algebra. We
refer to \cite{Ray} for more on ultraproducts of noncommutative
$L_p$ spaces and how they can be identified as a noncommutative
$L_p$ space $L_p((\prod_{\U} L_1(\M_n))^*)$. Let $\mu_s(x)$ stand
for the generalized $s$-numbers of $x$, see \cite{FK}. The
following will be a key result below.

\begin{lemma} \label{crit}
Let $(x_n)$ be a bounded sequence in $L_p(\mathcal{M}_n)$ such
that $$\lim_{\delta \to 0} \, \limm_{n,\U} \int_0^{\delta}
\mu_s(x_n)^p \, ds = 0 = \lim_{\gamma \to \infty} \limm_{n,\U}
\int_\gamma^{\infty} \mu_s(x_n)^p \, ds. $$ Then we have
$(x_n)^{\bullet} \in L_p(\mathcal{M}_{\U,sf})$. Moreover, the
converse is also true.
\end{lemma}

\dem Given $\eps>0$, we choose $\gamma,\delta$ such that
$$\max \Big\{ \int_0^{\delta} \mu_s(x_n)^p
\, ds, \int_\gamma^{\infty} \mu_s(x_n)^p \, ds \Big\} < \eps/2.$$
It is clearly no restriction to assume that the $x_n$'s are
positive elements. Let us set $a_n = \mu_{\gamma}(x_n)$ and $b_n =
\mu_{\delta}(x_n)$. If $q_n = 1_{[a_n,b_n]}(x_n)$, we observe that
$\tau(q_n) \le \gamma$ and that $z_n = q_n x_n q_n$ is bounded by
$b_n$. Note  that $\delta^p \mu_{\delta}(x_n)\le \|x_n\|_p^p$
implies that $$\lim_{n,\U} b_n \le
\delta^{-p}\lim_{n,\U}\|x_n\|_p^p$$ is well-defined. Therefore
$(z_n)^\bullet \in \mathcal{M}_{\U,sf}$. The first assertion then
follows from $\|x_n-y_n\|_p^p < \eps$. For the converse we observe
that $\mathcal{M}_{\U,sf}$ is norm dense in $L_p(\mathcal{M}_\U)$
and the assertion is trivially true for $(x_n)^\bullet$ in
$\mathcal{M}_{\U,sf}$. The proof is complete. \fin

\noindent {\bf Proof of Corollary B.} If $1 \le p < q \le 2$, we
shall prove:
\begin{itemize}
\item[a)] There is no cb-embedding of $R_q+C_q$ into semifinite
$L_p$.

\item[b)] Let $\mathcal{R}_0$ stand for the hyperfinite
$\mathrm{II}_1$ factor and assume that there exists a complete
embedding of $\ell_q$ into $L_p(\mathcal{M})$ with $\mathcal{M}$
semifinite. Then, $L_q(\mathcal{R}_0)$ cb-embeds into some
semifinite $L_p$ space.
\end{itemize}
The combination of both results gives rise to the assertion.
Indeed, we know from the noncommutative Khintchine inequality
\cite{LuP} that $R_q+C_q$ cb-embeds into $L_q[0,1]$ which also
cb-embeds into $L_q(\mathcal{R}_0)$. Therefore, we deduce from a)
that there is no cb-embedding of $L_q(\mathcal{R}_0)$ into
semifinite $L_p$. Apply b) to conclude.

\vskip5pt

\noindent \textbf{Step 1.} The proof of a) essentially reproduces
Xu's argument in \cite{X3}. Assume there exists a complete
embedding $j: R_q+C_q \to L_p(\mathcal{M})$ with $\mathcal{M}$
semifinite and equipped with a normal semifinite faithful trace
$\tau$. Let $$j^*: L_{p'}(\mathcal{M}) \to R_q \cap C_q$$ denote
the adjoint mapping. Since $R_q \cap C_q$ can be regarded as the
diagonal subspace of $R_q \oplus C_q$, we may write $j^* =
(\Lambda_1, \Lambda_2)$. Since $\Lambda_1: L_{p'}(\mathcal{M}) \to
R_q$ is completely bounded, we deduce that there exists a positive
unit functional $f \in L_{p'/2}(\mathcal{M})^*$ such that
$$\|\Lambda_1(x)\|_{R_q} \le c f(xx^*)^{\frac{1-\theta}{2}}
f(x^*x)^{\frac{\theta}{2}} \quad \mbox{with} \quad 1/q =
(1-\theta)/p + \theta/p'.$$ This was proved by Pisier \cite{P4}
for $p=1$ and by Xu \cite{X2} for $1 < p \le 2$. We also refer to
\cite[Lemma 5.8]{X3} for a precise statement. When $p>1$, $f$ can
be regarded as a positive element in the unit ball of
$L_{p/(2-p)}(\mathcal{M})$ while for $p=1$, $f$ can be taken as a
normal state on $\mathcal{M}$ since $\Lambda_1$ is normal, see
\cite{P5} for details. In particular, we deduce
$$\|\Lambda_1(x)\|_{R_q} \le c \, \tau(fxx^*)^{\frac{1-\theta}{2}}
\tau(fx^*x)^{\frac{\theta}{2}}.$$ Arguing as in the proof of
Theorem 5.6 of \cite{X3}, we may apply an approximation argument
which allows us to assume that $\mathcal{M}$ is finite and $f =
\mathbf{1}_\mathcal{M}$. In that case our estimate for $\Lambda_1$
becomes $\|\Lambda_1(x)\|_{R_q} \le c \, \tau(xx^*)^{1/2}$.
Moreover, the same argument for $\Lambda_2$ produces
\begin{equation} \label{Grothen}
\|j^*(x)\|_{R_q \cap C_q} = \max \Big\{ \|\Lambda_1(x)\|_{R_q},
\|\Lambda_2(x)\|_{C_q} \Big\} \le c \, \tau(xx^*)^{\frac12}.
\end{equation}
This provides a factorization $j^* = v^* u^*$ with $u^*:
L_{p'}(\mathcal{M}) \to L_2(\mathcal{M})$ the natural inclusion
map. Arguing (twice) as in \cite{X3}, we see that $u^*$ becomes a
complete contraction when we impose on $L_2(\mathcal{M})$ the
o.s.s. of $L_2^{c_{p'}}(\mathcal{M}) \cap
L_2^{r_{p'}}(\mathcal{M})$. On the other hand, it follows from
\eqref{Grothen} that $v^*$ is a bounded map between Hilbert spaces
so that $v \in S_\infty$. To conclude, we note that $j=uv$
provides a factorization
$$R_q+C_q \stackrel{v}{\longrightarrow} R_p+C_p
\stackrel{u}{\longrightarrow} j(R_q+C_q).$$ By a simple
modification of \cite[Lemma 5.9]{X3}, we deduce that $$u \in
\mathcal{CB}(R_p+C_p,R_q+C_q) = S_{2pq/|p-q|}.$$ Thus, the
identity on $R_q+C_q$ belongs to $S_{2pq/|p-q|}$ which contradicts
$1 \le p < q \le 2$.

\vskip5pt

\noindent \textbf{Step 2.} Assume that there exists a cb-embedding
$j_p$ of $\ell_q$ into $L_p(\mathcal{M})$ for some semifinite von
Neumann algebra $\mathcal{M}$ equipped with a normal faithful
semifinite trace $\tau$. According to Lemma \ref{discrteembedding}
and our assumption, we find a cb-embedding
\begin{eqnarray*}
u_{np}: x \in L_q(M_n) & \mapsto & \frac{1}{n^{2/q}}
\sum_{k=1}^{n^2} \pi_{tens}^k(x) \ten \delta_k \in
L_p(M_{n}^{\ten_{n^2}};\ell_q^{n^2}) \\ & \mapsto &
\frac{1}{n^{2/q}} \sum_{k=1}^{n^2} \pi_{tens}^k(x) \ten
j_p(\delta_k) \in L_p(M_n^{\ten_{n^2}} \ten \mathcal{M}).
\end{eqnarray*}
Taking ultraproducts, we find a cb-embedding $$w_p:
L_q(\mathcal{R}_0) \to \prodd_{n,\U} L_p \big( M_{n}^{\ten_{n^2}}
\ten \mathcal{M} \big) = L_p(\widehat{\mathcal{M}}_{\U}).$$ On the
other hand, by our assumption we may regard $\ell_q$ as an
infinite-dimensional subspace of $L_p(\mathcal{M})$ not containing
$\ell_p$. According to the noncommutative form \cite{JP3} of
Rosenthal's theorem, given any $p < r < q$ we may find a positive
density $d \in L_1(\mathcal{M})$ with $\tau(d)=1$ and a embedding
$j_r: \ell_q \to L_r(\mathcal{M})$ satisfying
\begin{equation} \label{Rosenthalth}
j_p(x) = d^{\frac1p - \frac1r} j_r(x) + j_r(x) d^{\frac1p -
\frac1r}.
\end{equation}
Let us consider the map
\begin{eqnarray*}
u_{nr}: x \in L_q(M_n) & \stackrel{u_{r}^1}{\longmapsto} &
\frac{1}{n^{2/q}} \sum_{k=1}^{n^2} \pi_{tens}^k(x) \ten \delta_k
\in L_r(M_{n}^{\ten_{n^2}};\ell_q^{n^2}) \\ &
\stackrel{u_{r}^2}{\longmapsto} & \frac{1}{n^{2/q}}
\sum_{k=1}^{n^2} \pi_{tens}^k(x) \ten j_r(\delta_k) \in
L_r(M_n^{\ten_{n^2}} \ten \mathcal{M}).
\end{eqnarray*}
The first half is a complete embedding by Lemma
\ref{discrteembedding}. The second one is not necessarily bounded
since $j_r$ is not necessarily completely bounded. However, it is
easily seen that the composition of both is an isomorphic
embedding. Namely, we may clearly assume that $x \in L_q(M_n)$ is
self-adjoint. In that case, $$\mathcal{A}_x = \big\langle
\pi_{tens}^k (x) \big\rangle_{1 \le k \le n^2}''$$ is a
commutative von Neumann algebra. Thus
\begin{eqnarray*}
\lefteqn{\|u_{nr}(x)\|_{L_r(M_{n^{n^2}} \ten \mathcal{M})}} \\ & =
& \Big\| \frac{1}{n^{2/q}} \sum_{k=1}^{n^2} \pi_{tens}^k(x) \ten
j_r(\delta_k) \Big\|_{L_r(\mathcal{A}_x \ten \mathcal{M})} \\ &
\sim & \Big\| \frac{1}{n^{2/q}} \sum_{k=1}^{n^2} \pi_{tens}^k(x)
\ten \delta_k \Big\|_{L_r(\mathcal{A}_x; \ell_q^{n^2})} =
\|u_r^1(x)\|_{L_r(\mathcal{A}_x;\ell_q^{n^2})} \sim
\|x\|_{L_q(M_{n^{n^2}})}.
\end{eqnarray*}
We may take ultraproducts again and consider $$w_r:
L_q(\mathcal{R}_0) \to \prodd_{n,\U} L_r \big( M_n^{\ten_{n^2}}
\ten \mathcal{M} \big) = L_r(\widehat{\mathcal{M}}_\U).$$ If
$\delta = (\delta_n)^\bullet$ with $\delta_n =
\mathbf{1}_{M_{n^{n^2}}} \ten d$ and according to
\eqref{Rosenthalth}, we have
\begin{equation} \label{Rosth2}
w_p(x) = \delta^{\frac1p - \frac1r} w_r(x) + w_r(x)
\delta^{\frac1p - \frac1r}.
\end{equation}
We claim that $w_p(x)$ belongs to
$L_p(\widehat{\mathcal{M}}_{\U,sf})$, the semifinite part of
$\widehat{\mathcal{M}}_{\U}$, for any $x \in L_q(\mathcal{R}_0)$.
It suffices to check that the limits of Lemma \ref{crit} are zero.
We do it only for $\delta^{1/p - 1/r} w_r(x)$, since the term
$w_r(x) \delta^{1/p - 1/r}$ is estimated similarly. Since we have
$\mu_s(ab) \le \mu_{s/2}(a) \mu_{s/2}(b)$ and $\mu_s(\delta_n) =
\mu_s(d)$ for all $n$, we set $\frac1t = \frac1p - \frac1r$ and
obtain
\begin{eqnarray*}
\Big( \int_0^{2\delta} \mu_s \big( \delta_n^{\frac1p-\frac1r}
u_{nr}(x) \big)^p \, ds \Big)^{\frac1p} & \le & 2^{\frac1p} \Big(
\int_0^{\delta} \mu_s (d^{\frac1p-\frac1r})^t \, ds
\Big)^{\frac1t} \Big( \int_0^{\delta} \mu_s (u_{nr}(x))^r \, ds
\Big)^{\frac1r} \\ & \le & 2^{\frac1p} \Big( \int_0^{\delta} \mu_s
(d^{\frac1p-\frac1r})^t \, ds \Big)^{\frac1t} \|u_{nr}(x)\|_r \\ &
\lesssim & 2^{\frac1p} \Big( \int_0^{\delta} \mu_s
(d^{\frac1p-\frac1r})^t \, ds \Big)^{\frac1t} \|x\|_q.
\end{eqnarray*}
Therefore, we deduce that $$\lim_{\delta \to 0} \limm_{n, \U}
\Big( \int_0^{\delta} \mu_s \big( \delta_n^{\frac1p-\frac1r}
u_{nr}(x) \big)^p \, ds \Big)^{\frac1p} = 0.$$ The argument to
estimate $\int_\gamma^\infty$ is exactly the same. This completes
the proof. \fin

\begin{remark}{\rm  Following a suggestion by G.
Pisier let us describe what goes `wrong' when considering a family
$v_n: \ell_q^n \to L_p(\M_n)$ of complete embeddings into a
semifinite von Neumann algebras $\M_n$. Note that by local
reflexivity such cb-isomorphism exists. Let us consider the
following conditions.

\vskip3pt

\begin{itemize}
\item[i)] There exits a cb-embedding of $\ell_q$ into $L_p(\M)$
with $\M$ semifinite.

\vskip3pt

\item[ii)] There exists a sequence $(\M_n)_{n \ge 1}$ of
semifinite von Neumann algebras and linear maps $v_n: \ell_q^n \to
L_p(\M_n)$ such that $\|v_n\|_{cb} \|v_n^{-1}\|_{cb} \le c$ for
all $n \ge 1$ and for every $f \in L_q(0,1)$, the sequence
$(f_n)^\bullet$ determined by
$$f_n = v_n \Big( n^{1- \frac1q} \sum_{k=1}^n
\big( \int_{\frac{k-1}{n}}^{\frac{k}{n}} f(x) \, dx \big) \delta_k
\Big)$$ belongs to the semifinite part $L_p(\widehat{\M}_{\U,sf})$
of $\prod_{n,\U} L_p(\M_n,\tau_n)$.

\vskip3pt

\item[iii)] There exists a sequence $(\M_n)_{n \ge 1}$ of
semifinite von Neumann algebras and linear maps $w_n: \ell_q^n \to
L_p(\M_n)$ such that $\|w_n\|_{cb} \|w_n^{-1}\|_{cb} \le c$ for
all $n \ge 1$ and there exists a sequence of densities $d_n \in
L_1(\M_n)$ such that $(d_n)^\bullet$ belongs to the semifinite
part of $\prod_{n,\U} L_1(\M_n, \tau_n)$ and
$$w_n(x) = d_n^{\frac1p - \frac1r} j_{n,r}(x) + j_{n,r}(x) \,
d_n^{\frac1p - \frac1r},$$ is Rosenthal's factorization \cite{JP3}
for some contractions $j_{n,r}: \ell_q^n \to L_r(\M_n)$.

\vskip3pt

\item[iv)] There exits a cb-embedding of $L_q(\mathcal{R}_0)$ into
$L_p(\M)$ with $\M$ semifinite.
\end{itemize}

\vskip5pt

\noindent We will show that the conditions above are equivalent.
Hence, even though a family of complete embeddings $v_n: \ell_q^n
\to L_p(\M_n)$ with uniformly controlled constants exists, the
uniform integrability condition in ii) or iii) is violated.
}\end{remark}

\dem The implication iv) $\Rightarrow$ i) is obvious and we have
seen in the proof of Corollary B above that i) $\Rightarrow$ iii),
just take the same $d_n$ all the time. The proof of ii)
$\Rightarrow$ iii) follows similarly. Indeed, by assumption we
obtain a continuous map $v: L_q(0,1) \to
L_p(\widehat{\M}_{\U,sf})$. We apply the noncommtuative Rosenthal
theorem \cite{JP3} and find $v(f) = d^{\frac1p - \frac1r} j(f) +
j(f) \, d^{\frac1p - \frac1r}$ for a bounded map $j: L_q(0,1) \to
L_{r}(\widehat{\M}_{\U,sf})$ and some density $d \in
L_1(\widehat{\M}_{\U,sf})$. By restricting $j$ to step functions
on the intervals $[\frac{k-1}{n}, \frac{k}{n}]$, we have found the
complete embeddings $w_n$ from condition iii). For the implication
iii) $\Rightarrow$ iv) we apply the argument from our proof of
Corollary B. Indeed, it suffices to check that for every
self-adjoint $x \in \mathcal{R}_0$, the sequence
\begin{eqnarray*}
u_n(x) & = & \frac{1}{n^{2/q}} \sum_{k=1}^{n^2} \pi_k(x) \ten
w_{n^2}(\delta_k) \\ & = & \frac{1}{n^{2/q}} \sum_{k=1}^{n^2}
\pi_k(x) \ten \big( d_{n^2}^{\frac1p-\frac1r} j_{n^2,r}(\delta_k)
+ j_{n^2,r}(\delta_k) \, d_{n^2}^{\frac1p - \frac1r} \big)
\end{eqnarray*}
belongs to the semifinite part of $\prod_{n,\U}
L_p(\mathcal{R}_0^{\ten_{n^2}}\ten \M_n)$. Referring to the
argument after \eqref{Rosth2}, it suffices to note that
\begin{eqnarray*}
\lefteqn{\hskip-10pt \Big\| \frac{1}{n^{2/q}} \sum_{k=1}^{n^2}
\pi_k(x) \ten j_{n^2,r}(\delta_k)
\Big\|_{L_r(\mathcal{R}_0^{\ten_{n^2}} \ten \M_n)}} \\ & \le &
\|j_{n^2,r}\| \, \frac{1}{n^{2/q}} \, \Big( \int_{\mathcal{A}_x}
\Big[ \sum_{k=1}^{n^2} |\pi_k(x)|^q \Big]^{\frac{r}{q}} d\mu_x
\Big)^{\frac1r} \le \|j_{n^2,r}\| \|x\|_{\mathcal{R}_0}
\end{eqnarray*}
is uniformly bounded in $n$ for every $x\in \mathcal{R}_0$. The
proof is complete. \fin

\vskip5pt

\noindent \textbf{Acknowledgement.} We would like to thank the
referee for a careful reading of this paper. His/her comments on
an earlier version of this paper have led us to improve the
presentation at several points.

\bibliographystyle{amsplain}


\vskip15pt


\noindent \textbf{Marius Junge} \\
\textsc{Department of Mathematics} \\
\textsc{University of Illinois at Urbana-Champaign} \\ 273 Altgeld
Hall, 1409 W. Green Street. Urbana, IL 61801. USA \\
\texttt{junge@math.uiuc.edu}

\vskip15pt

\noindent \textbf{Javier Parcet} \\
\textsc{Instituto de Ciencias Matem\'aticas} \textsc{CSIC-UAM-UC3M-UCM} \\
\textsc{Consejo Superior de Investigaciones Cient\'ificas} \\
C/ Serrano 121. 28006, Madrid. Spain \\
\texttt{javier.parcet@uam.es}

\end{document}